\input amstex

\loadeufm
\loadmsbm
\loadeufm

\documentstyle{amsppt}
\input amstex
\catcode `\@=11
\def\logo@{}
\catcode `\@=12
\magnification \magstep1
\NoRunningHeads
\NoBlackBoxes
\TagsOnLeft

\def \={\ = \ }
\def \+{\ +\ }
\def \-{\ - \ }

\def \b|{\big |}

\def \g1{\Gamma_1}

\def \nfp{\demo\nofrills{Proof:\usualspace\usualspace }}

\def\rarr#1#2{\smash{\mathop{\hbox to .5in{\rightarrowfill}}
 	 \limits^{\scriptstyle#1}_{\scriptstyle#2}}}

\def\larr#1#2{\smash{\mathop{\hbox to .5in{\leftarrowfill}}
	  \limits^{\scriptstyle#1}_{\scriptstyle#2}}}

\def\swarr#1#2 {\llap{$\scriptstyle #1$}  \swarrow
  	\vcenter to .5in{}\rlap{$\scriptstyle #2$}}

\topmatter
\title Non-Commutative Algebraic Geometry I: 
\centerline{Monomial Equations with a Single Variable}
\endtitle
\author
\centerline{ 
Z. Sela${}^{1,2}$}
\endauthor
\footnote""{${}^1$Hebrew University, Jerusalem 91904, Israel.}
\footnote""{${}^2$Partially supported by an Israel academy of sciences fellowship.} 
\abstract\nofrills{}
This paper is the first in a sequence 
on the structure of
sets of solutions to systems of equations over a free associative algebra. 
We start by constructing a Makanin-Razborov diagram that encodes all the homogeneous solutions to a
homogeneous monomial system of equations. Then we analyze the set of solutions to monomial systems of equations with a single variable.
\endabstract
\endtopmatter

\document

\baselineskip 12pt

Algebraic geometry studies the structure of sets of solutions to systems of equations usually
over fields or commutative rings.
The developments and the considerable abstraction that currently exist in the study of varieties over commutative rings,
still resist to apply to the study of varieties over non-abelian rings or over other non-abelian algebraic structures.

Since about 1960 ring theorists, P. M. Cohn [Co], G. M. Bergman [Be-Ha], and others, have tried to study varieties over
non-abelian rings, notably free associative algebras (and other free rings). However, the pathologies that they tackled and
the lack of unique factorization that they study in detail ([Co], chapters 3-4), prevented any attempt to prove or even speculate
what can be the structure of varieties over free associative algebras. 

In this sequence of papers we suggest to study  varieties over free associative algebras using techniques and analogies of structural results from
the study of varieties over free groups and semigroups. Over free groups and semigroups geometric techniques as well as low dimensional topology
play an essential role in the structure of varieties. These include Makanin's algorithm for solving equations, Razborov's analysis of sets of solutions over
a free group, the concepts and techniques that were used to construct and analyze the JSJ decomposition, and the applicability of the JSJ machinery to
study varieties over free groups and semigroups ([Se1],[Se2]). Our main goal is
 to demonstrate  that these techniques and concepts can be modified to be applicable over free associative algebras as well.

Furthermore, we believe that the concepts and techniques that proved to be successful over free groups and semigroups, can be adapted to  
analyze varieties over free objects in other non-commutative and at least "partially" associative algebraic structures. In that respect, we hope that it will be possible
to use or even axiomatize the properties of varieties over the free objects in these algebraic structures, 
in order to set dividing lines between non-commutative algebraic structures, in analogy with classification theory (of first order theories) in model theory [Sh].

We start the analysis of systems of equations over a free associative algebra with what we call monomial systems of equations. These are systems of equations 
over a free associative algebra, in which every polynomial in the system contains  two monomials. 
In the first section of this paper we analyze the case of homogeneous solutions to homogeneous monomial systems of equations. In this case it is possible
to apply the techniques that were used in analyzing varieties over free semigroups [Se2], and associate a Makanin-Razborov diagram that encodes all the 
homogeneous solutions to homogeneous monomial system of equations.

In the second section we introduce $limit$ $algebras$ that are a natural analogue of a $limit$ $group$, and prove that such algebras are always embedded in 
(limit) division algebras.  The automorphism (modular) groups of these division algebras are what is needed in the sequel in order to modify and shorten 
solutions to monomial systems of equations.

In the third chapter we present a combinatorial approach to (cases of) the celebrated Bergman's centralizer theorem [Be]. Finally, in the fourth section we use this 
combinatorial approach to analyze the set of solutions to a monomial system of equations with a single variable. The results that we obtain are analogous
to the well known structure of the set of solutions to systems of equations with a single variable over a free group or semigroup. We prove all our results under the assumption 
that the top homogeneous parts of the coefficients in the equations are monomials with no periodicity, in order to simplify our arguments, but we believe 
that eventually this assumption can be dropped.

In the next paper in the sequence we use the techniques that are presented in this paper to analyze monomial systems of equations that have no quadratic (or surface)
parts. In the third paper in the sequence we analyze the quadratic parts of monomial systems of equations. Eventually, we hope to use our analysis of sets of solutions
to monomial system of equations to the analysis of general varieties.

\smallskip
I would like to thank Eliyahu Rips who inspired me years ago to try and study  varieties over non-commutative algebras, and David Kazhdan for his continuing interest.

\vglue 1.5pc
\centerline{\bf{\S1. Homogeneous solutions of monomial equations }}
\medskip
For simplicity,  we will always assume that the free algebras that we consider are over the field with two
elements $GF_2$. 
Let $FA$ be a free associative algebra over $GF_2$: $FA=GF_2<a_1,\ldots,a_k>$. In order to study the structure of general
varieties over the associative algebra $FA$, we start with varieties that are defined by monomial systems of equations.
A system of equations, $\Phi$, is called $monomial$ if it is defined using a finite set of
unknowns:
$x_1,\ldots,x_n$, and a finite set of   equations:
$$\align
u_1(c_1,\ldots,c_{\ell},  x_1,\ldots,x_n) &= 
v_1(c_1,\ldots,c_{\ell},  x_1,\ldots,x_n)  \\
& \vdots \hfill \\
u_s(c_1,\ldots,c_{\ell},  x_1,\ldots,x_n) &= 
v_s(c_1,\ldots,c_{\ell},  x_1,\ldots,x_n)  \\
\endalign$$
where the words $u_i$ and $v_i$ are monomials in the free algebra generated by  the variables: 
$x_1,\ldots,x_n$ and coefficients from the algebra $FA$: $c_1,\ldots,c_{\ell}$, i.e., a word in the free semigroup generated by these elements (note that
the coefficients, $c_1,\ldots,c_{\ell}$, are general elements and not necessarily monomials). 
A monomial system of equations is called $homogeneous$ if all the coefficients, $c_1,\ldots,c_{\ell}$,  in the system are homogeneous elements in the free
associative algebra $FA$.

We start by analyzing all the homogeneous solutions of a homogeneous monomial system, i.e., all the assignments
of homogeneous elements in $FA$ to the variables $x_1,\ldots,x_n$, such that the equalities in
a homogeneous monomial  system of equations are valid.

Let $x_1^0,\ldots,x_n^0$ be a homogeneous solution of the monomial system $\Phi$. Substituting the elements 
$x_1^0,\ldots,x_n^0$ in the monomials $u_i$ and $v_i$, $1 \leq i \leq s$, we get a finite set of
equalities in the free algebra $FA$. Since all the elements that appear in each of these equalities are
homogeneous, for each index $i$ we can associate a segment, $J_i$, of length that is equal to the degree of
$u_i$ and $v_i$ after the substitution of $x^0_1,\ldots,x^0_n$. We further add notation on  the segment $J_i$
for the beginnings and the ends of each of the elements $x^0_1,\ldots,x^0_n$ and the coefficients of the system: $c_1,\ldots,c_{\ell}$.

With the segments $J_1,\ldots,J_s$, and the notation for the beginnings and ends of $x^0_1,\ldots,x^0_n$ and 
$c_1,\ldots,c_{\ell}$, we can naturally associate a band complex (bands are added for different appearances
of the same variable). All the lengths that appear in the complex are integers, so the band complex must be 
simplicial. Note that all the operations that are used in the Rips machine, or in the Makanin procedure, to transfer
the original complex into a standard simplicial complex are valid in our context, i.e., it is possible to cut the 
elements $x^0_1,\ldots,x^0_n$ and $c_1,\ldots,c_{\ell}$ and represent them as multiplication of
new elements, according to the operations that are performed in modifying the band complexes along the procedure.

The simplicial band complex (that is equivalent to the band complex we started with), divides the segments
$J_1,\ldots,J_s$, and the segments  of the elements, $x^0_1,\ldots,x^0_n$ and $c_1,\ldots,c_{\ell}$, into 
finitely many oriented subsegments
that belong to finitely many equivalence classes. With each equivalence class of subsegments we associate a label.
Note that the labels are divided into labels that get (fixed) values in the free algebra $FA$, and are part of the
coefficients, $c_1,\ldots,c_{\ell}$, and "free" labels.

Given  the labels of the oriented subsegments, with each of the variables $x_1,\ldots,x_n$ and the coefficients
$c_1,\ldots,c_{\ell}$ we naturally associate a word in the free semigroup that is generated by the labels.
By construction for each index $i$, the words that are associated with $u_i$ and $v_i$ in this free semigroup are
identical. Therefore, we have constructed a homomorphism from the semigroup that is formally generated by
elements that are associated with $x_1,\ldots,x_n$ and $c_1,\ldots,c_{\ell}$ to the free semigroup that is generated by
the labels, and this homomorphism (of semigroups) satisfies the equalities in the system $\Phi$, when we interpret
the monomials in these equalities as equations in a semigroup.   

Let $G_{\Phi}$ be the semigroup that is generated by copies of $x_1,\ldots,x_n$ and the coefficients: $c_1,\ldots,c_{\ell}$, modulo
the relations: 
$$u_i(x_1,\ldots,x_n,c_1,\ldots,c_{\ell})=v_i(x_1,\ldots,x_n,c_1,\ldots,c_{\ell})$$ 
for $1 \leq i \leq s$, where the monomials 
$u_i$ and $v_i$ are interpreted as words in a free semigroup. By the
previous argument with  a (homogeneous)
solution of the system $\Phi$ it is possible to associate a homomorphism from $G_{\Phi}$ into a free semigroup.

Conversely, given a semigroup homomorphism of $G_{\Phi}$ into a free semigroup, that fits with a decomposition of the 
constants $c_1,\ldots,c_{\ell}$ into a product of homogeneous elements (there are finitely many possible ways
to represent each of the coefficients $c_1,\ldots,c_{\ell}$ as such a product), it is possible to associate with 
such a product a family of solutions of the systems $\Phi$. 

Therefore, the study of homogeneous solutions of a homogeneous monomial system of equations  over an associative free algebra is reduced
to the study of a collection of semigroup homomorphisms from a given f.g.\ semigroup into a free semigroup. By $[Se2]$
with the collection of semigroup homomorphisms it is possible to associate canonically a finite collection
of pairs: $(S_1,L_1),\ldots,(S_m,L_m)$, where each of the groups $L_j$'s  is a limit group, and each of
the semigroups $S_j$ is a f.g.\ subsemigroup that generates $L_j$. Furthermore, with $G_{\Phi}$ and its collection
of homomorphisms it is possible to associate (non-canonically) a Makanin-Razborov diagram that encodes
all its homomorphisms into free semigroups. By our observation, this Makanin-Razborov
diagram of pairs encodes all the homogeneous solutions of the system $\Phi$. Indeed, every solution of the system $\Phi$ is
obtained from a homomorphism of $G_{\Phi}$ into a free semigroup by substituting arbitrary homogeneous elements from
$FA$ instead of the terminal (non-constant) free elements that appear in the terminal level of a Makanin-Razborov diagram of
$G_{\Phi}$.

\vglue 1.5pc
\proclaim{Theorem 1.1} With a homogeneous monomial system of equations over the free associative algebra $FA$ it
is possible to associate (non-canonically) a Makanin-Razborov diagram that encodes all its homogeneous solutions.
\endproclaim

As a corollary of the encoding of homogeneous solutions of system of homogeneous monomial equations by pairs of limit groups and
their subsemigroups we get the following.

\vglue 1.5pc
\proclaim{Corollary 1.2} The collections of sets of homogeneous solutions to homogeneous monomial systems of 
equations is Noetherian. i.e., every descending sequence of such sets terminates after a finite 
time. 
\endproclaim

\nfp Follows immediately from the descending chain condition for limit groups [Se1], or the Noetherianity of 
varieties over free groups
and semigroups [Gu].

\line{\hss$\qed$}

Theorem 1.1 associates a Makanin-Razborov diagram with the set of homogeneous solutions to a homogeneous monomial systems of equations. Our main
goal in this sequence of papers is to associate a Makanin-Razborov diagram with the set of (not necessarily homogeneous) solutions of a general
monomial system of equations, at least in the minimal rank case. i.e., in the case in which the Makanin-Razborov diagram that is
associated with the homogeneous system that is associated with top homogeneous part of the non-homogeneous system contains no free products.

\vglue 1.5pc
\centerline{\bf{\S2. Limit algebras, their division algebras and modular groups}}
\medskip
The construction of the Makanin-Razborov diagram of a system of equations over a free group uses extensively the (modular)
automorphism groups of the limit groups that are associated with its nodes. These modular groups enable one to proceed from a limit group to 
maximal shortening quotients of it, that are always proper quotients.

The semigroups that appear in the construction of the Makanin-Razborov diagram of a system of equations over a free semigroup do not have a large automorphism group in general.
 e.g., a finitely generated free semigroup have a finite automorphism group. Hence, to study homomorphisms from a given f.g.\ semigroup $S$ to the free semigroup $FS_k$ we did
the following in [Se2]. 

Given a f.g.\ semigroup, $S$, we can naturally associate a group with it. Given a presentation
of $S$ as a semigroup, we set the f.g.\ group $Gr(S)$ to be the group with the presentation of $S$
interpreted as a presentation of a group. Clearly, the semigroup $S$ is naturally mapped into the
group, $Gr(S)$, and the image of $S$ in $Gr(S)$ generates $Gr(S)$. We set  $\eta_S: S \to Gr(S)$ to
be this natural homomorphism  of semigroups. 

\noindent
The free semigroup, $FS_k$, naturally embeds into a free group, $F_k$. By the construction of the group, $Gr(S)$,
 every homomorphism of semigroups,
$h:S \to FS_k$, extends to a  unique homomorphism of groups, $h_G: Gr(S) \to F_k$, so that: $h=h_G \circ \eta_S$.

By construction, every
homomorphism (of semigroups), $h:S \to FS_k$, extends to a homomorphism (of groups), $h_G:Gr(S) \to F_k$. Therefore,
the study of the structure of $Hom(S,FS_k)$, is equivalent to the study of the structure of the collection of
homomorphisms of groups, $Hom(Gr(S),F_k)$, that restrict to homomorphisms of (the semigroup) $S$ into the
free semigroup (the $positive$ $cone$), $FS_k$.

By (canonically) associating a finite collection
of maximal limit quotients with the set of homomorphisms, $Hom(Gr(S),F_k)$,  that restrict to 
(semigroup) homomorphisms from $S$ to $FS_k$, we are able to (canonically) replace the pair $(S,Gr(S))$ with a finite collection of limit quotients:
$(S_1,L_1),\ldots,(S_m,L_m)$, where each of the groups $L_i$ is a limit group. Limit groups have rich modular groups, and these are later used to proceed to the
next levels of the Makanin-Razborov diagram of the given system of equations over the free semigroup $FS_k$.
   
\smallskip
In studying sets of solutions to systems of equations over a free associative algebras, we need to study homomorphisms: $h:A \to FA_k$, where $A$ is a f.p.\ algebra, and
$FA_k$ is the free associative algebra of rank $k$. As in the case of groups and semigroups, to study such homomorphisms we pass to convergent sequences of homomorphisms:
$\{h_n:A \to FA_k\}$, and look at the $limit$ $algebras$ $LA$ that are associated with such  convergent sequences. 
Algebras, and in particular limit algebras, have automorphisms,
but these are not the automorphisms that will be needed in the sequel to modify and shorten homomorphisms.

By a classical construction of Malcev [Mal] and Neumann [Ne], and by different constructions of Amitsur [Am] and others, the free associative algebra $FA_k$ can be embedded
into a division algebra $Div(FA_k)$ (note that there are various different division algebras into which $FA_k$ embeds). Given a convergent sequence $\{h_n:A \to FA_k\}$
with an associated limit algebra $LA$, it is  straightforward to get an embedding: $LA \to Div(LA)$, where $Div(LA)$ is a division algebra, that is also obtained
from the convergent sequence and from the embedding $FA_k \to Div(FA_k)$. 

In the sequel we will use (a subgroup of) the group of automorphisms of the division algebra $Div(LA)$, in order to modify (shorten) the homomorphisms $h: A \to FA_k$ that
we need to study. These will be the modular groups that are associated with limit algebras that appear along the nodes of the 
Makanin-Razborov diagrams of the given systems of equations
over the free associative algebra $FA_k$.

An important example is (a special case of) what we call $surface$ (or $quadratic$)  $algebra$: $$SA=<x_1,\ldots,x_n \, | \, x_1 \ldots x_n=x_{\sigma(1)} \ldots x_{\sigma(n)}>$$ 
for an appropriate  permutation
$\sigma \in S_n$. Such a surface algebra is a limit algebra. Hence, it is embedded in a division algebra, $Div(SA)$. For appropriate convergent sequences,
the modular group of $Div(SA)$ contains the automorphism group of a corresponding surface. Therefore, we  call the modular group of $Div(SA)$, 
the $Bergman$ $modular$ $group$
of a surface algebra, since it contains (or is generated by) generalized Dehn twists that are inspired by Bergman's centralizer theorem [Be]. 
These modular groups play an essential role in constructing Makanin-Razborov diagrams for monomial systems of equations over a free associative algebra in the sequel.

\vglue 1.5pc
\centerline{\bf{\S3. A combinatorial approach to Bergman's theorem}}
\medskip
In the first section we studied homogeneous solutions to homogeneous monomial systems of equations.
In this section we start the study of non-homogeneous solutions to arbitrary monomial systems of equations.
We start by studying the centralizers of elements in a free associative algebra, i.e., 
we give combinatorial proof to Bergman's theorem, and then use these techniques to study related system of
equations. We start with the following theorem, that can be proved easily by a direct induction, but we
present a proof that uses techniques that we will use in the sequel.

\vglue 1.5pc
\proclaim{Theorem 3.1} Let $u \in  FA$ be an element for which its top degree homogeneous part is a
monomial, and suppose that this top degree monomial has no non-trivial roots. Then the centralizer of
$u$ in $FA$ is precisely the elements in the (one variable) algebra that is generated by $u$.
\endproclaim

\nfp 
Suppose that $x$ is a (non-trivial) element that satisfies $xu=ux$. By our analysis of homogenous elements, the
top degree homogenous part of $x$ must be a monomial, which is a power of the top degree monomial in $u$.
Hence, the top dgree monomial of $x$ has to be identical to the top degree monomial of $u^m$ for some $m$.
Therefore, $deg(x+u^m) < deg(x)$ and $u(x+u^m)=(x+u^m)u$, so the theorem follows by induction on the degree of $x$.

For later applications we present a different proof. First, note that 
$xu=ux$ if and only if $x(u+1)=(u+1)x$. Hence, we may assume that the monomials in $u$ do not include the one corresponding
to the identity.

\vglue 1.5pc
\proclaim{Lemma 3.2} Suppose that $deg(u) \geq 2$ and that $deg(x) \geq 2deg(u)$. There exists an element $w \in FA$, such that: $x=uw=wu$ in the abelian group:
$G^{deg(x)}/G^{deg(x)-2}$, and $ux=xu=uwu$ in the abelian group: 
$G^{deg(xu)}/G^{deg(xu)-2}$, where $G^m$ is the abelian group that is generated by the monomials of degree $m$ in $FA$.  
\endproclaim

\nfp We analyze the codegree 1 monomials in the two sides of the equation: $xu=ux$. If there are no cancellations between the codegree 1 monomials
(that are obtained using the distributive law) in each side of the equation, the claim follows in the way that it follows in the homogeneous case.  

Suppose that there are cancellations between codegree 1 monomials in ithe product $ux$. In that case a codegree 1 monomial of $x$ can be written as
$w_1u_0$, where $u_0$ is the top monomial of $u$, and $deg(w_1)=deg(x)-deg(u)-1$. In this case we add the monomial $w_1$ as a codegree 1 monomial to $w$.
iNote that in this case $u_0w_1$ and $u_1w_0$ ia a pair of cancelling codegree 1 monomials,
 where $u_1$ is a codegree 1 monomial of $u$, that satisfies: $u_1x_0=u_0w_1u_0$.

Continuing the analysis for all the codegree 1 products in the two sides of the equation, we get the statement of the theorem.

\line{\hss$\qed$}


Since $x=t_1u$ in the abelian group 
$G^{deg(x)}/G^{deg(x)-2}$, if $deg(x)\geq 2deg(u)$ 
 it follows that $xu=t_1u^2$ and $ux=ut_1u$ in the abelain group
$G^{deg(xu)}/G^{deg(xu)-2}$.
 Hence, 
if $deg(t_1) \geq 2 deg(u)$, $t_1=t_2u=ut_2$ in 
$G^{deg(t_1)}/G^{deg(t_1)-2}$. Applying these arguments iteratively we get that $x=tu^m$ in 
the abelian group: $G^{deg(x)}/G^{deg(x)-2}$, for some $t$ that satisfies: $deg(t)=deg(u)$. 

Therefore, $xu=tu^mu=ux=utu^m$ in 
$G^{deg(xu)}/G^{deg(xu)-2}$, which means that: $tu=ut$ in the abelian group
$G^{deg(tu)}/G^{deg(tu)-2}$,
 and $deg(t)=deg(u)$.

Let $t$ satisfy $tu=ut$ in 
$G^{deg(tu)}/G^{deg(tu)-2}$,
and $deg(t)=deg(u)$. In this case the top degree monomials of $t$ and $u$ are identical,
and we denote this monomial $u_0$.
Suppose that $s,v$ are monomials of codegree 1 in either $t$ or $u$, 
and suppose that $su_0=u_0v$. In that case $v$ is the suffix of $u_0$ and $s$ is the prefix of $u_0$. Since
$u_0$ is not a proper power, $vu_0$ can not be presented as $u_0w$ for any codegree 1 monomial $w$,
and $u_0s$ can not be presented as $wu_0$ for any codegree 1 monomial $w$.

Hence, if $s,v$ are codegree 1 monomials in either $t$ or $u$, and $su_0=u_0v$, then both $u_0s$ and 
$vu_0$ can be presented uniquely in each of the two products $tu$ and $ut$, that imply that 
$s$ and $v$ must be codegree 1 monomials in both $u$ and $t$. Therefore, the codegree 1 monomials of
$t$ and $u$ must be identical, so $t=u$ in the abelian group:
$G^{deg(tu)}/G^{deg(tu)-2}$, and $x=u^{m+1}$ in the abelian group 
$G^{deg(xu)}/G^{deg(xu)-2}$,  
for some non-negative integer $m$.

We use (a finite) induction and assume that $x=u^{m+1}$ in the abelian group:
$G^{deg(x)}/G^{deg(x)-c}$, for some positive integer $c < deg(u)$. i.e., we assume that the equality holds for all 
the monomials in $x$ and $u$ of codegree smaller than $c$. To complete the proof of the theorem, we need to
prove the same equality for all the monomials of codegree at most $c$.

By the inductive hypothesis $x=u^{m+1}$ in the abelian group
$G^{deg(x)}/G^{deg(x)-c}$. Hence, $x=x_{c-1}+v$ where $x_{c-1}$ is the sum of all the monomials of codegree smaller than $c$
in $x$ and $deg(v) \leq deg(x)-c$. Furthermore, $x_{c-1}$ is precisely the sum of all the monomials of codegree smaller than $c$ in $u^{m+1}$.

Let $u_{c-1}$ be the sum of the monomials of codegree less than $c$ in $u$. We set $s_c$ to be the sum of all the monomials of 
codegree $c$ in $u_{c-1}^{m+1}$. By construction: $u_{c-1} (x_{c-1}+s_c)= (x_{c-1}+s_c) u_{c-1}=u_{c-1}^{m+2}$ in the abelian
group: 
$G^{deg(xu)}/G^{deg(xu)-(c+1)}$. i.e., the monomials of codegree at most $c$ are identical for the 3 different products.

Recall that $x=x_{c-1}+v$, where $deg(v) \leq deg(x)-c$. We set $x=x_{c-1} +s_c + r$, where $deg(s_c)=deg(x)-c$ and $deg(r) \ leq deg(x)-c$.
Let $q_c$ be the sum of the monomials of codegree $c$ in $u$. Then:
$$ux=(u_{c-1}+q_c)(x_{c-1}+s_c+r)=xu=(x_{c-1}+s_c+r)(u_{c-1}+q_c)$$ in the abelian group: 
$G^{deg(xu)}/G^{deg(xu)-(c+1)}$. Since: $u_{c-1}(x_{c-1}+s_c)=(x_{c-1}+s_c)u_{c-1}$ in that abelian group, it follows
that: $u_{top}r+q_cx_{top}=ru_{top}+x_{top}q_c$ in that abelian group, where $u_{top}$ and $x_{top}$ are the top monomials
in $u$ and $x$ in correspondence. Therefore, all these monomials are products of a top degree monomial with a codegree $c$
monomial, and these can be broken precisely as in the codegree 1 case, assuming $deg(x) \geq 2deg(u)$.

We are left with the case in which $deg(x)=deg(u)$. In that case we write: $u=u_{c-1}+q_c$ and $x=u_{c-1}+r_c$ in
the abelian group: $G^{deg(x)}/G^{deg(x)-(c+1)}$, where $q_c$ and $r_c$ are the codegree $c$ monomials in $u$ and $x$ in correspondence.
Since the contributions of products of monomials of codegree smaller than $c$ in $xu$ and in $ux$ are identical, we need to look only
at the equation: $u_{top}r_c+q_cx_{top}=x_{top}q_c+r_cu_{top}$ for the monomials of codegree $c$, where $x_{top}=u_{top}$ are identical monomials. 
By the argument that was used in the codegree 1 case (when $deg(x)=deg(u)$), it follows that $q_c=r_c$, and the general step of the induction is proved.

So far we may conclude that $x=u^{m+2}$ in the abelian group: 
$G^{deg(x)}/G^{deg(x)-(deg(u)+1)}$. Therefore, $x+u^{m+2}$ commutes with $u$ and $deg(x+u^{m+2}) \leq deg(x) - deg(u)$, and the theorem follows.

\line{\hss$\qed$}

So far we assumed that the top homogeneous element of $u$ is a monomial and that its top monomial  doesn't have a proper root. We continue by allowing $u$ to be a proper
power.

\vglue 1.5pc
\proclaim{Theorem 3.3} Let $u \in  FA$ be an element for which its top degree homogenous part is a
monomial, and suppose that $u=p(v)$ and the top degree monomial of $v$ does not have
a proper root. Then the centralizer of
$u$ in $FA$ is precisely the elements in the algebra that is generated by $v$.
\endproclaim

\nfp 
Suppose that $x$ is a (non-trivial) element that satisfies $xp(v)=p(v)x$. First, note that like theorem 3.1, theorem 3.3 can be proved easily by
replacing $x$ by $x+v^m$, for an appropriate $m$,  such that $deg(x+v^m)<deg(x)$, and $(x+v^m)p(v)=p(v)(x+v^m)$.
However, as in the proof of theorem 3.1 and for future purposes, we prefer to present a different proof. For that proof we assume that $deg(v)>1$.

As in theorem 3.1, by our analysis of homogenous elements, the
top degree homogenous part of $x$ must be a monomial, which is a power of the top degree monomial in $v$.

As in the proof of theorem 3.1, if $deg(x) > deg(u)$ the arguments that were used in
the proof of lemma 3.2, that remain valid under the assumptions of the theorem,
 enable us to analyze the codegree 1 monomials in $x$. In that case, as in the proof of theorem 3.1, 
there exists an element $t_1$,
that contains a top degree monomial and a homogenous part of codegree 1, such that $x=ut_1$ and $x=t_1u$,
i.e., that $xu=ux=ut_1u$,
in the ablian group: 
$G^{deg(x)}/G^{deg(x)-2}$. 

Applying these arguments iteratively, as in the proof of theorem 3.1,  we get that $x=tu^m$ in 
the abelian group: $G^{deg(x)}/G^{deg(x)-2}$, for some $t$ that satisfies: $deg(t) \leq deg(u)$, 
which means that: $tu=ut$ in the abelian group:
$G^{deg(tu)}/G^{deg(tu)-2}$.
In particular, the top degree monomial of $t$ must be a power
of the top degree monomial of $v$. 

In case $deg(t) \leq deg(u)$ and $tu=ut$ in the abelian group:
$G^{deg(tu)}/G^{deg(tu)-2}$,
we apply the same argument that we used in case $deg(t)=deg(u)$ in the 
proof of theorem 3.1. By these arguments, if $u=v^b$ in the ableian group:
$G^{deg(u)}/G^{deg(u)-2}$, then $t=v^s$ in the abelian group:
$G^{deg(t)}/G^{deg(t)-2}$, where $s$ is an integer, $1 \leq s \leq b$. This implies that $x=v^{\ell}$,
for some positive integer $\ell$, in the abelian group:
$G^{deg(x)}/G^{deg(x)-2}$.


We further use (a finite) induction and assume that $x=v^{\ell}$ in the abelian group:
$G^{deg(x)}/G^{deg(x)-c}$, for some positive integer $c < deg(v)$. i.e., we assume that the equality holds for all 
the monomials in $x$ and $v$ of codegree smaller than $c$. To complete the proof of the theorem, we need to
prove the same equality for all the monomials of codegree at most $c$.

By the inductive hypothesis $x=v^{\ell}$ in the abelian group
$G^{deg(x)}/G^{deg(x)-c}$. Hence, $x=x_{c-1}+h$ where $x_{c-1}$ is the sum of all the monomials of codegree smaller than $c$
in $x$ and $deg(h) \leq deg(x)-c$. Furthermore, $x_{c-1}$ is precisely the sum of all the monomials of codegree smaller than $c$ in $v^{\ell}$.

Let $u_{c-1}$ be the sum of the monomials of codegree less than $c$ in $u$, 
and let $v_{c-1}$ be the sum of the monomials of codegree less than $c$ in $v$. We set $s_c$ to be the sum of all the monomials of 
codegree $c$ in $v_{c-1}^{\ell}$. 

$u=p(v)$, so we set $d_c$ to be the sum of all the codegree $c$ monomials in $p(v_{c-1})$.
By construction: $(u_{c-1}+d_c) (x_{c-1}+s_c)= (x_{c-1}+s_c) (u_{c-1}+d_c)=(v_{c-1})^{\ell+b}$ in the abelian
group: 
$G^{deg(ux)}/G^{deg(ux)-(c+1)}$. i.e., the monomials of codegree at most $c$ are identical for the 3 different products.

Recall that $x=x_{c-1}+h$, where $deg(h) \leq deg(x)-c$. We set $x=x_{c-1} +s_c + r$, where $deg(s_c)=deg(x)-c$ and $deg(r) \ leq deg(x)-c$.
Similarly, we set $u=u_{c-1}+d_c+q$ where $deg(q) \leq deg(u)-c$. Then:
$$ux=(u_{c-1}+d_c+q)(x_{c-1}+s_c+r)=xu=(x_{c-1}+s_c+r)(u_{c-1}+d_c+q)$$ in the abelian group: 
$G^{deg(xu)}/G^{deg(xu)-(c+1)}$. Since: $(u_{c-1}+d_c)(x_{c-1}+s_c)=(x_{c-1}+s_c)(u_{c-1}+d_c)$ in that abelian group, it follows
that: $u_{top}r+qx_{top}=ru_{top}+x_{top}q$ in that abelian group, where $u_{top}$ and $x_{top}$ are the top monomials
in $u$ and $x$ in correspondence. Therefore, all these monomials are products of a top degree monomial with a codegree $c$
monomial, and these can be broken precisely as in the codegree 1 case, assuming $deg(x) > deg(u)$.

As in the codegree 1 case, we are left with the case in which $deg(x) \leq deg(u)$. In that case 
we write: $u=u_{c-1}+d_c+q_c$ and $x=x_{c-1}+s_c+r_c$a as above. 
By the same argument that was used in that case in analyzing the codegree 1 monomials, the monomials of codegree $c$ in $r_c$
are precisely the monomials of codegree $c$ in $v^{\ell}+s_c$, and the induction follows for $c \leq deg(v)$. Hence,
$x=v^{\ell}$ in the abelian group: 
$G^{deg(x)}/G^{deg(x)-deg(v)}$. Therefore, $x+v^{\ell}$ commutes with $v$ and $u$ and $deg(x+v^{m}) < deg(x)$, and the theorem follows.

\line{\hss$\qed$}

It is possible to use the techniques that we used in this section analyze centralizers of general elements with monomial top
homogeneous part, and centralizers of general elements, but we won't need to apply these techniques in this generality in the sequel,
so we omit these generalizations.

\vglue 1.5pc
\centerline{\bf{\S4. Equations with a single variable}}
\medskip
In the previous section we gave combinatorial proofs to special cases of Bergman's theorem on the
structure of centralizers in free associative algebras. Such  combinatorial proofs are needed in
order to study the set of solutions to related systems of equations that play a central role
in understanding the set of solutions to a general monomial system of equations.

In this section we  study   the set of solutions to monomial systems of equations with a single variable.
As will be demonstrated in the sequel the techniques that are used in this section, play an essential role in
studying monomial equations with no quadratic nor free parts.

Recall that over free groups and semigroups  equations with a single variable were analyzed long before the analysis
of general systems of equations by Lorenc [Lo] and Appel [Ap]. The approach we use combines the technique and results
for studying equations with a single variable over a free group and semigroup with the combinatorial
approach that we used in analyzing centralizers and quadratic equations.

\vglue 1.5pc
\proclaim{Lemma 4.1} Let $u,v \in  FA$ and suppose that the top homogeneous parts of $u$ and $v$ are monomials
that are not proper powers.

If the equation $ux=xv$ has a non-trivial solution, then the set of solutions to the equation $ux=xv$ is a set $\{wp(v)\}$, where
$uw=wv$ and $p$ is an arbitrary polynomial in a single variable. Furthermore, the element $w$, which is the solution of minimal degree of 
the equation,  is unique. 
\endproclaim

\nfp First, the set of solutions of the equation $ux=xv$ is a linear subspace of $FA$. If $w_1$ and $w_2$ are solutions to the
equation $ux=xv$, and they are of the same degree, then their top homogeneous monomials are identical. Hence, $w_1+w_2$, which is also a solution
of that equation has strictly smaller degree than $w_1$ and $w_2$. Therefore, if the equation $ux=xv$ has a solution, then it has a unique solution
of minimal degree, that we denote $w$.

If $x_0$ is an arbitrary solution of $ux=xv$, then there exists some non-negative integer $b$, such that $wv^b$ and $x_0$ have the same top monomial.
Since both $x_0$ and $wv^b$ are solutions of the equation $ux=xv$, $x_0+wv^b$ is a solution of this equation and: $deg(x_0+v^b) < deg(x_0)$. Hence,
the proof of  the lemma follows by induction on the degree of the solution $x_0$. 

\line{\hss$\qed$}

Unlike the case of free groups or semigroups, the equation $ux=xv$ may have a solution, and still it can be
that there are no solutions with $deg(x) \leq deg(u)=deg(v)$. 

Let $t$, $\mu$, and $\rho$  be arbitrary elements in the algebra $FA$. Let $w=t \mu t \rho t \mu t$,
$v=(\rho t \mu + \mu t \rho t \mu) t$ and $u=t(\mu t \rho + \mu t \rho t \mu)$. Then: $uw=wv$ and in general
there
is no element $y \in FA$, such that $deg(y) \leq deg(u)=deg(v)$ and $uy=yv$. 

To bound the degree of a minimal degree solution we need the following lemma.

\vglue 1.5pc
\proclaim{Lemma 4.2} Let $u,v \in  FA$ be as in lemma 4.1, and suppose that the equation $ux=xv$ has a non-trivial solution.
Then there exists a solution, $w$, $uw=wv$, and $deg(w) \leq deg(u) \cdot (2^{deg(u)}+2)$.
\endproclaim

\nfp Suppose that $x_1 \neq 0$ satisfies $ux_1=x_1v$. If $deg(x_1) \leq deg(u)\cdot (2^{deg(u)}+2) $ the lemma follows. Hence, we may assume
that $deg(x_1)>deg(u)\cdot (2^{deg(u)}+2) $.

We use the analysis that was applied in analyzing centralizers in the previous section.  
By the analysis of homogeneous elements,
the top degree homogeneous part of $x_1$ must be a monomial. Let $u_0$, $v_0$, and $x_0$ be the top monomials of $u$, $v$, and $x_1$. Then they must satisfy:
$u_0x_0=x_0v_0$. Therefore, there exists a monomial $z_0$ such that: $x_0=u_0z_0=z_0v_0$.

As in analyzing centralizers, we continue the analysis of $x_1$ by analyzing its codegree 1 monomials. We examine the codegree 1 monomials in the products:
$ux_1$ and $x_1v$. We get an element $z$, such that: $x_1=zv=uz$ in the abelian group: 
$G^{deg(x_1)}/G^{deg(x_1)-2}$. 

We continue iteratively analyzing products of codegree 2. Note that products of codegree 2 that include monomials of codegree 0 and 1
of $u$, $v$ and $z$, that correspond to codegree 1 monomials of $u$ and $v$ and codegree 1 monomials of $x_1$ 
(from the two sides of the equation), cancel in pairs. A product of a codegree 1 monomial of $z$ with a codegree 1
monomial of $v$, that corresponds to a products of a codegree 1 monomial of $x_1$ with a codegree 1 monomial of $v$, cancels with either: 
\roster
\item"{(1)}" a product of a codegree 1 monomial of $u$ with a codegree 1 monomial of $x_1$.

\item"{(2)}" a product of a codegree 2 monomial of $x_1$ with the top monomial of $v$.

\item"{(3)}" a product of the
top monomial of $u$ with a codegree 2 monomial of $x_1$. 
\endroster

Suppose that the given codegree 2 product, which is a product of codegree monomials of $z$ and $v$, is equal only to a codegree 2 product
of type (3). In that case we do not add anything to $z$. If it is equal only to a codegree 2 product of type (1) we add a codegree 2 monomial
to $z$.   
If it is equal only to a codegree 2 product of type (2) we add a codegree 2 monomial
to $z$ as well.   

There must be an even number of codegree 2 products that are equal to the given codegree 2 product of a codegree 1 monomial of $x_1$
with a codegree 1 monomial of $v$. Hence, the only remaining possibility is that the given codegree 2 product is equal to codegree
n that case we do not add anything to $z$.

A pair of canceling codegree 2 products, that are products of codegree 2 monomial of $u$ and $v$ with codegree 2 monomials of $x_1$ (from the two
sides of the equation), that are not equal to any other codegree 2 products, contribute a codegree 2 monomial to $z$.

After possibly adding codegree 2 monomials to $z$, the equation that was valid for codegree 1 products is now valid for codegree 2 products, i.e.
$x_1=uz=zv$ in the abelian group:
$G^{deg(x_1)}/G^{deg(u)-3}$. 

We continue iteratively to construct the element $z$ by adding higher codegree monomials, so that the constructed element $z$ satisfies the equation: $x_1=uz=zv$ 
for products of higher and higher codegree.
Suppose that: $x_1=uz=zv$ in the abelian group:
$G^{deg(x_1)}/G^{deg(x_1)-d}$, i.e., that the equation holds for all products of codegree at most $d-1$,
where $d$ is a positive integer, $d \leq (deg(x_1)-deg(u))$. We iteratively add codegree $d$ monomials to
$z$ so that the equalities hold for all codegree $d$ products as well.

As in analyzing codegree 2 products, products of codegree d that include monomials of codegree smaller than $d$
of $u$, $v$ and $z$, that correspond to smaller codegree monomials of  $x_1$ 
(from the two sides of the equation), cancel in pairs. 

Suppose that  a codegree $d$ product, can be presented as either:
\roster
\item"{(1)}"  products of the top monomial of $u$ with
codegrees $m_i$ monomials of $z$ with codegrees $\ell_i$ monomials of $v$, for
some subset of tuples $(m_i,\ell_i)$, where $m_i+\ell_i=d$ for every index $i$.  

\item"{(2)}" 
products of 
codegrees $s_j$ monomials of $u$ with codegrees $t_j$ monomials of $z$ and with the top monomial of $v$, for
some subset of tuples $(s_j,t_j)$, where $s_j+t_j=d$ for every index $j$.  

\item"{(3)}" a product of a codegree d monomial of $x_1$ with the top monomial of $v$.

\item"{(4)}" a product of the
top monomial of $u$ with a codegree d monomial of $x_1$. 
\endroster

Suppose that a given codegree $d$ product, can be presented in the form (1) for an odd set of tuples $(m_i,\ell_i)$. If this codegree $d$
product can be presented only in the form (2), and this must be for an odd set of tuples $(s_j,t_j)$, we add a codegree $d$ monomial to $z$.
If this codegree $d$ product can be presented  in the form (4), and possibly in the form (2) in an even number of ways, we add nothing to $z$. 
If this codegree $d$ product can be presented  in the 
form (3), and possibly in the form (2) in an even number of ways, we add a codegree $d$ monomial to $z$.
If this given codegree $d$ product can be presented in an odd number of ways in form (2), and in forms (3) and (4), we add nothing to $z$. 

Suppose that  a given codegree $d$ product can be presented in an even (or $0$) number of ways in both forms (1) and (2).
If it can not be presented in form (3) and (4) we add nothing to $z$. If it can be presented in forms (3) and (4) we add a codegree $d$
monomial to $z$.

After possibly adding these codegree $d$ monomials to $z$, the equation that was valid for all products up to codegree $d-1$  is now valid for codegree $d$ products, i.e.,
$x_1=uz=zv$ in the abelian group:
$G^{deg(x_1)}/G^{deg(x_1)-(d+1)}$. 

Finally, we get an element $z$ that satisfies $x_1=uz=zv$ in the abelian group:
$G^{deg(x_1)}/G^{deg(u)-1}$, i.e., the equalities hold for all products up to degree $deg(u)=deg(v)$. 

We continue by looking at the equality: $uz=zv$ in the abelian group:
$G^{deg(x_1)}/G^{deg(u)-1}$. Repeating the same argument we can find an element $z_2$, such that:
$z=uz_2=z_2v$ in the abelian group: 
$G^{deg(z)}/G^{deg(u)-1}$. Continuing inductively, we get an element $z_{r+1}$, such that: 
$z_r=uz_{r+1}=z_{r+1}v$ in the abelian group: 
$G^{deg(z_r)}/G^{deg(u)-1}$. 

Therefore, there exist elements of distinct degree: $\{s_m, \ m=1,\ldots,2^{deg(u)}+1\}$, such that $deg(s_m) \leq (1+m) deg(u)$ and:
$us_m=s_mv$ in the abelian group:
$G^{deg(us_m)}/G^{deg(u)-1}$. 

By a simple pigeon hole argument there exist a subcollection of the indices: $1 \leq i_1 < \ldots < i_f \leq 2^{deg(u)}+1$, 
such that: $s=s_{i_1}+\ldots+s_{i_f}$ and $us=sv$. Hence, $s$ is a solution of the given equation,
and $deg(s) \leq deg(u) \cdot (2+2^{deg(u)})$.

\line{\hss$\qed$}

So far we assumed that the top degree elements of $u$ and $v$ are monomials that are not a proper power. First, we omit
the periodicity assumption,
and allow the top degree monomials of $u$ and $v$ to have non-trivial roots.

\vglue 1.5pc
\proclaim{Lemma 4.3} Let $u,v \in  FA$ and suppose that the top homogeneous parts of $u$ and $v$ are monomials.
Suppose that the top degree monomials of $u$ have non-trivial roots of degree bounded by $q$.

Suppose that the equation $ux=xv$ has a non-trivial solution. Then there exists elements $w_1,\ldots,w_d$, $1 \leq d \leq q$, such
that the set of solutions to the equation $ux=xv$ is a set of the form: $\{w_1p_1(v)+ \ldots +w_dp_d(v)\}$, where 
$uw_i=w_iv$ and $p_1,\ldots,p_d$ are arbitrary polynomials in $v$. Furthermore, the elements $w_i$ satisfy:
$deg(w_i) \leq deg(u) \cdot (2+2^{deg(u)})$.
\endproclaim

\nfp Let $u_0$ and $v_0$ be the top monomials of $u$ and $v$, and let $x_0$ be the top monomial of a solution $x$.
Let $t_0$ be a primitive root of $v_0$. Then there exists some fixed element $s_0$, $deg(s_0) < deg(t_0)$, such
that: $x_0=s_0t_0^m$ for some non-negative integer $m$.

Suppose that $t_0^q=v_0$. The top monomial of a solution $x$ is of the form $x_0=s_0t_0^m$, so we can divide the
solutions $x_1$ according to the residue classes of the non-positive integers $m$ modulo $q$. For each residue class
for which there is a solution, we fix one of the shortest solutions in the class. We denote these shortest solutions,
$w_1,\ldots,w_d$, for some $d$, $1 \leq d \leq q$.

Let $x$ be a solution. $x$ must be in the same class as one of the fixed shortest solutions $w_i$. Hence, for some
non-negative integer $b$, $x$ and $w_iv^b$ are both solutions and they have the same top monomial. Therefore,
if $x \neq w_i$, $x+w_iv^b$ is a non-trivial solution and: $deg(x+w_iv^b) < deg(x)$. By a finite induction
$x=w_1p_1(v)+ \ldots +w_dp_d(v)$, for some polynomials $p_1,\ldots,p_d$. 

The bound on the lengths of the elements $w_1,\ldots,w_d$, that were assumed to be shortest in their residue classes,
follows by the same argument that was used to prove lemma 4.2.

\line{\hss$\qed$}

So far we analyzed the equation $ux=xv$. We use similar methods to analyze the more general
equation: $u_1xu_2=v_1xv_2$.

\vglue 1.5pc
\proclaim{Theorem 4.4} Let $u_1,u_2,v_1,v_2 \in  FA$ and suppose that the top homogeneous parts of $u_i$ and $v_i$ are monomials
with no periodicity (i.e., the top monomials in $u_i$ and $v_i$ contain no subwords $\alpha^2$ for some non-trivial word $\alpha$),
and that $deg(u_1) > deg(v_1)$.

Suppose that the equation $u_1xu_2=v_1xv_2$ has a solution of degree bigger than 
$2 \cdot (deg(u_1)+deg(v_2))^2$.  Then:

\roster
\item"{(1)}" there exist elements $s,t \in FA$ such that: $u_1=v_1s$ and $v_2=tv_2$.

\item"{(2)}" an element $x \in FA$ is a solution to the equation: $u_1xu_2=v_1xv_2$ iff it is a solution of the equation $sx=xt$. 
\endroster
\endproclaim

\nfp First, note that if (1) is true and $x$ satisfies: $u_1xu_2=v_1xv_2$, then: $v_1sxu_2=v_1xtu_2$. Hence, $sx=xt$. Conversely,
every solution of the equation: $sx=xt$ satisfies $u_1xu_2=v_1xv_2$, so (2) is true.

As we did in analyzing centralizers and analyzing the equation $ux=xv$, we analyze the homogeneous parts in $x$ and in $u_i$ and $v_i$
going from top to bottom. Let $u_i^0$ and $v_i^0$, $i=1,2$, be the top monomials in $u_i$ and $v_i$. Let $x_1$ be a solution
of the equation: $u_1xu_2=v_1xv_2$, and suppose that 
$deg(x_1)>
max(deg(u_1),deg(v_2)) +2(deg(u_1)-deg(v_1))$.  
By our analysis of homogeneous
solutions, the top homogeneous part of the solution $x_1$ must be a monomial as well that we denote $x_0$.

Since $u_1^0x_0u_2^0=v_1^0x_0v_2^0$, there exists monomials $s_0,t_0$, $deg(s_0)=deg(t_0)$, such that: $u_1^0=v_1^0s_0$,
$v_2^0=t_0u_2^0$, and $x_0=f_0t_0^b=s_0^be_0$, for some positive integer $b$, and: $deg(f_0)=deg(e_0)<deg(s_0)$.

We continue by analyzing monomials of codegree 1 in $u_i$, $v_i$ and $x_1$. By the same analysis that was used in
analyzing centralizers and in lemma 4.2, there exist elements $s,t$ with top monomials $s_0$ and $t_0$, and an element
$w$ with top monomial $w_0$, $w_0t_0=s_0w_0=x_0$, such that:
\roster
\item"{(i)}" $u_1=v_1s$ in the abelian group: 
$G^{deg(u_1)}/G^{deg(u_1)-2}$. 

\item"{(ii)}" $v_2=tu_2$ in the abelian group: 
$G^{deg(v_2)}/G^{deg(v_2)-2}$. 

\item"{(iii)}" $x_1=sw=wt$ in the abelian group: 
$G^{deg(x_1)}/G^{deg(x_1)-2}$. 
\endroster

By iteratively applying the same construction, the above 3 equalities imply that $w=s^mf=et^m$ in the abelian group:
$G^{deg(w)}/G^{deg(w)-2}$, for some positive integer $m$, and elements $e,f$, $deg(s) \leq deg(e)=deg(f) < 2deg(s)$. 

\smallskip
We continue by analyzing products of codegree 2. First, note that as in analyzing centralizers, if we look at codegree 2 products that involve only
top monomials and codegree 1 monomials from $s$, $t$, $u_i$, $v_i$ and $w$, such that the products restrict to codimension 0 or 1 monomials of $x_1$, $u_i$ and $v_i$, then
such codegree 2 products cancel in pairs from the two sides of the equation.

We further look at codegree 2 products that contain a codegree 1 monomial of $u_2$. If the codegree 2 product contains the top monomial of $t$, then such a codegree 2 product
cancels with a corresponding codegree 2 product from the other side of the equation, since all the corresponding monomials of $u_i$, $v_i$, and $x_1$ (from the two sides of the
equation) are either codegree 0 or codegree 1. Hence, we look at codegree 2 products that contain codegree 1 monomials of $t$ and $u_2$, and, therefore, top monomials
of $u_1$ and $w$. Such a codegree 2 product, which is a product of the top monomial of $u_1$, a codegree 1 monomial of $x_1$ and a codegree 1 monomial of
$u_2$,   cancels with either: 
\roster
\item"{(1)}" a product of the top monomial of $v_1$, the top monomial of $x_1$, and a codegree 2 monomial of $v_2$.

\item"{(2)}" a product of the top monomial of $v_1$, a codegree 1 monomial of $x_1$ and a codegree 1 monomial of $v_2$.

\item"{(3)}" a product of the top monomial of $v_1$, a codegree 2 monomial of $x_1$ and the top monomial of $v_2$.

\item"{(4)}" a product of the top monomial of $u_1$, a codegree 2 monomial of $x_1$ and the top monomial of $u_2$.

\item"{(5)}" a product of the top monomial of $u_1$, the top monomial of $x_1$ and a codegree 2 monomial of $u_2$.
\endroster

If the given codegree 2 product cancels only with a product of type (1) we don't add anything to $w$ nor to $t$. Suppose that the given
codegree (2) product cancels only with a product of type (2). If the codegree 1 monomial of $v_2$ equals the top monomial of $t$ with a codegree
1 monomial of $u_2$, then the codegree 2 product of type (2) cancels with a codegree 2 product from the other side of the equation that
contains only codegree 0 and 1 monomials of $u_1$, $x_1$ and $u_2$. Hence, we can assume that the codegree 1 monomial of $v_2$ is a product of
a codegree 1 monomial of $t$ with the top monomial of $u_2$. In that case we add a codegree 2 monomial to $t$ and leave $w$ unchanged.

If the given codegree 2 product cancels only with a codegree 2 product of type (3), we add a codegree 2 monomial to $t$ and a codegree 2 monomial
to $w$. If the given codegree 2
product cancels only with a product of type (4) we add a codegree 2 monomial to $t$. In case the given product cancels only with a codegree 2 product of type
(5) we don't add anything (apart from the codegree 2 monomial of $u_2$).

A products of type (3) can not cancel with a product of type (5). Hence, the only left possibilities are a collection of products of 3 different types that
cancel with the given codegree 2 product. We list the various possibilities for the collections of codegree 2 products of 3 different types  that cancel
with the given codegree 2 product and indicate what we add in each possibility:
\roster
\item{(i)}" products (1)-(3) cancel. We add a codegree 2 monomial to $w$, apart from an existing codegree 2 monomial of $v_2$ (that is equal to the
products of the codegree 1 monomials of $t$ and $u_2$ in the given codegree 2 product).  

\item"{(ii)}" products (1),(2), and (4). In that case we don't add anything to $w$ and $t$. A monomial of codegree 2 that already appears in $v_2$,
and is equal to the product of the given codegree 1 monomials of $t$ and $u_2$.

\item"{(iii)}" products (1),(2) and (5). We add a codegree 2 monomial to $t$, in addition to the codegree 2 monomials that already appear in 
$u_2$ and $v_2$.

\item"{(iv)}" products (1),(3) and (4). We add a codegree 2 monomial to $w$, and the existing codegree 2 monomial to $v_2$.

\item"{(v)}" products (2),(3) and (4). We add a codegree 2 monomial to $w$ and a codegree 2 monomial to $t$.

\item"{(vi)}" products (1), (4) and (5). We add a codegree 2 monomial to $t$, and the existing codegree 2 monomials to $u_2$ and $v_2$.

\item"{(vii)}" products (2),(4) and (5). We just add the existing codegree 2 monomial to $u_2$.
\endroster

So far we analyze codegree 2 products that cancel with a given codegree 2 product that is a product of the top monomial of $u_1$ and $w$ and codegree
1 monomials of $t$ and $u_2$. We continue by analyzing those codegree 2 products that cancel with a given codegree 2 product of the
top monomial of $v_1$, a codegree 1 monomial of $x_1$ that is a product of the top monomial of $s$ and a codegree 1 monomial of $w$, and a codegree 1 monomial of $v_2$,
 that is equal to a product of a codegree 1
monomial of $t$ with the top monomial of $u_2$. Such a given codegree 2 product can cancel with either:
\roster
\item"{(1)}" a product of the top monomial of $v_1$, a codegree 2 monomial of $x_1$, and the top monomial of $v_2$.

\item"{(2)}" a product of the top monomial of $v_1$, the top  monomial of $x_1$ and a codegree 2 monomial of $v_2$.

\item"{(3)}" a product of the top monomial of $u_1$, a codegree 2 monomial of $x_1$ and the top monomial of $u_2$.

\item"{(4)}" a product of the top monomial of $u_1$, the top monomial of $x_1$ and a codegree 2 monomial of $u_2$.

\item"{(5)}" a product of the top monomial of $u_1$, a codegree 1 monomial of $x_1$ and a codegree 1 monomial of $u_2$.

\item"{(6)}" a product of a codegree 1 monomial of $u_1$, which is a product of the top monomial of $v_1$ with a codegree
1 monomial of $s$, with a codegree 1 monomial of $x_1$, which is the product of a codegree 1 monomial of $w$ with the top monomial of $t$, and with the
top monomial of $u_2$.
\endroster

If the given codegree 2 product equals only to a codegree 2 product of type (1), we add a codegree 2 monomial to $w$. If it equals only to a codegree 2
product of type (2), we add a codegree 2 monomial to $t$, apart from the existing codegree 2 monomial of $v_2$. 
If it equals only to a codegree 2 product of type (3), we do not add anything. If it equals only
to a codegree 2 product of type (4) we add a codegree 2 monomial to $t$, apart from the existing codegree 2 monomial of $u_2$. 
We already analyze all the codegree 2 products that cancel with a codegree 2 product of
type (5), so we omit this case. If it equals only to a product of type (6) we add a codegree 2 monomial to $w$.

A codegree 2 product of types (1) or (6)  can not cancel with a codegree 2 product of type (4). A monomial of type (5) that cancels with a monomial
of type (6) is a product of lower codegree monomials of $u_i$, $v_i$ and $x_1$ from the two sides of the equation, so we omit this case. Hence, there are 3 cases left: 
\roster
\item"{(i)}" products (1), (2) and (3) cancel with the given codegree 2 product. In that case we add a codegree 2 monomial to $w$ and a codegree 2 monomial to $t$, apart
from the existing codegree 2 monomial of $v_2$.

\item"{(ii)}" products (2), (3), and (4). In that case we only add the already existing codegree 2 monomials of $u_2$ and $v_2$.

\item"{(iii)}" products (1), (3) and (6). In that case we do not add anything.
\endroster

Codegree 2 products that contain codegree 1 monomials of $v_1$ or $u_1$ are treated exactly in the same way. Hence, we are left with sets of codegree 2 products
that cancel, and each of these codegree 2 products is a product of top monomials with codegree 2 monomials of one of the $u_i$, $v_i$ or $x_1$. These are analyzed precisely as they
treated in the proof of lemma 4.4 and in analyzing codegree 1 products, and in each such cancellation codegree 2 monomials may be added to either $s$, $t$ or $w$, apart from existing
codegree 2 monomials of $u_i$ and $v_i$. Finally, we (possibly) added codegree 2 monomials to $s$, $w$ and $t$, such that: 
\roster
\item"{(i)}" $u_1=v_1s$ in the abelian group: 
$G^{deg(u_1)}/G^{deg(u_1)-3}$. 

\item"{(ii)}" $v_2=tu_2$ in the abelian group: 
$G^{deg(v_2)}/G^{deg(v_2)-3}$. 

\item"{(iii)}" $x_1=sw=wt$ in the abelian group: 
$G^{deg(x_1)}/G^{deg(x_1)-3}$. 
\endroster

We continue iteratively with products with higher codegree.
Let $d=min(deg(v_1),deg(u_2),deg(u_1)-deg(v_1))$. Let $r \leq d-1$ and suppose that we added codegree $r$ monomials to $s$, $w$ and $t$ such that the equations above hold for all coproducts 
of codegree bounded by $r-1$.

We analyze codegree $r$ products in the same way we analyzed codegree 2 products. First, note that if a codegree $r$ product is a product of monomials of $u_i$, $v_i$, $s$, $t$ and $w$,
that correspond to products of monomials of codegree smaller than $r$ of $u_i$, $v_i$, and $x_1$ from the two sides of the equation, then such codegree $r$ products cancel in pairs.

Suppose that a codegree $r$ product is a product of the top monomials of $u_1$ and $w$, and monomials of codegree $q_i$ of $t$ and codegree $m_i$ of $u_2$, such that
$q_i+m_i=r$  and  $q_i,m_i$ are positive integers, and there are odd number of such pairs $(q_i,m_i)$. We treat this case in the same way we treated the case of a codegree 2
product that includes a codegree 1 monomial of $t$ and a codegree 1 monomial of $u_2$.
This odd set of codegree $r$ products (that are all equal)
cancels with either: 
\roster
\item"{(1)}" a product of the top monomial of $v_1$, the top monomial of $x_1$, and a codegree r monomial of $v_2$.

\item"{(2)}" an odd set of codegree $r$  products of the top monomial of $v_1$,  codegree $e_j$ monomial of $x_1$ and  codegree $p_j$ monomial of $v_2$, for
some positive set of pairs $(e_j,p_j)$ that satisfy $e_j+p_j=r$, and such that the codegree $p_j$ monomial of $v_2$ is a product of a codegree $p_j$ monomial of
$t$ with the top monomial of $u_2$.

\item"{(3)}" a product of the top monomial of $v_1$, a codegree r monomial of $x_1$ and the top monomial of $v_2$.

\item"{(4)}" a product of the top monomial of $u_1$, a codegree r monomial of $x_1$ and the top monomial of $u_2$.

\item"{(5)}" a product of the top monomial of $u_1$, the top monomial of $x_1$ and a codegree r monomial of $u_2$.
\endroster

The treatment of the various cases is identical to what we did in analyzing codegree 2 products (cases (i)-(vii)), just that instead of adding codegree
2 monomials to the various elements we add  codegree $r$ monomials. The case in which a codegree $r$ product is obtained in an odd number of ways
as the product of the top monomials $v_1$ and $s$, codegree $e_j$ monomial of $w$ and codegree $p_j$ of $v_2$ is treated in an identical way to the analysis
of a codegree 2 product that contains a codegree 1 monomial of $w$ and a codegree 1 monomial of $v_2$, which is the product of a codegree 1 monomial of $t$ with 
the top monomial of $u_2$.  

As in the case of codegree 2 products, we still need to treat all the codegree $r$ products that cancel with a product that can be presented in an odd number of
ways as the product of
the top monomial of $v_1$,  codegree $e_j$ monomial of $x_1$ and  codegree $p_j$ monomial of $v_2$, for
some positive set of pairs $(e_j,p_j)$ that satisfy $e_j+p_j=r$, and such that the codegree $p_j$ monomial of $v_2$ is the 
product of a codegree $p_j$ monomial of $t$ with the top monomial of $u_2$. The treatment of this case is also identical to the treatment of such a codegree
2 product.




Therefore, we constructed elements $s,t,w$ for which:
\roster
\item"{(i)}" $u_1=v_1s$ in the abelian group: 
$G^{deg(u_1)}/G^{deg(u_1)-d}$. 

\item"{(ii)}" $v_2=tu_2$ in the abelian group: 
$G^{deg(v_2)}/G^{deg(v_2)-d}$. 

\item"{(iii)}" $x_1=sw=wt$ in the abelian group: 
$G^{deg(x_1)}/G^{deg(x_1)-d}$. 
\endroster

We divide the continuation according to minimum between $deg(v_1)$, $deg(u_2)$ and $deg(u_1)-deg(v_1)$. 
First we assume that: 
$d=min(deg(v_1),deg(u_2),deg(u_1)-deg(v_1))=deg(u_1)-deg(v_1)$.

In analyzing codegree $d$ products, there are special codegree $d$ products that we need to single out and treat separately,
 as they may involve cancellations between codegree $d$ products that contain  codegree $d$   monomials of $u_1$ or $v_1$ and
those that contain codegree $d$ monomials of $u_2$ or $v_2$.

As in analyzing smaller codegree products, note that codegree $d$ products that are products of smaller codegree monomials of
the $u_i$, $v_i$, $s$, $w$ and $t$, and correspond to smaller codegree monomials of $u_i$, $v_i$ and $x_1$ from the two
sides of the equation cancel in pairs.

We continue by analyzing codegree $d$ products that are products of top degree monomials of $u_1$ and $w$ codegree $q_i$ monomials of $t$ and 
codegree $m_i$ monomials of $u_2$, such that $q_i+m_i=d$, there are odd number of such pairs $(q_i,m_i)$, and the product of these monomials of $t$ and $u_2$
is not equal to $u_2^0$, the top monomial of $u_2$. 

Such codegree $d$ products are analyzed exactly in the same way they were analyzed in codegree $r$ products for $r<d$. Similarly, we analyze codegree $d$
products that are obtained odd number of times as the product of the top monomials of $v_1$ and $s_1$, a codegree $e_j$ monomial of $w$ and a codegree $p_j$
monomial of $v_2$, such that the product of a codegree $e_j$ monomial of $w$ and a codegree $p_j$ monomial of $v_2$ is not $w_0u_2^0$, i.e., the product
of the top monomials of $w$ and $u_2$.     

In a a similar way we analyze codegree $d$ products that are products of smaller codegree monomials of $u_1$ and $w$ and the top monomials of $t$ and $u_2$, and 
products of smaller codegree monomials of $v_1$ and $s$ and the top monomials of $w$ and $v_2$, assuming the products of these smaller degree monomials are
not equal to $v_1^0w_0$ or to $v_1^0$.  

We continue by analyzing canceling pairs of codegree $d$ products that are products of top monomials of $v_i$, $u_i$ and $x_1$,  with one codegree $d$ monomials 
of these elements, such that this codegree $d$ monomial of $u_1$ is not $v_1^0$, the codegree $d$ monomial of $x_1$ is not $w_0$ and the codegree $d$ monomial of
$v_2$ is not $u_2^0$. These codegree $d$ products are analyzed in the same way they were analyzed for smaller codegree products.

We are left with codegree $d$ products that are either:
\roster
\item"{(1)}" a product of a codegree $d$ monomial of $u_1$, that is equal to $v_1^0$, with top monomials of $x_1$ and $u_2$.

\item"{(2)}" a product of the top monomial of $v_1$ with top monomials of $x_1$ and a codegree $d$ monomial of $v_2$, that is equal to $u_2^0$.

\item"{(3)}" a product of the top monomial of $u_1$ with a codegree $d$ monomial of $x_1$, that is equal to the top monomial of $w$, with
the top monomial of $u_2$.

\item"{(4)}" a product of the top monomial of $v_1$ with a codegree $d$ monomial of $x_1$, that is equal to the top monomial of $w$, with
the top monomial of $v_2$.

\item"{(5)}" 
 an odd set of codegree $d$  products of the top monomial of $v_1$,  codegree $e_j$ monomial of $x_1$ and  codegree $f_j$ monomial of $v_2$, for
some positive pairs $(e_j,p_j)$ that satisfy $e_j+p_j=d$, such that the monomial of $x_1$ is a product of the top monomial
of $s$ with a codegree $e_j$ monomial of $w$, and the product of each codegree $e_j$ monomial of $x_1$ with a codegree $p_j$
monomial of $v_2$ is equal to $w_0v_2^0$ (the product of the top monomials of $w$ and $v_2$).

\item"{(6)}" 
 an odd set of codegree $d$  products of the top monomial of $u_1$,  codegree $q_i$ monomials of $x_1$, that are products of the top monomial of $w$ with
a codegree $q_i$ monomial of $t$, with  codegrees $m_i$ monomials of $u_2$, for
some positive set of pairs $(q_i,m_i)$ that satisfy $q_i+m_i=d$, such that the product of each codegree $q_i$ monomial of $t$ with a codegree $m_i$
monomial of $u_2$ is equal to $u_2^0$. 

\item"{(7)}" 
 an odd set of codegree $d$ products of codegree $f_i$ monomials of $u_1$ and  codegree $g_i$ monomials of $x_1$ with the top monomial of $u_2$, for
some positive  pairs $(f_i,g_i)$ that satisfy $f_i+g_i=d$, such that the monomial of $x_1$ is a product of a codegree $g_i$ monomial of $w$ with the top monomial of $t$,
and the product of each codegree $f_i$ monomial of $u_1$ with a codegree $g_i$
monomial of $x_1$ is equal to $u_1^0w_0$ (the product of the top monomials  of $u_1$ and $w$).

\item"{(8)}" 
 an odd set of codegree $d$ products of codegree $h_j$ monomials of $v_1$ and  codegree $k_j$ monomials of $x_1$, that are products of
a codegree $k_j$ monomial of $s$ with the top monomial of $w$, with the top monomial of $v_2$, for
some positive  pairs $(h_j,k_j)$ that satisfy $h_j+k_j=d$, such that the product of each codegree $h_j$ monomial of $v_1$ with a codegree $k_j$
monomial of $s$ is equal to $v_1^0$.
\endroster

First note that (3) exists if and only if (4) exists and they cancel each other. If (3) and (4) are the only existing possibilities, we add a codegree $d$ monomial
to $w$, which is the codegree $d$ prefix or suffix of the top monomial of $w$. Also note that if cases (1) or (2) exist, codegree $d$ monomials that already appear in 
$u_1$ or $v_2$ are added to them.
Suppose that only 2 of the possibilities (1), (2) and (5)-(8) exist, possibly in addition to (3) and (4). We go over the various alternatives:
\roster
\item"{(i)}" If only (1) and (2) exist, we add the constant element 1 to $s$ and $t$, and the codegree $d$ prefix of $w_0$ to $w$, where $w_0$ is the top monomial of $w$.
If (3) and (4) exists as well, we only add 1 to $s$ and $t$.

\item"{(ii)}" If only (5) and (6) exist, we just add 1 to $t$. If (3) and (4) exist as well, we add the codegree $d$ prefix of $w_0$ to $w$. The case in which
only (7) and (8) exist is identical.

\item"{(iii)}" If only (5) and (8) exist, we add 1 to $s$ and the codegree $d$ prefix of $w_0$ to $w$. If (3) and (4) exist as well, we only add 1 to $s$. The case in which
only (6) and (7) exist is treated identically.

\item"{(iv)}" If only (5) and (7) exist, we add the codegree $d$ prefix of $w_0$ to $w$. If (3) and (4) exist as well, we do not add anything to any of the variables. 

\item"{(v)}" If only (6) and (8) exist, we add $1$ to $s$ and $t$ and the codegree $d$ prefix of $w_0$ to $w$. If (3) and (4) exist as well, we only add 1 to $s$ and $t$.

\item"{(vi)}" If only (1) and (8) exist, we add 1 to $s$ and the codegree $d$ prefix of $w_0$ to $w$. If (3) and (4) exist as well we only add 1 to $s$.

\item"{(vii)}" If only (1) and (7) exist, we do the same as in (v), we add 1 to $s$ and the codegree $d$ prefix of $w_0$ to $w$. If (3) and (4) exist as well we only add 1 to $s$.

\item"{(viii)}" If only (1) and (6) exist, we add 1 to $s$ and $t$, and the codegree $d$ prefix of $w_0$ to $w$. If (3) and (4) exist as well we only add 1 to $s$ and $t$.

\item"{(ix)}" If only (1) and (5) exist, we add 1 to $s$ and the codegree $d$ prefix of $w_0$ to $w$. If (3) and (4) exist as well we only add 1 to $s$.
\endroster

The cases in which only case (2) and one of the cases  (5)-(8) exist, are treated according to cases (vi)-(ix). Suppose that exactly four of the cases (1)-(2) and (5)-(8) exist,
possibly in addition to (3) and (4). We go over the alternatives: 
\roster
\item"{(i)}" If only (1),(2),(5) and (6) exist, we add 1 to $s$ and the codegree $d$ prefix of $w_0$ to $w$. If (3) and (4) exist as well, we just add 1 to
$s$. The case in which
only (1),(2),(7) and (8) exist is identical.

\item"{(ii)}" If only (1),(2),(5) and (8) exist, we add 1 to $t$. If (3) and (4) exist as well, we add 1 to $t$ and the codegree $d$ prefix of $w_0$ to $w$.
The case in which
only (1),(2),(6) and (7) exist is identical.

\item"{(iii)}" If only (1),(2),(5) and (7) exist, we add 1 to $s$ and $t$. If (3) and (4) exist as well, we  add 1 to $s$ and $t$ and
the codegree $d$ prefix of $w_0$ to $w$. 

\item"{(iv)}" If only (1),(2),(6) and (8) exist, we do not change any of the variables. If (3) and (4) exist as well, we  add 
the codegree $d$ prefix of $w_0$ to $w$. 

\item"{(v)}" If only (5),(6),(7) and (8) exist, we add 1 to $s$ and $t$. If (3) and (4) exist as well, we add 1 to $s$ and $t$ and
the codegree $d$ prefix of $w_0$ to $w$. 

\item"{(vi)}" If only (1),(5),(6) and (7) exist, we add 1 to $s$ $t$. If (3) and (4) exist as well, we add 1 to $s$ and $t$ and
and the codegree $d$ prefix of $w_0$ to $w$. 

\item"{(vii)}" If only (1),(5),(6) and (8) exist, we add 1 to $t$. If (3) and (4) exist as well, we add 1 to $t$ and the codegree $d$ prefix of $w_0$ to $w$.

\item"{(viii)}" If only (1),(5),(7) and (8) exist, we add the codegree $d$ prefix of $w_0$ to $w$. If (3) and (4) exist as well, we do not add anything to any of the variables.

\item"{(ix)}" If only (1),(6),(7) and (8) exist, we add 1 to $t$ and the codegree $d$ prefix of $w_0$ to $w$. If (3) and (4) exist as well, we  just add 1 to $t$. 
\endroster

The cases in which only case (2) and 3 of the cases  (5)-(8) exist, are treated according to cases (vi)-(ix). Suppose that cases (1),(2) and (5)-(8) exist.
In that case we add the codegree $d$ prefix of $w_0$ to $w$. If cases (3) and (4) exist as well, we do not add anything to any of the variables.

\medskip
This completes the analysis of codegree $d$ products. 
We continue with the analysis of codegree $d+1$ products. First, as in analyzing smaller codegree products, codegree $d+1$ products that are products of
smaller codegree monomials of $u_i$, $v_i$, $s$, $t$, $w$, that correspond to products of smaller degree monomials of $u_i$, $v_i$ and $x_1$ from the two sides of the equation, 
cancel in pairs.  

\vglue 1.5pc
\proclaim{Lemma 4.5} 
Suppose that a codegree $d+1$ product is a product of the top monomials of $u_1$ and $w$, and monomials of codegrees $q$ of $t$ and codegree $m$ of $u_2$, such that
$q \geq 0$ and $m>0$ and $q+m=d+1$.

Such a codegree $d+1$ product can not be:
\roster
\item"{(1)}" a product of
the top monomial of $v_1$, a codegree $e$ monomial of $x_1$ and a codegree $f$ monomial of $v_2$, for $e > 0$ and $f \geq 0$,
that satisfy $e+f=d+1$, and the codegree $f$ monomial of $v_2$ is a product of a codegree $f$ monomial
of $t$ with the top monomial $u_2^0$ of $u_2$.

\item"{(2)}" a product of
the top monomial of $u_1$, a codegree $d+1$ monomial of $x_1$ and the top monomial of $u_2$.
\endroster

\endproclaim

\nfp If such a codegree $d+1$ product can be presented as a product in the forms (1) or (2), $u_2^0$ has a prefix which is a suffix of $t_0$. 
Hence, $v_2^0$ has non-trivial periodicity that contradicts
our assumptions.

\line{\hss$\qed$}

Suppose that a codegree $d+1$ product can be presented as  a product of the top monomials of $u_1$ and $w$, and monomials of codegrees $q_i$ of $t$ and codegree $m_i$ of $u_2$, 
such that $q_i \geq 0$ and $m_i>0$  and $q_i+m_i=d+1$, 
 and there are odd number of such pairs $(q_i,m_i)$. 

By lemma 4.5, the  same codegree $d+1$ product is the product of the top monomial of $v_1$, the top  monomial of $x_1$, and a codegree $d+1$ monomial of $v_2$.

Furthermore, by the same argument that was used in the proof of lemma 4.5, if a codegree $d+1$ product is the product of the top monomials of $v_1$ and $x_1$ and
a codegree $d+1$ monomial of $v_2$, then it must be the product of an odd number of products
of the top monomials of $u_1$ and $w$, and monomials of codegrees $q_i$ of $t$ and codegree $m_i$ of $u_2$, 
such that $q_i \geq 0$ and $m_i>0$  and $q_i+m_i=d+1$. 

\vglue 1.5pc
\proclaim{Lemma 4.6} 
Suppose that a codegree $d+1$ product can be presented in
an odd number of ways as products of 
the top monomial of $v_1$,  a codegree $e_j$ monomial of $x_1$, which is the product of the top monomial of $s$ with a codegree $e_j$ monomial of
$w$, and  a codegree $f_j$ monomial of $v_2$, for positive $e_j$ and $f_j$  and
$e_j+f_j=d+1$, and the codegree $f_j$ monomials of $v_2$ are products of  codegree $f_j$ monomials
of $t$ with the top monomial $u_2^0$ of $u_2$.

\roster
\item"{(1)}" Suppose that this  codegree $d+1$ product can not be presented in
an odd number of ways as products of 
a codegree $g_j$ monomial of $u_1$,
and  a codegree $h_j$ monomial of $x_1$, 
 which is a codegree $h_j$ monomial of $w$ with the top monomial of $t$,
for positive $g_j$ and $h_j$  and
$g_j+h_j=d+1$, and the codegree $g_j$ monomials of $u_1$ are products of $v_1^0$ with  codegree $g_j$ monomials
of $t$.


Then the same codegree $d+1$ product is either the product of the top monomial of $u_1$, a codegree $d+1$ monomial of $x_1$, and the top monomial of $u_2$,
or the product of the top monomial of $v_1$, a codegree $d+1$ monomial of $x_1$ and the top monomial of $v_2$, and exactly one of the two occurs. 

\item"{(2)}" Suppose that this  codegree $d+1$ product can  be presented in
an odd number of ways as products of 
a codegree $g_j$ monomial of $u_1$,
and  a codegree $h_j$ monomial of $x_1$, 
 which is a codegree $h_j$ monomial of $w$ with the top monomial of $t$,
for positive $g_j$ and $h_j$  and
$g_j+h_j=d+1$, and the codegree $g_j$ monomials of $u_1$ are products of $v_1^0$ with  codegree $g_j$ monomials
of $t$.

Then either the same codegree $d+1$ product is both the product of the top monomial of $u_1$, a codegree $d+1$ monomial of $x_1$, and the top monomial of $u_2$,
and the product of the top monomial of $v_1$, a codegree $d+1$ monomial of $x_1$ and the top monomial of $v_2$, or none of these two possibilities occur. 
\endroster
\endproclaim

\nfp Such a codegree $d+1$ products does not cancel only with the product of monomials of codegree less that $d$ of $u_i$, $v_i$ and $x_1$, from the two sides of the equation. By
lemma 4.5 such codegree $d+1$ products can not be equal to the following products:
\roster
\item"{(1)}"   top degree monomials of $u_1$ and $x_1$ with monomials of $t$ and $u_2$. 

\item"{(2)}" 
top degree monomials of $v_1$ and $x_1$ and a codegree $d+1$ monomial of  $v_2$.

\item"{(3)}"   monomials of $v_1$ and $s$ with top degree monomials of $x_1$ and $v_2$.

\item"{(4)}" 
top degree monomials of $v_1$ and $x_1$ and a codegree $d+1$ monomial of  $v_2$.
\endroster

Therefore, such a codegree $d+1$ product must be equal to an odd number of products in the forms that are listed in the statement of the lemma.

\line{\hss$\qed$}

Lemmas 4.5 and 4.6 enable us to treat codegree $d+1$ products in a similar way to the analysis of codegree $r$ products for $r<d$.

Suppose that a codegree $d+1$ product is obtained in an odd number of ways as the product of the
top monomials of $u_1$ and $x_1$ with a codegree $q_i$ monomial of $t$, and a codegree $m_i$ monomial of $u_2$, such
that $q_i \geq 0$ and $m_i>0$  and $q_i+m_i=d+1$. By lemma 4.5, such a product must be equal to a product of the top monomials of
$v_1$ and $x_1$ and a codegree $d+1$ monomial of $v_2$. 

An analogous conclusion holds if 
a codegree $d+1$ product is obtained in an odd number of ways as the product of 
a codegree $m_i$ monomial of $v_1$ with a codegree $q_i$ of $s$ with the top monomials of $x_1$ and $v_2$, such
that $q_i \geq 0$ and $m_i>0$  and $q_i+m_i=d+1$. 

Suppose that a codegree $d+1$ product can be presented in
an odd number of ways as products of 
the top monomial of $v_1$,  a codegree $e_j$ monomial of $x_1$, which is the product of the top monomial of $s$ with a codegree $e_j$ monomial of
$w$, and  a codegree $f_j$ monomial of $v_2$, for positive $e_j$ and $f_j$  and
$e_j+p_j=d+1$, and the codegree $f_j$ monomials of $v_2$ are products of  codegree $f_j$ monomials
of $t$ with the top monomial $u_2^0$ of $u_2$.

Suppose that this  codegree $d+1$ product can not be presented in
an odd number of ways as products of 
a codegree $g_j$ monomial of $u_1$,
and  a codegree $h_j$ monomial of $x_1$, 
 which is a codegree $h_j$ monomial of $w$ with the top monomial of $t$,
for positive $g_j$ and $h_j$  and
$g_j+h_j=d+1$, and the codegree $g_j$ monomials of $u_1$ are products of $v_1^0$ with  codegree $g_j$ monomials
of $t$.

If the same codegree $d+1$ product is the product of the top monomial of $u_1$, a codegree $d+1$ monomial of $x_1$, and the top monomial of $u_2$, we do
not add anything.
If it is the  product of the top monomial of $v_1$, a codegree $d+1$ monomial of $x_1$ and the top monomial of $v_2$, we add a codegree $d+1$ monomial to $w$. 

Suppose that this  codegree $d+1$ product can  be presented in
an odd number of ways as products of 
a codegree $g_j$ monomial of $u_1$,
and  a codegree $h_j$ monomial of $x_1$, 
 which is a codegree $h_j$ monomial of $w$ with the top monomial of $t$,
for positive $g_j$ and $h_j$  and
$g_j+h_j=d+1$, and the codegree $g_j$ monomials of $u_1$ are products of $v_1^0$ with  codegree $g_j$ monomials
of $t$.

If the same codegree $d+1$ product is both the product of the top monomial of $u_1$, a codegree $d+1$ monomial of $x_1$, and the top monomial of $u_2$,
and the product of the top monomial of $v_1$, a codegree $d+1$ monomial of $x_1$ and the top monomial of $v_2$, then we do not add anything.
If none of these two possibilities occur, we add a codegree $d+1$ monomial to $w$ (by lemma 4.6 either both or none occur). 

Suppose that a codegree $d+1$ product can be presented only as the product of the top monomial of $u_1$, a codegree $d+1$ monomial of $x_1$ and the top monomial
of $u_2$, and as the product of the top monomial of $v_1$, a codegree $d+1$ monomial of $x_1$, and the top monomial of $v_2$. In that case we add a codegree $d+1$ monomial to
$w$.

This concludes the analysis of codegree $d+1$ products.
The analysis of codegree $d+r$ products, $r<d$,  is identical to the analysis of codegree $d+1$ products. Hence, we (possibly) finally add codegree $d+r$ monomials to $w$,
and the existing codegree $d+r$ monomials to $u_i$ and $v_i$
for $1 \leq r <d$, and do not change $s$ and $t$, such that:  
\roster
\item"{(i)}" $u_1=v_1s$ in the abelian group: 
$G^{deg(u_1)}/G^{deg(u_1)-2d}$. 

\item"{(ii)}" $v_2=tu_2$ in the abelian group: 
$G^{deg(v_2)}/G^{deg(v_2)-2d}$. 

\item"{(iii)}" $x_1=sw=wt$ in the abelian group: 
$G^{deg(x_1)}/G^{deg(x_1)-2d}$. 
\endroster

In analyzing codegree $2d$ products, as in analyzing codegree $d$ products, there are special codegree $2d$ monomials that we need to single out and treat separately,
 as they may involve cancellations between codegree $2d$ products that contain codegree $d$ or $2d$ monomials of $u_1$ or $v_1$ and
those that contain codegree $d$ or $2d$ monomials of $u_2$ or $v_2$.

As in analyzing smaller codegree products, note that codegree $2d$ products that are products of smaller codegree monomials of
the $u_i$, $v_i$, $s$, $w$ and $t$, and correspond to smaller codegree monomials of $u_i$, $v_i$ and $x_1$ from the two
sides of the equation cancel in pairs.

As in analyzing codegree $d$ products, we continue by analyzing codegree $2d$ products that are products of top degree monomials of $u_1$ and $w$, codegree $q_i$ monomials of $t$ and 
codegree $m_i$ monomials of $u_2$, such that $q_i+m_i=2d$, there are odd number of such pairs $(q_i,m_i)$, and the product of these monomials of $t$ and $u_2$
is not equal to a codegree $d$ suffix of $u_2^0$, the top monomial of $u_2$ (which is a codegree $2d$ suffix of $v^2_0$, the top monomial of $v_2$). 
Such codegree $2d$ products must cancel with the product of the top monomials of $v_1$ and $x_1$ and a codegree $2d$ monomial of $v_2$. In this case we only add the already
existing codegree $2d$ monomial to $v_2$. 

Similarly, we analyze codegree $2d$
products that are obtained odd number of times as the product of the top monomials of $v_1$ and $s_1$, a codegree $e_j$ monomial of $w$ and a codegree $p_j$
monomial of $v_2$, such that the product of a codegree $e_j$ monomial of $w$ and a codegree $p_j$ monomial of $v_2$ does not have a suffix which is the codegree $d$
suffix of  $u_2^0$.
We analyze codegree $2d$ products that contain similar monomials of $v_1$, $u_1$, and  $x_1$ in a similar way.

Suppose that a codegree $2d$ product is obtained in an odd number of ways as the product of the
top monomials of $u_1$ and $x_1$, a codegree $q_i$ monomial of $t$, and a codegree $m_i$ monomial of $u_2$, such
that $q_i$ and $m_i$ are positive and $q_i+m_i=2d$, and the product of the monomial of $t$ and the monomial of $u_2$ is a codegree $d$ monomial of $u_2^0$.

Because we assumed that the coefficients do not have any periodicity, such a codegree $2d$ product must cancel with either a product of the top monomials of
$u_1$ and $x_1$ and a codegree $2d$ monomial of $u_2$, or a product of the top monomials of $v_1$ and $x_1$ and a codegree $2d$ monomial of $v_2$. In both of these cases we only add
the already existing codegree $2d$ monomials to $u_2$ or $v_2$. 

If $x_1$ contains a monomial that is equal to the $2d$ prefix (or suffix) of $x_1^0$, then the codegree $2d$ product that contains the top monomials of the $u_i$'s and
this codegree $2d$ monomial of $x_1$ cancels with the codegree $2d$ product of the top monomials of the $v_i$'s with that codegree $2d$ monomial of $x_1$.

As in analyzing codegree $d$ products, we continue by analyzing canceling pairs of codegree $2d$ products that are products of top monomials of $v_i$, $u_i$ and $x_1$,  
with one codegree $2d$ monomials 
of these elements, such that this codegree $2d$ monomial of $u_1$ is not the codegree $d$ prefix of $v_1^0$, the codegree $2d$ monomial of $x_1$ is not the codegree $d$ 
prefix (or suffix)
of $w_0$, and the codegree $2d$ monomial of
$v_2$ is not the codegree $d$ suffix of $u_2^0$. These codegree $2d$ products are analyzed in the same way they were analyzed for smaller codegree products.

We are left with codegree $2d$ products that are either:
\roster


\item"{(1)}" a product of the top monomial of $u_1$ with a codegree $2d$ monomial of $x_1$, that is equal to a codegree $d$ prefix (or suffix) of the
top monomial of $w$, with
the top monomial of $u_2$.

\item"{(2)}" a product of the top monomial of $v_1$ with a codegree $2d$ monomial of $x_1$, that is equal to the codegree $d$ prefix of the top  monomial of $w$, with
the top monomial of $v_2$.

\item"{(3)}" 
 an odd set of codegree $2d$  products of the top monomial of $v_1$,  codegree $e_j$ monomial of $x_1$ and  codegree $p_j$ monomial of $v_2$, for
some positive pairs $(e_j,p_j)$ that satisfy $e_j+p_j=2d$, such that the codegree $e_j$ monomial of $x_1$ is the product of the top monomial of $s$ with a codegree
$e_j$ monomial of $w$, and the product of each codegree $e_j$ monomial of $x_1$ with a codegree $p_j$
monomial of $v_2$ is equal to the product of a codegree $d$ prefix of $w_0$ with $v_2^0$. 


\item"{(4)}" 
 an odd set of codegree $2d$ products of codegree $f_i$ monomials of $u_1$ and  codegree $g_i$ monomials of $x_1$ with the top monomial of $u_2$, for
some positive  pairs $(f_i,g_i)$ that satisfy $f_i+g_i=2d$, such that the codegree $g_i$ monomial of $x_1$ is the product of a codegree $g_i$ monomial of $w$
with the top monomial of $t$, and the product of each codegree $f_i$ monomial of $u_1$ with a codegree $g_i$
monomial of $x_1$ is equal to the product of $u_1^0$ with the codegree $d$ prefix of $w_0$. 

\endroster

Note that (1) exists if and only if (2) exists and they cancel each other.  
If (1) and (2) are the only existing possibilities, we add a codegree $2d$ monomial
to $w$, which is the codegree $d$ prefix or suffix of the top monomial of $w$. 

If only possibilities (3) and (4) exist, we add the codegree $2d$ prefix of $w_0$ to $w$. If (1)-(4) do all exist, $w$ remains unchanged.

This completes the analysis of codegree $2d$ products. 
codegree $2d+r$ products, for $1 \leq r < d$, are treated in the same way we treated codegree $d+r$ products. Codegree $3d$ products
are treated in the same we treated $2d$ products, and so on. Finally, in case $d=deg(u_1)-deg(v_1)=deg(v_2)-deg(u_2)$, we obtained the conclusion of the theorem.

\smallskip
Suppose that:
$d=min(deg(v_1),deg(u_2),deg(u_1)-deg(v_1))=deg(v_1)$. In that case we continue the analysis of codegree $r$ homogeneous parts in
$u_i,v_i$ and $x_1$, $d \leq r < deg(u_1)-deg(v_1)$, precisely as we analyzed the codegree $r$ homogeneous parts for $1 \leq r \leq d-1$.
For $r = deg(u_1)-deg(v_1)$, we use the same analysis that we apply for codegree $d$ products, in case $d=deg(u_1)-deg(v_1)$.
For $r > deg(u_1)-deg(v_1)$, we continue the analysis of codegree $r$ homogeneous parts according to the analysis of codegree higher than $d$ in case $d=deg(u_1)-deg(v_1)$.
The analysis in the case: 
$d=min(deg(v_1),deg(u_2),deg(u_1)-deg(v_1))=deg(u_2)$ is identical.

\line{\hss$\qed$}

Theorem 4.4 reduces the analysis of solutions to the equation: $u_1xu_2=v_1xv_2$ to the equation $xt=sx$, in case the equation
$u_1xu_2=v_1xv_2$ has a long enough solution, and the coefficients have no periodicity. The same techniques allows one to reduce 
a general equation with one variable, in case the coefficients have no periodicity.

\vglue 1.5pc
\proclaim{Theorem 4.7} Let $u_1,\ldots,u_n,v_1,\ldots,v_n \in  FA$ and suppose that the top homogeneous parts of $u_i$ and $v_i$ are monomials
with no periodicity, and that for at least one index $i$, $1 \leq i \leq n$, $u_i \neq v_i$.
Suppose that the equation:
$$u_1xu_2xu_3 \ldots u_{n-1}xu_n \, = \, v_1xv_2xv_3 \ldots v_{n-1}xv_n$$ has a solution $x_1$ of degree bigger than 
$2(deg(u_1)+ \ldots + deg(u_n))^2$. Suppose further, 
that all the periodicity
in the top monomials, that are associated with  the top monomials of the  two sides of the equation after substituting the solution $x_1$, 
is contained in the periodicity of the top
monomial of the solution $x_1$.

Then there exist some elements $s,t \in FA$, $deg(s)=deg(t) \leq max deg(u_i)$, such that:
\roster
\item"{(1)}" every solution of the equation $sx=xt$ is a solution of the given equation.

\item"{(2)}" every solution  $x_2$ of the given equation,  that satisfies: 
$$deg(x_2)>2(2+2^{deg(s)+2})+
(2(deg(u_1)+ \ldots + deg(u_n))^2$$ 
is a solution of the equation: $sx=xt$.
\endroster
\endproclaim

\nfp Let $x_1$ be a solution of the given equation that satisfies: $deg(x_1) > 2(deg(u_1)+ \ldots +deg(u_n))$. We start by looking at the top homogeneous part of $x_1$ that we
denote $x_1^0$. Clearly $x_1^0$ satisfies the homogeneous equation:
$$u_1^0x_1^0u_2^0x_1^0u_3^0 \ldots u_{n-1}^0x_1^0u_n^0 \, = \, v_1^0x_1^0v_2^0x_1^0v_3^0 \ldots v_{n-1}^0x_1^0v_n^0$$ 
where $u_i^0$ and $v_i^0$ are the top monomials in $u_i$ and $v_i$.

We start the analysis of the given equation under the assumption that there exists an index $i$, for which: $deg(u_i) \neq deg(v_i)$. In that case 
there is a shift between the appearances of some of the (homogeneous) elements $x_1^0$ in the two sides of the 
(homogeneous) equation. Let $i_1$ be the first index $i$, for which $deg(u_i) \neq deg(v_i)$. The next appearances of $x_1^0$ in the two sides of the equation must have a shift of
$|deg((u_{i_1})-deg(v_{i_1})|$. Since the top homogeneous parts of $u_i$ and $v_i$ are monomials, it follows that the top
homogeneous part of $x_1$ is a monomial as well. We denote it $x_1^0$.

Let $d$ be the minimum positive shift between pairs of appearances of $x_1^0$ in the top degree equation. Then $x_1^0=e_0(t_0)^b=(s_0)^bf_0$, for some elements $t_0,s_0$ in
the semigroup generated by the free generators of the algebra $FA$, $a_1,\ldots,a_k$. $deg(s_0)=deg(t_0)=d$, $e_0$ is a prefix of $s_0$ and suffix of $t_0$, and $f_0$
is a suffix of $t_0$ and prefix of $s_0$. Since the top monomial $u_i^0$ and $v_i^0$ have no periodicity, $t_0$ and $s_0$ have no periodicity as well.
 
\smallskip
In the sequel, we further assume 
that 
the degrees of the coefficients satisfy: $deg(u_i),deg(v_i)>d$ for
all the indices $i$.

Since we assumed that the length of $x_1^0$ is bigger than the sum of the lengths of the degrees, $deg(u_i)$, an appearance of $x_1^0$ in the product that
is associated with the top monomial in the left side of the equation overlaps with the corresponding appearance of $x_1^0$ in the right
side of the equation, and may overlap with the previous or the next appearance of $x_1^0$ of the right side of the equation as well. Our assumptions that
$deg(u_i),deg(v_i) > d$ together with the assumption that the coefficients have no periodicity imply that 
an appearance of $x_1^0$ in the product that is associated with top monomial in one side of the equation may overlap only with the corresponding 
appearance of $x_1^0$ in the other side of the equation (and not with the previous or the next one). 

Let $1 \leq i_1 < \ldots < i_{\ell} \leq n-1$ be the indices for which there is a (non-trivial) 
shift between the appearances of $x_1^0$ in the two sides of the equation. Let $1 \leq j_1 < \ldots < j_{n-1-\ell} \leq n-1$ be 
the complementary indices, i.e., those indices for which there is no shift between the corresponding appearances of the monomials
$x_1^0$ in the two sides of the equation.

We start by analyzing the codegree 1 monomials in the products that are associated with the two sides of the equation. We further assume that the length of the period in
$x_1^0$. i.e. $d=|deg(u_{i_1})-deg(v_{i_1})|>1$. Note that
any codegree 1 monomial in the two products is a product of top monomials with a single codegree 1 monomial from one of the $u_i$'s, $v_i$'s 
or one of the appearances of $x_1$.

We set $i_1$ to be the first index for which $|u_{i_1}-v_{i_1}| > 0$. Suppose that $i_1 > 1$. $x_1^0$ is periodic, its period is of length at least 2, and
$x_1^0$ contains at least 2 periods. Hence,  a codegree 1 product that contains a codegree 1 monomial of $u_1$ can cancel with either a codegree 1 product that contains
a codegree 1 monomial in $v_1$ or a codegree 1 product that contains a codegree 1 monomial in the first appearance of $x_1$. 

If the two canceling codegree 1 products contain
codegree monomials of $u_1$ and $v_1$, then these two codegree 1 monomials must be equal. Otherwise, the codegree 1 product that contains a codegree 1 monomial of
$u_1$ cancels with a codegree 1 product that contains a codegree 1 monomial of the first appearance of $x_1$. Now, this last codegree 1 monomial appears in the other side
of the equation as well, and it can cancel only with a codegree 1 product that contains a codegree 1 monomial of $v_1$, that must be identical to the
codegree 1 monomial of $u_1$ that we started with. Therefore, the codegree 1 homogeneous parts of $v_1$ and $u_1$ are equal.
Continuing with the same argument iteratively, the codegree 1 homogeneous parts of the elements $u_i$ and $v_i$ are equal for all $i<i_1$, and $i>i_{\ell}$.

Let $j_s$ be one of the indices for which there is no shift between the corresponding appearances of $x_1^0$ in the two sides of the equation. We look at the
codegree 1 products in the two sides of the equation. Each such codegree 1 product is a product of a single codegree 1 monomial from a single appearance of $x_1^0$ or
exactly one of the coefficients $u_i$ or $v_i$, with top degree monomials. Note that the codegree 1 products that involve codegree 1 monomials of the $j_s$ appearance of
$x_1$ in the two sides of the equation (and top degree monomials from all the coefficients and the other appearances of $x_1$), are precisely the same codegree 1 
products in the two sides of the equation. Hence, these do cancel. All the other codegree 1 products in the two sides of the equation contain
$x_1^0$ in the $j_s$ appearance of $x_1$. Since $x_1^0$ is periodic, and the length of the period is bigger than 1, a codegree 1 product that includes
a codegree 1 monomial to the left of the $j_s$ appearance of $x_1$ can not be equal to a codegree 1 product that contains a codegree 1 monomial to the right
of the $j_s$ appearance of $x_1$. 

Therefore, the left codegree 1 products (with respect to the $j_s$ appearance of $x_1$) 
from the two sides of the equation have to cancel and the right codegree 1
products have to cancel as well.
In particular, if for some index $i$, both the $i-1$ and the $i$ appearances of $x_1$ in the two sides of the equation have no shift, then $u_i$ and $v_i$
have the same codegree 1 homogeneous parts.

At this point we need to examine the appearances of the variables $x_1$ in which there is a shift between the two sides of the equation,
i.e., in place $i_1,\ldots,i_{\ell}$, and the coefficients, $u_i$ and $v_i$, that are connected to these appearances. To do that we break
the appearances of the variables $x_1$, and the coefficients $u_i$,$v_i$ in the two sides of the equations into regions (or intervals).

We look at the top monomial in the two sides of the equation.
For each index $i$ we add a breakpoint at the left point of the pair $u_i,v_i$, and to the right of that pair. We denote the variable that
is associated with the region (interval) between the
right point of the pair $u_{i},v_{i}$ and the left point of the pair $u_{i+1},v_{i+1}$, $w_i$. The top monomial of $w_i$ is a prefix or a suffix 
of the top monomial $x_1^0$ of $x_1$. We denote by $q_i$ the variable that is associated with the region between $w_{i-1}$ and $w_{i}$. Note
that the region that is associated with $q_i$ contains the support of $u_i$ and $v_i$. If the region that is associated with $q_i$ contains the
right part of the $i-1$ appearance of $x_1$, we denote the variable that is associated with that right part $t_{i-1}$. If it contains the
left part of the $i$ appearance of $x_1$, we denote the variable that is associated with that left part, $s_{i}$. 

As in our previous arguments, we intend to break the solution $x_1$, so that $x_1=s_iw_i=w_it_i$, whenever the variables $w_i,s_i,t_i$ are
defined and in appropriate abelian (quotient) groups. Furthermore, each of the elements $q_i$ can be broken according to the two sides of
the equation. Hence, we intend to show that $q_i=t_{i-1}u_is_i$ or $q_i=t_{i-1}u_i$ or $q_i=u_is_i$ or $q_i=u_i$, and correspondingly for the elements $v_i$ 
(instead of the $u_i$'s), depending on the way the elements $q_i$ are broken in the two sides of the equation.

Because of the periodicity of $x_1^0$, and since we assume that the length of the period of $x_1^0$ is bigger than $1$, a
codegree 1 product that contains a codegree 1 monomial of $v_i$ can not cancel with a codegree 1 product that contains a
codegree 1 monomial in $v_{i'}$ or $u_{i'}$ for $i \neq i'$, and likewise for the $u_i$'s.

Suppose that $q_i=v_i=t_{i-1}u_is_i$. In that case two codegree 1 products that contain codegree 1 monomials of the $i-1$ and $i$
appearances of $x_1$ that are both from the $v_i$ side, or both from the $u_i$ side, can not cancel. 
Furthermore, two codegree 1 products that cancel and belong to the two sides of the equation,
can not contain codegree 1 monomials from both appearances $i-1$ and $i$ of $x_1$.

Hence, in that case a pair of canceling codegree 1 products may either be:
\roster
\item"{(1)}" codimension 1 monomials of the same appearance of $x_1$ from the two sides of the equation.

\item"{(2)}" codimension 1 monomial of either the $i$ or $i-1$ appearance of $x_1$ for one product, and a codimension 1 monomial
of $u_i$ for the second product.

\item"{(3)}" codimension 1 monomial of either the $i$ or $i-1$ appearance of $x_1$ for one product, and a codimension 1 monomial
of $v_i$ for the second product.

\item"{(4)}" codimension 1 monomial of $u_i$ in one product, and codimension 1 monomial of $v_i$ in the second product.
\endroster

If case (1) occurs we add a codegree 1 monomial to $w_i$ or $w_{i-1}$ (depending on the appearance of $x_1$). In case (2) we add a codegree
1 monomial to $t_{i-1}$ or $s_i$, and the existing one to $u_i$. In case (3) we add a codegree 1 monomial to $w_i$ or $w_{i-1}$, the existing
codegree 1 monomial to $v_i$ and a codegree 1 monomial to $t_{i-1}$ or $s_i$ (depending on the appearance of $x_1$). In case (4) we add only the existing 
codegree 1 monomials to $v_i$ and $u_i$.

Suppose that  $q_i=t_{i-1}v_i=u_is_i$. Two codegree 1 products that contain codegree 1 monomials of the $i-1$ and $i$ appearances of $x_1$ from
the same side of the equation can not be equal. Furthermore, a codegree 1 product that contains a codegree 1 monomial of the $i-1$ appearance
of $x_1$ in the $u_i$ side can not cancel with a codegree 1 product that contains a codegree 1 monomial in the $i$ appearance of
$x_1$ from the $v_i$ side. Since we assumed that $deg(u_i), deg(v_i)>d$, and the coefficients have no periodicity, 
a codegree 1 product that contains a codegree 1 monomial of the $i-1$ appearance
of $x_1$ in the $v_i$ side can not cancel with a codegree 1 product that contains a codegree 1 monomial in the $i$ appearance of
$x_1$ from the $u_i$ side.

Like in the case $q_i=v_i=t_{i-1}u_is_i$, in that case a pair of canceling codegree 1 products may either be:
\roster
\item"{(1)}" codimension 1 monomials of the same appearance of $x_1$ from the two sides of the equation.

\item"{(2)}" codimension 1 monomial of either the $i$ or $i-1$ appearance of $x_1$ for one product, and a codimension 1 monomial
of $u_i$ for the second product.

\item"{(3)}" codimension 1 monomial of either the $i$ or $i-1$ appearance of $x_1$ for one product, and a codimension 1 monomial
of $v_i$ for the second product.

\item"{(4)}" codimension 1 monomial of $u_i$ in one product, and codimension 1 monomial of $v_i$ in the second product.
\endroster

If case (1) occurs we add a codegree 1 monomial to $w_i$ or $w_{i-1}$ (depending on the appearance of $x_1$). In case (2), if a codimension 1 of the $i-1$
appearance of $x_1$ is part of the canceling pair, we add a codimension 1 monomial to $w_{i-1}$, a codimension 1 monomial to $t_{i-1}$, and the existing
codimension 1 monomial to $u_i$. If a codimension 1 of the $i$ appearance of $x_1$ is part of the cancelling pair, we add codimension 1
monomial to $s_i$, and the existing codimension 1 monomial to $u_i$. In case (3) we do the equivalent additions for $v_i$, $w_i$, $t_{i-1}$ and $s_i$. 
In case (4) we just add the codimension 1 existing monomials to $u_i$ and $v_i$.

So far we have constructed elements $w_i$, $t_i$, and $s_i$, such that the equations: $x_1=w_it_i$, $x_1=s_iw_i$, $q_i=u_is_i$ or
$q_i=t_{i-1}u_i$, or $q_i=t_{i-1}u_is_i$ or $q_i=u_i$ (and correspondingly for the $v_i$'s) hold for products of codegree at most 1.
We continue by analyzing products of codegree $r$, $r<d$, assuming that we analyzed all the products of smaller codegree, constructed the
elements $w_i$, $s_i$ and $t_i$, and they satisfy the last equations for products of codegree at most $r-1$. 

\medskip
We analyze codegree $r$ products in a similar way to their analysis in the proof of theorem 4.4. 
First, note that if a codegree $r$ product is a product of monomials of $u_i$, $v_i$, $s_i$, $t_i$ and $w_i$,
that correspond to products of monomials of codegree smaller than $r$ of $u_i$, $v_i$, and all the appearances of $x_1$ from the two sides of the equation, 
then such codegree $r$ products cancel in pairs.

Let $i$ be an index for which $deg(u_i)=deg(v_i)$, and there is no shift between the $i-1$ and $i$ appearances of $x_1$. By our analysis of codegree 1 monomials,
the top monomials, and the codegree 1 homogeneous parts of $u_i$ and $v_i$ are identical. codegree $r$ products from one side of the equation, 
that contain codegree $r$ monomials of the
$i$ or $i-1$ appearances of $x_1$, cancel with corresponding codegree $r$ products from the other side of the equation. Hence, a codegree $r$ product that contains a 
codegree $r$ monomial of $u_i$ can cancel only with a codegree $r$ product that contains a codegree $r$ monomial of $v_i$. Therefore, the codegree $r$ homogeneous part
of $u_i$ is identical to the codegree $r$ homogeneous part of $v_i$. Furthermore,
for the purpose of analyzing codegree $r$ products, the given equation can be broken into finitely many equations by taking out such pairs of coefficients $u_i,v_i$,
and the appearances of the solution $x_1$ that are adjacent to them.

Suppose that  for some index $i$ there is no shift between the appearances of $x_1$ in the two sides of the equation. In that case codegree $r$ products that contain
codegree $r$ monomials of the $i$  appearance of $x_1$ from one side of the equation cancel with codegree $r$ products that contain codegree $r$ monomials of that $i$
appearance of $x_1$ from the other side of the equation. Hence, for the purpose of analyzing codegree $r$ products, the given equation breaks into several
equations, by taking out all the appearances of $x_1$ that have no shift. Therefore, for the continuation of the analysis of codegree $r$ products,
we may assume that there are no appearances of $x_1$ with no shift.

Since we assumed that the equation does not contain appearances of $x_1$ in the two sides of the equation with no shift between them, the
analysis of codegree $r$ products that contain positive codegree monomials of either $u_1$ or $v_1$, or positive codegree monomials of either
$u_2$ or $v_2$, is identical to the analysis of codegree $r$ monomials in theorem 4.4, i.e., in the equation: $u_1xu_2=v_1xv_2$. Hence, we only need to analyze
codegree $r$ products that contain positive codegree monomials from some element $q_i=u_is_i=t_{i-1}v_i$ or from an element $q_i=t_{i-1}u_is_i=v_i$.

Let $q_i=u_is_i=t_{i-1}v_i$. Since $u_i$ and $v_i$ have no periodicity, a codegree $r$ product that contains a codegree $r$ monomial of $x_1$ in its $i-1$ appearance
can not cancel with a codegree $r$ from the same side of the equation that contains a codegree $r$ monomial of $x_1$ in its $i$ appearance. Furthermore, a codegree
$r$ product that contains a codegree $r$ monomial of $x_1$ in its $i-1$ appearance from the $u_i$ side of the equation, can not cancel with a codegree $r$ product
that contains a codegree $r$ monomial of $x_1$ in its $i$ appearance from the $v_i$ side of the equation.  

Suppose that a codegree $r$ product, can be expressed as products of codegree $q_j$ monomials of $t_{i-1}$ and codegree $m_j$ monomials of $v_i$ with top monomials
of the other elements in the $v_i$ side of the equation, such that $q_j \geq 0$ and $m_j$ is positive  and $q_j+m_j=r$, in an odd number of ways. Such
codegree $r$ products can cancel with either:  
\roster
\item"{(1)}" 
an odd number of products of codegree $f_j$ monomials of $u_i$ and codegree $g_j$ monomials of $s_i$ with top monomials
of the other elements in the $u_i$ side of the equation, such that $g_j \geq 0$ and $f_j$ is positive  and $f_j+g_j=r$. 

\item"{(2)}" a product of a codegree $r$ monomial of $x_1$ in its $i-1$ appearance with other top monomials in the $v_i$ side of the equation.

\item"{(3)}" a product of a codegree $r$ monomial of $x_1$ in its $i$ appearance with other top monomials in the $v_i$ side of the equation.

\item"{(4)}" a product of a codegree $r$ monomial of $x_1$ in its $i-1$ appearance with other top monomials in the $u_i$ side of the equation.

\item"{(5)}" a product of a codegree $r$ monomial of $x_1$ in its $i$ appearance with other top monomials in the $u_i$ side of the equation.

\item"{(6)}" an odd number of  products of a codegree $b_j$ monomial of $w_{i-1}$ with a codegree $a_j$ monomial of $t_{i-1}$ for positive $a_j,b_j$, $a_j+b_j=r$, 
with top monomials
of the other elements from the $v_i$ side.

\item"{(7)}" an odd number of  products of a codegree $c_j$ monomial of $w_{i}$ with a codegree $h_j$ monomial of $s_{i}$ for positive $c_j,h_j$, $c_j+h_j=r$, with top monomials
of the other elements from the $u_i$ side.
\endroster

If only case (1) occurs we don't add anything to any of the elements except the existing codegree $r$ monomials of $u_i$ and $v_i$. If only case (2) occurs we add
a codegree $r$ monomial to $t_{i-1}$. If only case (3) occurs we add a codegree $r$ monomial to $w_i$ and a codegree $r$ monomial to $s_i$. If only case (4) occurs
we add a codegree $r$ monomial to $w_{i-1}$ and to $t_{i-1}$. If only case (5) occurs we add a codegree $r$ monomial to $s_i$. 

Cases (2) and (3) can not occur together, nor cases (4) and (5), nor cases (3) and (4).
If only cases (1), (2) and (4) occur, we add a codegree $r$ monomial to $w_{i-1}$. 
If only cases (1), (2) and (5) occur, we add codegree $r$ monomials to $t_{i-1}$
and $s_i$. 
If only (1), (3) and (5) occur, we add a codegree $r$ monomial to $w_i$. 

We still need to treat cases (6) and (7). Note that the existence of these cases means that  codegree $r$ products that were supposed to exist given the
smaller codegree monomials of the various elements, may or may not exist, depending on the existence of codegree $r$ monomials in the
various appearances of the element $x_1$. Also, note that case (6) can not occur with case (3), and case (7) can not occur with case (4). 

If only case (6) appears, we add a codegree $r$ monomial to $t_{i-1}$. If only case (7) appears we add a codegree $r$ monomial to $s_i$. If only cases (1), (2) and (6) appear,
we do not add anything. 
If only cases (1), (2) and (7) appear, we add a codegree $r$ monomial to $t_{i-1}$ and to $s_i$. 
If only (1), (3) and (7) appear, we do not add anything. 
If only (1), (4) and (6) appear, we add a codegree $r$ monomial to $w_{i-1}$. If only (1), (5) and (6) appear, we add  codegree $r$ monomials to $t_{i-1}$ and to $s_i$. If
only (1), (5) and (7) appear, we don't add anything.

If only (2), (4) and (6) appear, we add  codegree $r$ monomials to $w_{i-1}$ and to $t_{i-1}$.
If only (2), (5) and (6) appear, we add a codegree $r$ monomial to $s_i$.
If only (2), (5) and (7) appear, we add codegree $r$ monomial to $t_{i-1}$.
If only (3), (5) and (7) appear, we add  codegree $r$ monomials to $w_i$ and to $s_i$.  

If only (1), (6) and (7) appear, we add codegree $r$ monomials to $t_{i-1}$ and to $s_i$. If only (2), (6) and (7) appear, we add a codegree $r$ monomial to $s_i$.
If only (5), (6) and (7) appear, we add a codegree $r$ monomial to $t_{i-1}$.
If only (1), (2), (5), (6) and (7) appear, we do not add anything.

The case in which the codegree $r$ product is a product of case (1) is dealt with in a symmetric way. Hence, suppose that the codegree $r$ product is not a 
product of case (1) and can not be expressed in an odd number of ways as
products of codegree $q_j$ monomials of $t_{i-1}$ and codegree $m_j$ monomials of $v_i$ with top monomials
of the other elements in the $v_i$ side of the equation, such that $q_j \geq 0$ and $m_j$ is positive  and $q_j+m_j=r$.

If only (2) and (4) appear, we add a codegree $r$ monomial to $w_{i-1}$. If only (2) and (5) appear, we add codegree $r$ monomials to $t_{i-1}$ and $s_i$.
If only (2) and (6) appear, we do not add anything. If only (2) and (7) appear, we add codegree $r$ monomials to $t_{i-1}$ and $s_i$. If only (3) and (5) appear,
we add a codegree $r$ monomial to $w_{i}$. If only (3) and (7) appear, we add a codegree $r$ monomial to $w_i$. If only (4) and (6) appear, we
add a codegree $r$ monomial to $w_{i-1}$. If only  (5) and (6) appear, we add codegree $r$ monomials to $t_{i-1}$ and $s_i$. If only  (5) and (7)
appear, we do not add anything. If only (6) and (7) appear, we add codegree $r$ monomials to $t_{i-1}$ and $s_i$. 
Finally, if (2), (5), (6) and (7) appear, we do not add anything.

As in the proof of theorem 4.4, it can still  be that a codegree $r$ product is of type (7) and  can also be presented 
in an odd number of ways as  products of a codegree $b_j$ monomial of $w_{i}$ with a codegree $a_j$ monomial of $t_{i}$ for positive $a_j,b_j$, $a_j+b_j=r$, 
with top monomials
of the other elements from the $v_i$ side. 

In that case it can either be presented only in these two forms or also in both forms (3) and (5). If it can be presented in forms (3) and (5) we do not add
anything. If it can not we add a codegree $r$ monomial to $w_i$.

\smallskip
This concludes the construction of the elements $s_i,t_i,w_i$ for codegree $r$ products that involve $q_i=u_is_i=t_{i-1}v_i$.
Suppose that $q_i=v_i=t_{i-1}u_is_i$. 
Since $u_i$ and $v_i$ have no periodicity, a codegree $r$ product that contains a codegree $r$ monomial of $x_1$ in its $i-1$ appearance
can not cancel with a codegree $r$ product that contains a codegree $r$ monomial of $x_1$ in its $i$ appearance. 

Suppose that a codegree $r$ product, can be expressed as products of codegree $q_j$ monomials of $t_{i-1}$, codegree $m_j$ monomials of $u_i$, and codegree $p_j$ of
$s_i$ with top monomials
of the other elements in the $u_i$ side of the equation, such that $q_j,m_j,p_j \geq 0$,  either  $m_j>0$ or $q_j,p_j>0$, and $q_j+m_j+p_j=r$, in an odd number of ways. Such
codegree $r$ products can cancel with either:  
\roster
\item"{(1)}" a codegree $r$ monomial of $v_i$.

\item"{(2)}" a product of a codegree $r$ monomial of $x_1$ in its $i-1$ appearance with other top monomials in the $u_i$ side of the equation.

\item"{(3)}" a product of a codegree $r$ monomial of $x_1$ in its $i$ appearance with other top monomials in the $u_i$ side of the equation.

\item"{(4)}" a product of a codegree $r$ monomial of $x_1$ in its $i-1$ appearance with other top monomials in the $v_i$ side of the equation.

\item"{(5)}" a product of a codegree $r$ monomial of $x_1$ in its $i$ appearance with other top monomials in the $v_i$ side of the equation.

\item"{(6)}" an odd number of  products of a codegree $b_j$ monomial of $w_{i-1}$ with a codegree $a_j$ monomial of $t_{i-1}$ for positive $a_j,b_j$, $a_j+b_j=r$, 
with top monomials
of the other elements from the $u_i$ side.

\item"{(7)}" an odd number of  products of a codegree $c_j$ monomial of $w_{i}$ with a codegree $h_j$ monomial of $s_{i}$ for positive $c_j,h_j$, $c_j+h_j=r$, with top monomials
of the other elements from the $u_i$ side.
\endroster

According to the various cases, we add monomials to the variables $t_i,s_i,w_i$, in a similar way to what we did in case $q_i=u_is_i=t_{i-1}v_i$. 
If only case (1) occurs we don't add anything to any of the elements except the existing codegree $r$ monomials of $u_i$ and $v_i$. If only case (2) occurs we add
a codegree $r$ monomial to $t_{i-1}$. If only case (3) occurs we add a codegree $r$ monomial to $s_i$. If only case (4) occurs
we add a codegree $r$ monomial to $w_{i-1}$ and to $t_{i-1}$. If only case (5) occurs we add a codegree $r$ monomials to $w_i$ and to $s_i$. 

Cases (2) and (3) can not occur together, nor cases (4) and (5), nor cases (3) and (4), nor (2) and (5).
If only cases (1), (2) and (4) occur, we add a codegree $r$ monomial to $w_{i-1}$. 
If only (1), (3) and (5) occur, we add a codegree $r$ monomial to $w_i$. 

As in the case in which: $q_i=u_is_i=t_{i-1}v_i$, the existence of  cases (6) and (7) means that  codegree $r$ products that were supposed to exist given the
smaller codegree monomials of the various elements, may or may not exist, depending on the existence of codegree $r$ monomials in the
various appearances of the element $x_1$. Also, note that case (6) can not occur with cases (3) or (5), and case (7) can not occur with cases (2) or  (4). 

If only case (6) appears, we add a codegree $r$ monomial to $t_{i-1}$. If only case (7) appears we add a codegree $r$ monomial to $s_i$. If only cases (1), (2) and (6) appear,
we do not add anything. If only (1), (3) and (7) appear, we do not add anything. 
If only (1), (4) and (6) appear, we add a codegree $r$ monomial to $w_{i-1}$. If
only (1), (5) and (7) appear, we add a codegree $r$ monomial to $w_i$.

If only (2), (4) and (6) appear, we add  codegree $r$ monomials to $w_{i-1}$ and to $t_{i-1}$.
If only (3), (5) and (7) appear, we add  codegree $r$ monomials to $w_i$ and to $s_i$.  

The case in which  case (1) occurs is dealt with in an analogous way.
Hence, suppose that the codegree $r$ product is not a 
product of case (1) and can not be expressed in an odd number of ways as
products of codegree $q_j$ monomials of $t_{i-1}$ and codegree $m_j$ monomials of $u_i$ and codegree $p_j$ monomials of $s_i$ with top monomials
of the other elements in the $u_i$ side of the equation, such that $q_j,m_j,p_j \geq 0$, either  $m_j$ is positive or both $q_j,p_j$ are positive,  and $q_j+m_j+p_j=r$.

If only (2) and (4) appear, we add a codegree $r$ monomial to $w_{i-1}$. 
If only (2) and (6) appear, we do not add anything. If only (3) and (5) appear,
we add a codegree $r$ monomial to $w_{i}$. If only (3) and (7) appear, we do not add anything. If only (4) and (6) appear, we
add a codegree $r$ monomial to $w_{i-1}$. If only  (5) and (7)
appear, we add a codegree $r$ monomial to $w_{i}$.  

It can still  be that a codegree $r$ product is of type (7) and  can also be presented 
in an odd number of ways as  products of a codegree $b_j$ monomial of $w_{i}$ with a codegree $a_j$ monomial of $t_{i}$ for positive $a_j,b_j$, $a_j+b_j=r$, 
with top monomials
of the other elements from the $v_i$ side. 
We treat this case precisely as we treated it  in the case $q_i=u_is_i=t_{i-1}v_i$.

This concludes the construction of the elements $s_i,t_i,w_i$ for codegree $r$ products when $r < d$.
The elements $w_i$, $t_i$, and $s_i$ that we constructed so far satisfy the equations: $x_1=w_it_i$, $x_1=s_iw_i$, $q_i=u_is_i$ or
$q_i=t_{i-1}u_i$, or $q_i=t_{i-1}u_is_i$ or $q_i=u_i$ (and correspondingly for the $v_i$'s)  for products of codegree smaller than $d$.

\medskip
To continue we need to analyze products of codegree $d$ and higher. For presentation purposes we start this analysis under the additional 
assumption that all the appearances of $x_1$ in the two sides of the equation have non-trivial shifts, i.e., the appearances of the top
monomial of the solution $x_1^0$ in the two sides of the equality for the top monomials are shifted. This assumption enables us to analyze
the higher codegree products using the arguments that were used in the proof of theorem 4.4 and in analyzing smaller codegree products.
Afterwards we drop this assumption.
  
As in theorem 4.4, in analyzing codegree $d$ products, there are special codegree $d$ products that we need to single out and treat separately,
 as they may involve cancellations between codegree $d$ products that contain codegree $d$ monomials of $u_i$ or $v_i$  and
those that contain codegree $d$ monomials of $u_{i+1}$ or $v_{i+1}$.

As in analyzing smaller codegree products, note that codegree $d$ products that are products of smaller codegree monomials of
the $u_i$, $v_i$, $s_i$, $w_i$ and $t_i$, and correspond to smaller codegree monomials of $u_i$, $v_i$ and $x_1$ from the two
sides of the equation cancel in pairs.

In analyzing codegree $r$ products for $r<d$, there is no interaction between elements in $q_i$ and $q_j$ for $i \neq j$. As in the
proof of theorem 4.4, in analyzing
codegree $d$ products such interaction may happen if $i$ and $j$ are consecutive indices. Hence, in analyzing codegree $d$ products we need
to go over the various possibilities for $q_i$ and $q_{i+1}$.

Suppose that $q_i=u_is_i=t_{i-1}v_{i}$.
Suppose that a codegree $d$ product, can be expressed as products of codegree $q_j$ monomials of $t_{i-1}$ and codegree $m_j$ monomials of $v_i$ with top monomials
of the other elements in the $v_i$ side of the equation, such that $q_j \geq 0$ and $m_j$ is positive  and $q_j+m_j=d$, in an odd number of ways. If the
$q_i$ part of such a product
is not equal to $u_i^0$ nor to $v_i^0$ (the top monomials of $u_i$ and $v_i$), such codegree $r$ products are analyzed exactly 
in the same way they were analyzed in codegree $r$ products for $r<d$. 

$u_i^0 \neq v_i^0$ because we assumed that the top monomials of the coefficients have no periodicity. If the $q_i$ part of such a product equals to $v_i^0$,
the codegree $d$ product may be equal to a codegree $d$ product that contains positive codegree monomials in $q_{i-1}$.  
If the $q_i$ part of such a product equals to $u_i^0$,
the codegree $d$ product may be equal to a codegree $d$ product that contains positive codegree monomials in $q_{i+1}$.  

Suppose that the $q_i$ part of the codegree $d$ product equals $u_i^0$. Suppose further that $q_{i+1}=u_{i+1}s_{i+1}=t_iv_{i+1}$. In that case such a codegree $d$ product 
can cancel with codegree $d$ products that are either a subset of the  ones that were analyzed for products of smaller codegree, or products that include positive
codegree monomials of $q_{i+1}$:
\roster
\item"{(1)}" 
an odd number of products of codegree $f_j$ monomials of $u_i$ and codegree $g_j$ monomials of $s_i$ with top monomials
of the other elements in the $u_i$ side of the equation, such that $g_j \geq 0$ and $f_j$ is positive  and $f_j+g_j=d$. 


\item"{(2)}" a product of a codegree $d$ monomial of $x_1$ in its $i$ appearance with other top monomials in the $v_i$ side of the equation.


\item"{(3)}" a product of a codegree $d$ monomial of $x_1$ in its $i$ appearance with other top monomials in the $u_i$ side of the equation.


\item"{(4)}" an odd number of products of  $q_j$ monomials of $t_{i}$ and codegree $m_j$ monomials of $v_{i+1}$ with top monomials
of the other elements in the $v_i$ side of the equation, such that $q_j \geq 0$ and $m_j$ is positive  and the product of the monomial of $t_i$ with
the monomial of $v_{i+1}$ is $v_{i+1}^0$.

\item"{(5)}" 
an odd number of products of codegree $f_j$ monomials of $u_{i+1}$ and codegree $g_j$ monomials of $s_{i+1}$ with top monomials
of the other elements in the $u_i$ side of the equation, such that $g_j \geq 0$ and $f_j$ is positive  and  
the product of the monomial of $u_{i+1}$ with
the monomial of $s_{i+1}$ is $v_{i+1}^0$.

\item"{(6)}" an odd number of  products of a codegree $c_j$ monomial of $w_{i}$ with a codegree $h_j$ monomial of $s_{i}$ for positive $c_j,h_j$, $c_j+h_j=d$, with top monomials
of the other elements from the $u_i$ side.

\item"{(7)}" an odd number of  products of a codegree $b_j$ monomial of $w_{i}$ with a codegree $a_j$ monomial of $t_{i}$ for positive $a_j,b_j$, $a_j+b_j=d$, 
with top monomials
of the other elements from the $v_i$ side.
\endroster

Note that case (2) occurs if and only if case (3) occurs.
If only one of the cases (1) or (6) occurs, we treat them as they were treated in analyzing codegree $r$ products for $r<d$. If only case (4) or only
case (5)  occurs we add
1 (the identity) to $s_i$ and $t_i$, and the codegree $d$ prefix of $w_i^0$ to $w_i$. 
If only
case (7) occurs 
we add 1 to $s_i$ and the codegree $d$ prefix of $w_i^0$ to $w_i$.

If only cases (1)-(3) occur, or only cases (2), (3) and (6) occur, we treat them as in they were treated for codegree $r$ products, $r<d$. If only cases (2)
and (3) in addition to one of the cases cases (4) or (5) occur, we add 1 to $s_i$ and $t_i$. If only cases (2), (3) and (7) occur, we add 1 to $s_i$.
If only (1), (4) and (5) occur, we don't add anything. If only (1), (6) and one of (4) or (5) occur, we add 1 to $t_i$ and the codegree $d$ prefix of $w_i^0$ to $w_i$.
If only (1), (7) and one of (4) or (5)  occur, 
we add 1 to $t_i$. If only (1), (6) and (7) occur, we add the codegree $d$ prefix of $w_i^0$ to $w_i$. If only (6), (7)  and one of (4) or (5) occur, we add
1 to $s_i$ and $t_i$ and the codegree $d$ prefix of $w_i^0$ to $w_i$. If only (4), (5) and (6) occur, we add 1 to $s_i$. If only (4), (5) and (7) occur, we add
1 to $s_i$ and the codegree $d$ prefix of $w_i^0$ to $w_i$.

If only (1)-(5) occur, we add the codegree $d$ prefix of $w_i^0$ to $w_i$. If only (1)-(3) and (6)-(7) occur, we do not add anything. If only (1)-(3), (6) and
one of (4) or (5) occur, we add 1 to $t_i$.   
If only (1)-(3), (7) and
one of (4) or (5) occur, we add 1 to $t_i$ and the codegree $d$ prefix of $w_i^0$ to $w_i$. If only (1) and (4)-(7) occur, we add the codegree $d$ prefix of $w_i^0$ to
$w_i$.

If only (2)-(6) occur, we add 1 to $s_i$ and the prefix of codegree $d$ of $w_i^0$ to $w_i$. If only (2)-(5) and (7) occur, we add 1 to $s_i$. If only
(2)-(3), (6)-(7) and one of (4) or (5) occur, we add 1 to $s_i$ and $t_i$ and the codegree $d$ prefix of $w_i^0$ to $w_i$. If all the possibilities (1)-(7) occur,
we do not add anything.    

Suppose that a codegree $d$ product can   be expressed as a product in case (1), and can not be expressed as products of codegree $q_j$ monomials of $t_{i-1}$ 
and codegree $m_j$ monomials of $v_i$ with top monomials
of the other elements in the $v_i$ side of the equation, such that $q_j \geq 0$ and $m_j$ is positive  and the $q_i$ part of the product is $u_i^0$ in an
even number (possibly none) ways. In that case the analysis of such a product and the monomials that are added to the elements $t_i$, $s_i$ and $w_i$ are
analogous to the analysis described above. 

Suppose that  such a codegree $d$ product can not be expressed as a product in case (1), but it can be expressed as a product in case (6). If only (6) and (7) occur,
we add the codegree $d$ prefix of $w_i^0$ to $w_i$. If only (6) and one of (4) or (5) occur, we add 1 to $t_i$ and the prefix of codegree $d$ of $w_i^0$ to $w_i$. If
only (4)-(7) occur, we add the codegree $d$ prefix of $w_i^0$ to $w_i$. If only (2)-(3) and (6)-(7) occur, we do not add anything. If only (2)-(3), (6) and one of
(4) or (5) occur, we add 1 to $t_i$. If only (2)-(7) occur, we do not add anything.

This concludes the analysis of such codegree $d$ products in case $q_{i+1}=u_{i+1}s_{i+1}=t_iv_{i+1}$. Suppose that $q_{i+1}=u_{i+1}=t_iv_{i+1}s_{i+1}$.
As before, such a codegree $d$ product 
can cancel with codegree $d$ products that are either a subset of the  ones that were analyzed for products of smaller codegree, or products that include positive
codegree monomials of $q_{i+1}$:
\roster
\item"{(1)}" 
an odd number of products of codegree $f_j$ monomials of $u_i$ and codegree $g_j$ monomials of $s_i$ with top monomials
of the other elements in the $u_i$ side of the equation, such that $g_j \geq 0$ and $f_j$ is positive  and $f_j+g_j=d$. 


\item"{(2)}" a product of a codegree $d$ monomial of $x_1$ in its $i$ appearance with other top monomials in the $v_i$ side of the equation.


\item"{(3)}" a product of a codegree $d$ monomial of $x_1$ in its $i$ appearance with other top monomials in the $u_i$ side of the equation.


\item"{(4)}" an odd number of products of codegree $q_j$ monomials of $t_{i}$, codegree $m_j$ monomials of $v_{i+1}$, and codegree $p_j$ monomials of
$s_{i+1}$ with top monomials
of the other elements in the $v_i$ side of the equation, such that $q_j,m_j,p_j \geq 0$,  either  $m_j>0$ or $q_j,p_j>0$, and $q_j+m_j+p_j=d$, and the product
of the corresponding monomials of $t_i$, $v_{i+1}$ and $s_{i+1}$, is the codegree $d$ suffix of $u_{i+1}^0$.

\item"{(5)}" a product a monomial of $u_{i+1}$, which is the codegree $d$ suffix of $u_{i+1}^0$,
 with the top monomials of the all the other elements from the $u_i$ side of the equation.

\item"{(6)}" an odd number of  products of a codegree $c_j$ monomial of $w_{i}$ with a codegree $h_j$ monomial of $s_{i}$ for positive $c_j,h_j$, $c_j+h_j=d$, with top monomials
of the other elements from the $u_i$ side.

\item"{(7)}" an odd number of  products of a codegree $b_j$ monomial of $w_{i}$ with a codegree $a_j$ monomial of $t_{i}$ for positive $a_j,b_j$, $a_j+b_j=d$, 
with top monomials
of the other elements from the $v_i$ side.
\endroster

Analyzing the various possibilities in this case is identical to the case $q_{i+1}=u_{i+1}s_{i+1}=t_iv_{i+1}$. 

Recall that we assumed that:  $q_i=u_is_i=t_{i-1}v_{i}$, and suppose
that a codegree $d$ product, can be expressed as products of codegree $f_j$ monomials of $u_{i}$ and codegree $g_j$ monomials of $s_i$ with top monomials
of the other elements in the $u_i$ side of the equation, such that $f_j \geq 0$ and $g_j$ is positive  and $f_j+g_j=d$, in an odd number of ways, and such that
the product of the monomial of $u_i$ with the monomial of $s_i$ is $v_i^0$.
In that case,
the codegree $d$ product may be equal to a codegree $d$ product that contains positive codegree monomials in $q_{i-1}$.  
Such a codegree $d$ product 
can cancel with codegree $d$ products that are either a subset of the  ones that were analyzed for products of smaller codegree, or products that include positive
codegree monomials of $q_{i-1}$:
\roster
\item"{(1)}" 
an odd number of products of codegree $q_j$ monomials of $t_{i-1}$ and codegree $m_j$ monomials of $v_i$ with top monomials
of the other elements in the $v_i$ side of the equation, such that $q_j \geq 0$ and $m_j$ is positive  and $q_j+m_j=d$. 


\item"{(2)}" a product of a codegree $d$ monomial of $x_1$ in its $i-1$ appearance with other top monomials in the $u_i$ side of the equation.


\item"{(3)}" a product of a codegree $d$ monomial of $x_1$ in its $i-1$ appearance with other top monomials in the $v_i$ side of the equation.


\item"{(4)}" 
an odd number of products of codegree $f_j$ monomials of $u_{i-1}$ and codegree $g_j$ monomials of $s_{i-1}$ with top monomials
of the other elements in the $u_i$ side of the equation, such that $g_j \geq 0$ and $f_j$ is positive  and  
the product of the monomial of $u_{i-1}$ with
the monomial of $s_{i-1}$ is $u_{i-1}^0$.

\item"{(5)}" an odd number of products of  $q_j$ monomials of $t_{i-2}$ and codegree $m_j$ monomials of $v_{i-1}$ with top monomials
of the other elements in the $v_i$ side of the equation, such that $q_j \geq 0$ and $m_j$ is positive  and the product of the monomial of $t_{i-2}$ with
the monomial of $v_{i-1}$ is $u_{i-1}^0$.

\item"{(6)}" an odd number of  products of a codegree $c_j$ monomial of $w_{i-1}$ with a codegree $h_j$ monomial of $t_{i-1}$ for positive $c_j,h_j$, $c_j+h_j=d$, with top monomials
of the other elements from the $v_i$ side.

\item"{(7)}" an odd number of  products of a codegree $a_j$ monomial of $s_{i-1}$ with a codegree $b_j$ monomial of $w_{i-1}$ for positive $a_j,b_j$, $a_j+b_j=d$, 
with top monomials
of the other elements from the $u_i$ side.
\endroster

The analysis of this case is identical to the case in which the $q_i$ part of a codegree $d$ product is $v_i^0$, and there is a possible cancellation with
codegree $d$ products that contain positive codegree monomials of $q_{i+1}$. An identical analysis applies also when $q_{i-1}=v_{i-1}=t_{i-2}u_{i-1}s_{i-1}$.

Suppose that $q_i=u_i=t_{i-1}v_is_i$ and $q_{i+1}=v_{i+1}=t_iu_{i+1}s_{i+1}$. Suppose that a codegree $d$ product can be presented in an odd number of ways as 
products of codegree $q_j$ monomials of $t_{i-1}$, codegree $m_j$ monomials of $v_{i}$, and codegree $p_j$ monomials of
$s_{i}$ with top monomials
of the other elements in the $v_i$ side of the equation, such that $q_j,m_j,p_j \geq 0$,  either  $m_j>0$ or $q_j,p_j>0$, and $q_j+m_j+p_j=d$, and the product
of the corresponding monomials of $t_{i-1}$, $v_{i}$ and $s_{i}$, is the codegree $d$ prefix of $u_{i}^0$. 

Such a codegree $d$ product 
can cancel with codegree $d$ products that are either a subset of the  ones that were analyzed for products of smaller codegree, or products that include positive
codegree monomials of $q_{i+1}$:
\roster
\item"{(1)}" a product of a monomial of $u_i$, which is the codegree $d$ prefix of $u_i^0$, 
with the top monomials of all the other elements from the $u_i$ side of the equation.


\item"{(2)}" a product of a codegree $d$ monomial of $x_1$ in its $i$ appearance with other top monomials in the $v_i$ side of the equation.


\item"{(3)}" a product of a codegree $d$ monomial of $x_1$ in its $i$ appearance with other top monomials in the $u_i$ side of the equation.


\item"{(4)}" a product of a monomial of $v_{i+1}$, which is the codegree $d$ suffix of $v_{i+1}^0$, with the top monomials of the all the other elements from the $v_i$ 
side of the equation.

\item"{(5)}" an odd number of products of codegree $q_j$ monomials of $t_{i}$, codegree $m_j$ monomials of $u_{i+1}$, and codegree $p_j$ monomials of
$s_{i+1}$ with top monomials
of the other elements in the $u_i$ side of the equation, such that $q_j,m_j,p_j \geq 0$,  either  $m_j>0$ or $q_j,p_j>0$, and $q_j+m_j+p_j=d$, and the product
of the corresponding monomials of $t_i$, $v_{i+1}$ and $s_{i+1}$, is the codegree $d$ suffix of $v_{i+1}^0$.

\item"{(6)}" an odd number of  products of a codegree $c_j$ monomial of $w_{i}$ with a codegree $h_j$ monomial of $s_{i}$ for positive $c_j,h_j$, $c_j+h_j=d$, with top monomials
of the other elements from the $v_i$ side.

\item"{(7)}" an odd number of  products of a codegree $b_j$ monomial of $w_{i}$ with a codegree $a_j$ monomial of $t_{i}$ for positive $a_j,b_j$, $a_j+b_j=d$, 
with top monomials
of the other elements from the $u_i$ side.
\endroster

Analyzing the various possibilities in this case is identical to the case $q_{i}=t_{i-1}v_i=u_is_i$. The analysis of the remaining case, in which:
$q_i=u_i=t_{i-1}v_is_i$ and $q_{i-1}=v_{i-1}=t_{i-2}u_{i-1}s_{i-1}$ is identical to the previous cases as well.

This concludes the construction of the elements $s_i,t_i,w_i$ for codegree $r$ products when $r \leq  d$, in case all the pairs of appearances of the top monomial of
the solution $x_1$ in the two sides of the equation, have non-trivial shifts. 
The elements $w_i$, $t_i$, and $s_i$ that we constructed so far satisfy the equations: $x_1=w_it_i$, $x_1=s_iw_i$, $q_i=u_is_i$ or
$q_i=t_{i-1}u_i$, or $q_i=t_{i-1}u_is_i$ or $q_i=u_i$ (and correspondingly for the $v_i$'s)  for products of codegree smaller or equal to $d$.

\medskip
As in the proof of theorem  4.4, we continue with the analysis of codegree $d+r$ products for $r<d$. First, as in analyzing smaller codegree products, 
codegree $d+r$ products that are products of
smaller codegree monomials of $u_i$, $v_i$, $s_i$, $t_i$, and $w_i$, that correspond to products of smaller codegree monomials of $u_i$, $v_i$ and $x_1$ (in all
its appearances) from the two sides of the equation, 
cancel in pairs. We start with two lemmas that are the analogues of lemma 4.5 and 4.6.

\vglue 1.5pc
\proclaim{Lemma 4.8} 
Suppose that a codegree $d+r$ product can be presented both as: 
\roster
\item"{(1)}" a product of a codegree $h$ monomial of $s_i$ with a codegree $c$ monomial of $w_{i}$,  for positive $c,h$, $c+h=d+r$, with top monomials
of the other elements from the $u_i$ side.

\item"{(2)}" a product of a codegree $b$ monomial of $w_{i}$ with a codegree $a$ monomial of $t_{i}$ for positive $a,b$, $a+b=d+r$, 
with top monomials
of the other elements from the $v_i$ side.
\endroster

Such a codegree $d+r$ product may only be presented as a product of smaller codegree monomials or (only) in one of the following two products:  
\roster
\item"{(i)}" a product of a codegree $d+r$ monomial of $x_1$ in its $i$ appearance with other top monomials in the $v_i$ side of the equation.

\item"{(ii)}" a product of a codegree $d+r$ monomial of $x_1$ in its $i$ appearance with other top monomials in the $u_i$ side of the equation.
\endroster
\endproclaim

\nfp In case it can be presented as another product of a codegree $d+r$ monomial with top degree monomials,  either the top monomial of $s_i$ or the top monomial of $t_i$
overlap with themselves with a cyclic shift, hence, they must be periodic. A contradiction to the assumption that the coefficients
do not have non-trivial periodicity. 

\line{\hss$\qed$}

\vglue 1.5pc
\proclaim{Lemma 4.9} 
With the notation of lemma 4.8, if a codegree $d+r$ product
can be presented in an odd number of ways as a product in the form (1) and in an even or none ways as a product of form (2), then such a product can be presented
precisely in one of the forms (i) or (ii). If a codegree $d+r$ product can be presented precisely in one of the forms (i) or (ii), then it can be presented
precisely in one of the forms (1) or (2) in an odd number of ways.

If a codegree $d+r$ product can be presented in odd number of ways in both forms (1) and (2), then it can either be presented in both forms (i) and (ii) or in neither of them.
If a codegree $d+r$ product can be presented in both forms (i) and (ii) then it can either be presented in both forms (1) or (2) in an odd number of ways,
or in both of them in even or no ways.
\endproclaim

\nfp If a codegree $d+r$ product can be presented in both forms (1) and (2) (odd or even number of times), the conclusion follows from lemma 4.8. Suppose that
it can be presented in an odd number of ways in form (1) and none in form (2). If it can also be presented as a codegree $d+r$ product that involves a positive 
codegree monomials of $u_{j}, v_{j}, s_j, t_j$ or $x_j$, for $j>i$, the top monomial of $u_{i+1}$ must have non-trivial periodicity, a contradiction. If it
can be also presented as a codegree $d+r$ product from the $u_i$'s sides of the equation that involves monomials of positive
codegree monomials of $u_{j}, s_j, t_j$ or $x_j$, for $j<i$, the top monomial of $u_{i}$ must have non-trivial periodicity, a contradiction.  

Suppose that the given codegree $d+r$ product can be also presented as a product of either:
\roster
\item"{(1)}" a codegree $q$ of $t_{i-1}$ and a codegree $m$ of $v_i$ with other top monomials from the $v_i$'s side of the equation.

\item"{(2)}" a codegree $f$ of $u_i$ and a codegree $g$ of $s_i$ with other top monomials from the $u_i$'s side of the equation.

\item"{(3)}" a codegree $d+r$ product from the $v_i$'s side of the equation that involves monomials of positive 
codegree monomials of $v_{j}, s_j, t_j$ or $x_j$  for $j<i$. 
\endroster
In all these cases the suffix of length $r$ of the top monomial of $u_i$ is identical to  the prefix of length $r$ of the period of $x$. If $r \leq deg(v_i)-d$, then
$v_i$ has non-trivial periodicity, a contradiction. Otherwise, the top monomial in the two sides of the equation contains periodicity that is not part of the periodicity of 
the solution $x$, a contradiction to our assumptions.

\line{\hss$\qed$}

Suppose that $q_i=u_is_i=t_{i-1}v_i$, and let $r$ be an integer, $0 < r < d$. By lemma 4.9 if a codegree $d+r$ product can be presented in an odd number of ways 
in the form (1) of lemma 4.8 then 
either:
\roster
\item"{(1)}"  it can be also presented   in an odd number of ways as in form (2) of lemma 4.8 and either in both forms (i) and (ii) in lemma 4.8, or in neither of
them.

\item"{(2)}"  it can be presented in an even or no ways in form (2) of lemma 4.8, and it can also be presented precisely in one of the forms   (i) or (ii) in lemma 4.8.
\endroster

By lemma 4.9, if a codegree $d+r$ product can be presented in form (i) of lemma 4.8, and  in even or no ways in forms (1) or (2) of that lemma, 
then it can also be presented in form (ii) of lemma 4.8.

Hence, if a codegree $d+r$ product can be presented in an odd number of ways in one of the forms (1),(2),(i) or (ii), then the appearances 
of the codegree $d+r$ products in these forms cancel in pairs. If it appears in odd number of ways in forms (1) and (2), and in forms (i) and (ii), we do not add anything.
If it appears in odd number of ways in forms (1) and (2) and not in the forms (i) nor (ii), we add a codegree $d+r$ monomial to $w_i$. If it appears in an odd number of ways
in the form (1), in an even or no ways in the form (2),  and appears in the form (i) we add a codegree $d+r$ monomial to $w_i$. 
If it appears in an odd number of ways in the form (1), in an even or no ways in the form (2), and in the form (ii), we do not add anything. If it appears in an even or no
ways in the forms (1) and (2), and in both form (i) and (ii), we add a codegree $d+r$ monomial to $w_i$.

Therefore, if a codegree $d+r$ product can be presented in an odd number of ways as 
products of codegree $q_j$ monomials of $t_{i-1}$ and codegree $m_j$ monomials of $v_i$ with top monomials
of the other elements in the $v_i$ side of the equation, such that $q_j \geq 0$ and $m_j$ is positive  and $q_j+m_j=d+r$, then it must be presented in 
an odd number of ways as
products of codegree $f_j$ monomials of $u_i$ and codegree $g_j$ monomials of $s_i$ with top monomials
of the other elements in the $u_i$ side of the equation, such that $g_j \geq 0$ and $f_j$ is positive  and $f_j+g_j=r$. 

This concludes the construction of the elements $s_i,t_i,w_i$ in case $q_i=t_{i-1}v_i=u_is_i$ (note that the elements $s_i,t_i$ did not change), to ensure that the
equalities they are suppose to satisfy hold for products up to codegree $d+r$. 

Suppose that $q_i=u_i=t_{i-1}v_is_i$. Lemmas 4.8 and 4.9 and their proofs remain valid in this case. Hence,  
a codegree $d+r$ product can be expressed in an odd number of ways as products of codegree $q_j$ monomials of $t_{i-1}$, codegree $m_j$ monomials of $v_i$, and codegree $p_j$ of
$s_i$ with top monomials
of the other elements in the $v_i$ side of the equation, such that $q_j,m_j,p_j \geq 0$,  either  $m_j>0$ or $q_j,p_j>0$, and $q_j+m_j+p_j=d+r$, if and only if
it is equal to a codegree $d+r$ monomial of $u_i$.

This concludes our treatment of codegree $d+r$ products for $r<d$. We continue by analyzing codegree $2d$ products. Lemmas 4.8 and 4.9 remain valid for codegree $2d$ products.
Hence, the analysis of codegree $2d$ products is identical to the analysis of codegree $d+r$ products for $r<d$. The analysis of higher codegree products, for codegree up to twice
the  maximal degree of the elements $u_i,v_i$ is identical as well. 

Hence, 
in case $deg(u_i),deg(v_i)>d$ and all the appearances
of the elements $x_1$ in the two sides of the equation have non-trivial shifts, 
we finally constructed elements $s_i,t_i,w_i$ that satisfy the equations:
\roster
\item"{(i)}" $q_i=u_is_i=t_{i-1}v_{i}$ or $q_i=v_i=t_{i-1}u_is_i$ or with exchanging the appearances of $u_i$ and $v_i$ and the equations.

\item"{(ii)}" $x_1=s_iw_i=w_it_{i}$ in the abelian group:
$G^{deg(x_1)}/G^{deg(x_1)-2(deg(s_i))}$. 
\endroster

Therefore, $s_1$ and $t_{n-1}$ are uniquely defined, and (given $x_1$) $w_1$ and $w_{n-1}$ are uniquely defined in the abelian group:
$G^{deg(w_1)}/G^{deg(w_1)-2(deg(s_i))}$. Hence, $t_1$ and $s_2$ are uniquely defined, and $w_2$ is uniquely defined in the abelian group: 
$G^{deg(w_2)}/G^{deg(w_2)-2(deg(s_i))}$. Continuing iteratively, all the elements $s_i$, $t_i$ are uniquely defined, and the elements $w_i$ are uniquely defined in  
the abelian group:
$G^{deg(w_i)}/G^{deg(w_i)-2(deg(s_i))}$.   

Since $s_iw_i=w_it_i$ it follows that: $s_ix_1=x_1t_i$ in the abelian group:
$G^{deg(s_ix_1)}/G^{deg(s_ix_1)-2(deg(s_i))}$. This implies that for every pair $i,j$ , $1 \leq i,j \leq n$, $(s_i+s_j)x_1=x_1(t_i+t_j)$ in the abelian group:
$G^{deg(s_ix_1)}/G^{deg(s_ix_1)-2(deg(s_i))}$,
 so for every pair $i,j$
either $s_i=s_j$ and $t_i=t_j$ or $s_i=s_j+1$ and $t_i=t_j+1$.

Since every pair $(s_i,t_i)$ is either $(s_1,t_1)$ or $(s_1+1,t_1+1)$, it follows that every element $\hat x$, that satisfies $s_1 \hat x= \hat xt_1$, is a solution of the
given equation. It remains to prove that every long enough solution of the given equation is a solution of the equation: $s_1x=xt_1$.

Let $x_2$ be a solution of the given equation that satisfies: 
$$deg(x_2)>2(2+2^{deg(s)+2})+
(2(deg(u_1)+ \ldots + deg(u_n))^2.$$ 
By continuing the analysis of higher codegree monomials of the solution $x_2$, we get that there exist elements $w_i$, such that for every index $i$, $1 \leq i \leq n$,
 there exists an element $w_i$, that satisfies: $s_iw_i=w_it_i=x_2$ in the abelian group:
$G^{deg(x_2)}/G^{deg(s_1)-1}$. By the argument that was used to prove lemma 4.2, it follows that there exists a solution $\hat x$ to the equation: $s_1x=xt_1$.  

$x_2$ satisfies $s_1x_2=x_2t_1$ in the abelian group:  
$G^{deg(s_1x_2)}/G^{2deg(s_1)-1}$. Hence, there exists an element $\hat x_2$, which is a solution of the equation $s_1x=xt_1$, and $x_2+\hat x_2=r$, where
$deg(r) \leq   
2+2^{deg(s_1)+2}$.

Suppose that the given equation is $v_1xv_2xv_3=u_1xu_2xu_3$, where $deg(v_1)<deg(u_1)$, and $deg(v_2)=deg(u_2)$. In this case, $u_1=v_1s_1$, $t_1u_2=v_1s_2$ and
$v_3=t_2u_3$. Hence, $(\hat x_2+r)v_2(\hat x_2+r)t_2=s_1(\hat x_2+r)u_2(\hat x_2+r)$. $\hat x_2$ is a solution to the equation $s_1x=xt_1$, so it is a solution
to the given equation. Therefore:
$$\hat x_2v_2rt_2+
rv_2\hat x_2t_2=
s_1\hat x_2u_2r+ 
s_1ru_2\hat x_2$$ 
in the abelian group: 
$G^{deg(
\hat x_2v_2rt_2)}/
G^{deg(
r v_2rt_2)}$.

Hence:
$$\hat x_2(v_2rt_2+
t_1u_2r)= 
(rv_2s_2+ 
s_1ru_2)\hat x_2$$ 
in the same abelian group. This implies that: 
$v_2rt_2+
t_1u_2r=p(t_1)$ and  
$rv_2s_2+ 
s_1ru_2=p(s_1)$ for the same polynomial $p$, in the abelian group: 
$G^{deg(
v_2rt_2)}/
G^{deg(
v_2rt_2)+deg(r)-deg(x_2)}$.

$t_1u_2=v_2s_2$, so $v_2(rt_2+s_2r)=p(t_1)$ in the abelian group: 
$G^{deg(
v_2rt_2)}/
G^{deg(
v_2rt_2)+deg(r)-deg(x_2)}$. By our assumption on $deg(x_2)$ it follows that: $v_2(rt_2+s_2r)=p(t_1)$. Similarly, $(rt_1+s_1r)u_2=p(s_1)$.
Hence, $p(s_1)$ is either $0$, or its leading term is of degree at least 2.

Since $(s_1,t_1)$ equals $(s_2,t_2)$ or $(s_2+1,t_2+1)$, we get that: $v_2(rt_1+s_1r)u_2=v_2p(t_1)=p(s_1)u_2$. We look at the leading term in the two sides of the
last equality. Since we assumed that the top monomials of $u_2$ and $v_2$ are not periodic, the top monomial of $u_2$ must be $\beta s_0$, and the top monomial
of $v_2$ must be $t_0 \beta$ where $\beta$ is a prefix of $t_0$ and a suffix of $s_0$. Hence, $t_0=\beta \alpha$ and $s_0= \alpha \beta$. But this is a contradiction, since
we assumed that the periodicity in the top monomials in the two sides of the given equation is contained in the solution $x_2$. Therefore, $s_1r+rt_1=0$, so $r$ is a solution of
the equation $s_1x=xt_1$, which means that: $x_2=\hat x_2+r$ is a solution to $s_1x=xt_1$ as well.

If the equation is: $u_1xu_2xu_3=v_1xv_2xv_3$, and $deg(v_1)>deg(u_1)$, $deg(v_3)>deg(u_3)$, then by the same arguments we get that $r$ (the remainder) has
to satisfy the equation: 
$$(rt_1+s_1r)v_2s_2 \hat x_2= \hat x_2t_1v_2(rt_2+s_2r).$$
That implies that if $rt_1+s_1r \neq 0$, $u_2$ must contain periodicity, a contradiction to our assumptions. Therefore, $rt_1+s_1r=0$, and both $r$ and $x_2$
are solutions of the equations: $s_1x=xt_1$.

Suppose that the length of the equation is bigger. $x_2$ is a long solution, and $x_2=\hat x_2+r$, where $\hat x_2$ is a solution of the equation: $s_1x=xt_1$, and
$deg(r) \leq 
2+2^{deg(s)+2}$. In that case we get the equality:
$$(\hat x_2+r)v_2(\hat x_2+r)v_3 \ldots v_{n-1}(\hat x_2+r)t_{n-1}=
s_1(\hat x_2+r)u_2(\hat x_2+r)u_3 \ldots u_{n-1}(\hat x_2+r).$$ 
and since $\hat x_2$ is a solution of the equation: $s_1x=xt_1$, we get the equality:
$$rv_2\hat x_2v_3 \ldots v_{n-1}\hat x_2t_{n-1}+ \ldots +
\hat x_2 v_2\hat x_2v_3 \ldots \hat x_2v_{n-1}rt_{n-1}=$$
$$s_1ru_2\hat x_2u_3 \ldots u_{n-1}\hat x_2+ \ldots +
s_1\hat x_2 u_2\hat x_2u_3 \ldots \hat x_2u_{n-1}r$$
in the abelian group:
$G^{m_1}/G^{m_2}$, where $m_1=
deg(s_1\hat x_2 u_2\hat x_2u_3 \ldots \hat x_2u_{n-1}r)$, and $m_2=m_1-deg(\hat x_2)+deg(r)$.

That implies the equality:
$$(s_1r+rt_1)u_2\hat x_2u_3 \ldots u_{n-1}\hat x_2+ 
\hat x_2 v_2\hat x_2v_3 \ldots \hat x_2v_{n-1}(rt_{n-1}+s_{n-1}r)+
\hat x_2w\hat x_2x_2=0$$
for some element  $w \in FA$, and in the same abelian group: $G^{m_1}/G^{m_2}$.
 
Therefore, 
$(s_1r+rt_1)u_2\hat x_2u_3 \ldots u_{n-1}$ is a polynomial in $s_1$, and:
$v_2\hat x_2v_3 \ldots \hat x_2v_{n-1}(rt_{n-1}+s_{n-1}r)$ is a polynomial in $t_1$. If $s_1r+rt_1 \neq 0$ this implies that the top monomials in the two sides of
the equation has periodicity that is not contained in the top monomials of the appearances of the solution $x_2$, a contradiction to our assumptions.
Hence, $r$ is a solution of the equation: $s_1x=xt_1$, and so is $x_2$.



\medskip
This concludes the proof of theorem 4.7 in case all the appearances of the top monomial of a solution $x_1$ in the  two monomials, that are the top products  in the
two sides of the given equation, have non-trivial shifts. We still
 need to complete the proof in the cases in which there are appearances of the top monomial of a solution $x_1$ with zero shifts. 

\vglue 1.5pc
\proclaim{Lemma 4.10} Let $u_1,u_2,v_1,v_2 \in FA$ satisfy $u_1 \neq v_1$, $deg(u_i)=deg(v_i)$, $i=1,2$, and suppose that the top homogeneous parts of $u_i$ and $v_i$ are monomials
(for $i=1,2)$) with no non-trivial periodicity. Then, if there exists a solution $x_1$ to the equation $u_1xu_2=v_1xv_2$, and $deg(x_1)> 2(deg(u_1)+deg(u_2))$, then there
exist elements $s,t \in FA$, such that $x$ is a solution of the equation $u_1xu_2=v_1xv_2$ if and only if it is a solution of the equation: $sx=xt$.
\endproclaim

\nfp  The top monomials of $u_1$ and $v_1$, and of $u_2$ and $v_2$, have to be equal. We set: $u_1=v_1+\mu_1$, $v_2=u_2+\mu_2$. $deg(\mu_1) < deg(v_1)$ and $deg(\mu_2) < deg(u_2)$.  
Hence, $(v_1+\mu_1)xu_2=v_1x(u_2+\mu_2)$, that implies: $\mu_1xu_2=v_1x\mu_2$. Since the top homogeneous parts of $v_1$ and $u_2$ are monomials
with no periodicity, so are the top homogeneous parts of $\mu_1$ and $\mu_2$. Since $deg(\mu_1) < deg(v_1)$ and $deg(\mu_2) < deg (v_2)$, the conclusion of the
lemma follows from theorem 4.4.  

\line{\hss$\qed$}

\vglue 1.5pc
\proclaim{Proposition 4.11} Let $u_1,u_2,u_3,v_1,v_2,v_3 \in FA$ satisfy $u_1 \neq v_1$, $u_3 \neq v_3$, $deg(u_i)=deg(v_i)$, $i=1,2,3$, 
and suppose that the top homogeneous parts of $u_i$ and $v_i$ 
are monomials
(for $i=1,2,3$) with no non-trivial periodicity. Then, if there exists a solution $x_1$ to the equation $u_1xu_2xu_3=v_1xv_2xv_3$,
and the only non-trivial periodicity in the top monomials of the two sides of the equation is contained in the top monomial of the solution $x_1$,  
and $deg(x_1)> 2(deg(u_1)+deg(u_2)+deg(u_3))$, then there
exist elements $s,t \in FA$, such that up to a swap between the $u$'s and the $v$'s:
\roster
\item"{(1)}" there exists $\mu_1$ for which $u_1=\mu_1 (s+1)$ and $v_1=\mu_1 s$.

\item"{(2)}" there exists $\mu_2$ and $\tau_2$, for which $t \mu_2=\tau_2 s$. Furthermore, $u_2=\tau_2 (s+1)$ and $v_2= (t+1) \mu_2$.

\item"{(3)}" there exists $\mu_3$ for which $u_3=t \mu_3$ and $v_3=(t+1) \mu_3$.
\endroster

As in the conclusion of theorem 4.7, every solution of the equation $sx=xt$ is a solution of the given equation $u_1xu_2xu_3=v_1xv_2xv_3$. 
Every solution  $x_2$ of the given equation, $u_1xu_2xu_3=v_1xv_2xv_3$, for which: 
$deg(x_2)>2(2+2^{deg(s_1)+2}+deg(u_1)+deg(u_2)+deg(u_3))$ is a solution of the equation $sx=xt$.
\endproclaim

\nfp  Since the top homogeneous parts of the $u_i$'s and $v_i$'s are monomials, the top homogeneous part of a solution $x_1$ to the given equation:
$u_1xu_2xu_3=v_1xv_2xv_3$ has to be a monomial as well. Also, the top monomial in $u_i$ equals the top monomial in $v_i$, $i=1,2,3$.

We look at the highest degree for which for some index $i$, $u_i \neq v_i$. This can not occur for a single index $i$. If $u_2=v_2$ at that highest degree, then the top
monomial in $u_2$ (and $v_2$), must have periodicity, a contradiction to our assumptions. Let $d$ be the codegree of that degree, and 
suppose that up to this codegree, $u_3=v_3$. In that case, the equation for codegree $d$ products reduces to the equation $u_1xu_2=v_1xv_2$. If we set $u_1=v_1+\mu_1$ and 
$v_2=u_2+\mu_2$, then for the codegree $d$ products, we get the equation: $\mu_1xu_2=v_1x\mu_2$. This implies that the top part of $\mu_1$ and $\mu_2$ are monomials,
that are the codegree $d$ prefix and suffix of the top monomials of $v_1$ and $v_2$ in correspondence, and that the top monomial of $x$ has a period of length $d$.  

In that case, it must be that $u_3=v_3$ for all the homogeneous parts of codegree less than $2d$, and hence, $\mu_1xu_2=v_1x\mu_2$ for all the products up to codegree $d$.
Therefore, there exists an element $s$, and an element $t$, such that $v_1=u_1=\mu_1s$ in the abelian group: $G^{deg(u_1)}/G^{deg(u_1)-d}$,  and $v_2=u_2=t \mu_2$ in the abelian  
group: $G^{deg(u_2)}/G^{deg(u_2)-d}$.   

Since $u_3=v_3$ for all the homogeneous parts of codegree less than $2d$, and the top monomial of $u_3$ (and $v_3$) do not have non-trivial periodicity, it follows that
$u_3=v_3$. Hence, $\mu_1xu_2=v_1x\mu_2$, and the conclusion of the proposition follows from theorem 4.4 in this case.

Suppose that for the codegree $d$ homogeneous parts: $u_i \neq v_i$ for $i=1,2,3$. In that case, we get the equation $(v_1+\mu_1)xu_2xu_3=v_1xv_2x(v_3+\mu_3)$, and $u_i=v_i$, $i=1,2,3$,
for all the homogeneous parts of codegree smaller than $d$. Hence, the top homogeneous parts of $\mu_1$ and $\mu_3$ are monomials, which are the codegree $d$ prefix and
suffix of the top monomials of $u_1$ and $u_3$ in correspondence. The top monomial of $x_1$ (the given solution to the given equation)
 has to be periodic, with a period of length $d$. Furthermore, $v_2=b_2+\mu_2$
and $u_2=b_2+\tau_2$ in the abelian group: 
$G^{deg(u_2)}/G^{deg(u_2)-(d+1)}$, where the top homogeneous parts of $\mu_2$ and $\tau_2$ are the codegree $d$ prefix and suffix of the top monomial of $u_2$ (and $v_2$).   

We continue by looking at products of codegree $d+1$. Every such product that contains monomials in $u_i$ that appear also in $v_i$, for $i=1,2,3$, cancels with
a similar product from the other side of the equation. Hence, to analyze cancellations,
we need to consider codegree $d+1$ products that contain monomials from $\mu_1$ or $\mu_3$, or monomials
of codegree $d$ and $d+1$ of $u_2$ and $v_2$ that do not appear in both. 

Suppose that a codegree $d+1$ product contains a codegree $d+1$ monomial of $\mu_1$, i.e., a monomial in $u_1$ that is not in $v_1$. Such a codegree $d+1$ product must contain
the top monomial of $x_1$ in its two appearances, and the top monomial of $u_2$ and $u_3$. Since the top monomial of $v_1$ doesn't have
non-trivial periodicity, such a codegree $d+1$ product can not cancel with a codegree $d+1$ product
that contains the top monomial of $v_1$. Therefore, a codegree $d+1$ product that cancels with it must contain a codegree 1 monomial of either $u_1$ or $v_1$, or the top monomial of
$\mu_1$. Since the top monomial of $v_2$ contains no periodicity, if this codegree $d+1$ product  contains a codegree 1 monomial of $u_1$ or $v_1$ it must contain 
the top monomial of $\mu_2$. Hence, this codegree $d+1$ product has to be from the $v_i$ side of the equation, and the codegree $d+1$ monomial of $\mu_1$ is the codegree
$d$ prefix of a codegree 1 monomial in $v_1$, times the (prefix) period of the top monomial of $x_1$, which is the degree $d$ suffix of $v_1$.
If such a codegree $d+1$ product cancels with a codegree $d+1$ product that contains the top monomial of $\mu_1$, then it must contain a codegree 1 monomial of $x_1$.

By the techniques that we used in the proofs of theorem 4.4 and in the first part of theorem 4.7, there exists an element $s_1$, $deg(s_1)=d$, with a top monomial $\mu_1$,
such that $\mu_1 s_1=u_1=v_1$, in the abelian group: 
$G^{deg(u_1)}/G^{deg(u_1)-2}$.    

Suppose that a codegree $d+1$ product contains a codegree $d+1$ monomial of $u_2$ or $v_2$. Since the top monomial of $u_2$ (and $v_2$) contains no
periodicity, such a product can cancel only with either:
\roster
\item"{(1)}" a codegree $d+1$ product that contains a codegree 1
monomial of $u_2$ or $v_2$ and the top monomial of either $\mu_1$ or $\mu_3$.

\item"{(2)}"  a codegree $d+1$ product that contains the top monomial of $\mu_2$, and a codegree 1 monomial in the second appearance of $x_1$, and
the top monomial of $\mu_1$, the top monomial of $u_2$, and the same codegree 1 monomial in the second appearance of $x_1$.   

\item"{(3)}"  a codegree $d+1$ product that contains the top monomial of $\tau_2$, and a codegree 1 monomial in the first appearance of $x_1$, and
the top monomial of $\mu_3$, the top monomial of $v_2$, and the same codegree 1 monomial in the first appearance of $x_1$.   
\endroster

Note that the two products that appear in possibilities (2) and (3) cancel each other. Hence, a codegree $d+1$ product that contains a codegree $d+1$ product that appears in
$u_2$ or $v_2$, but not both, must cancel with a unique codegree $d+1$ product that is described in (1).
  
Suppose that a codegree $d+1$ product, contains the top monomial of $\mu_1$ and a codegree 1 monomial of $u_2$. Since the top monomial of $u_2$ (and $v_2$)
has no periodicity, it can cancel only with a codegree $d+1$ product that contains either:
\roster
\item"{(1)}" a codegree 1 monomial of $v_2$ and the top monomial of $\mu_3$.

\item"{(2)}" a codegree $d+1$ monomial of $u_2$ or $v_2$.

\item"{(3)}" a codegree 1 monomial of the first appearance of $x_1$ and the top monomial of $\mu_2$.
\endroster

Similarly, suppose that a codegree $d+1$ product, contains the top monomial of $\mu_3$ and a codegree 1 monomial of $v_2$. 
It can cancel only with a codegree $d+1$ product that contains either:
\roster
\item"{(1)}" a codegree 1 monomial of $u_2$ and the top monomial of $\mu_1$.

\item"{(2)}" a codegree $d+1$ monomial of $u_2$ or $v_2$.

\item"{(3)}" a codegree 1 monomial of the second appearance of $x_1$ and the top monomial of $\tau_2$.
\endroster
Furthermore, a codegree $d+1$ product that contains the top monomial of $\mu_2$ can not cancel with a codegree $d+1$ product that contains the top monomial of $\tau_2$.

Hence, we can look at the collection of codegree $d+1$ products that contain the top monomial of $\mu_1$ and the entire collection of codegree 1 monomials of
$u_2$. Each such product cancels with precisely one  product that contains either a codegree $d+1$ monomial of $u_2$ or $v_2$, or a codegree 1 monomial of the first
appearance of $x_1$ and the top monomial of $\mu_2$, or a codegree 1 monomial of $v_2$ and the top monomial of $\mu_3$. A similar statement holds for
codegree $d+1$ products that contain a codegree 1 monomial of $v_2$ and the top monomial of $\mu_3$.

Therefore, there exist elements $t_1,s_2,b,w_1,w_2,\tau_2,\mu_2$ such that:
\roster
\item"{(1)}" $t_1 \mu_2 =v_2$ and $\tau_2 s_2=v_2$ in the abelian group:
$G^{deg(u_2)}/G^{deg(u_2)-2}$. $deg(s_2)=deg(t_1)=d$.

\item"{(2)}" $b+\tau_2=u_2$ and $b+\mu_2=v_2$ in the abelian group:
$G^{deg(u_2)}/G^{deg(u_2)-(d+2)}$. 

\item"{(3)}" $x_1=s_1w_1=w_1t_1=s_2w_2=w_2t_2$ in the abelian group:
$G^{deg(x_1)}/G^{deg(x_1)-2}$. 
\endroster

\smallskip
We continue by induction for $1 \leq r \leq d$, and assume that for $r<d$
there exist elements $t_1,s_2,b,w_1,w_2,\tau_2,\mu_2$ such that the equalities that were true for the top 2 monomials and codegree $d$ and codegree $d+1$ monomials
hold for the top $r$ monomials, and for codegree $d$-$d+r-1$ monomials:
\roster
\item"{(1)}" $t_1 \mu_2 =v_2$ and $\tau_2 s_2=v_2$ in the abelian group:
$G^{deg(u_2)}/G^{deg(u_2)-r}$. $deg(s_2)=deg(t_1)=d$.

\item"{(2)}" $b+\tau_2=u_2$ and $b+\mu_2=v_2$ in the abelian group:
$G^{deg(u_2)}/G^{deg(u_2)-(d+r)}$. 

\item"{(3)}" $x_1=s_1w_1=w_1t_1=s_2w_2=w_2t_2$ in the abelian group:
$G^{deg(x_1)}/G^{deg(x_1)-r}$. 
\endroster

\smallskip
We continue by studying codegree $d+r$ products. All such products that involve only monomials of codegree less than $d$ of the $u_i$, $v_i$, $1 \leq i \leq 3$,
cancel
in pairs. 
All such products that involve only monomials of codegree less than $d+r$ of the $u_i$, $v_i$, $1 \leq i \leq 3$, and codegree less than $r$ of $x_1$ (in its two
appearances from both sides of the equation) cancel
in pairs by the induction hypothesis.

Hence, to analyze the structure of $u_1$ and $v_1$ (and hence, of $\mu_1$ and $s_1$) we only need to consider codegree $d+r$  products that contain:
\roster
\item"{(i)}"  a codegree $d+r$ monomial of $u_1$  that does not appear in $v_1$ and vice versa. 

\item"{(ii)}" a codegree $r$ monomial of $v_1$ and the top monomial of $\mu_2$. 

\item"{(iii)}" a codegree $d+q_j$ monomial of $\mu_1$, $q_j <r$, and a codegree $r-q_j$ monomial of the first appearance of $x_1$.

\item"{(iv)}" a codegree $p_j$ monomial of $v_1$, $p_j < r$, and a codegree $r-p_j$ monomial of the first appearance of $x_1$ and
the top monomial of $\mu_2$.
\endroster

A product of type (iv)  that cancels with products of type (i) or (ii) must cancel with a corresponding product of type (iii) by our induction hypothesis. A product of type
(iii) that cancels with a product of type (i) or (ii) and in which $q_j$ is positive, and the codegree $r-q_j$ monomial of the first appearance of $x_1$ is obtained as 
a product of a codegree  $r-m_j$ monomial of $s_1$ with
a codegree $m_j-q_j$ monomial of $w_1$, for $q_j<m_j<r$, cancels with a product of type (iv). 

Therefore, to analyze the structure of $u_1$, $v_1$, $s_1$ and $w_1$, we consider only those codegree $d+r$ products that can be presented  either in 
form (i) or (ii), that we denote (1) and (2) in the sequel, or in the form:
\roster
\item"{(3)}" a product of the top monomial of $\mu_1$, and a codegree $r$ monomial of the first appearance of $x_1$.
\endroster

A codegree $d+r$ product that can be presented in one of the forms (1)-(3) can cancel with either:
\roster
\item"{(4)}" an odd number of products of a codegree $d+q_j$ monomial of $\mu_1$, $0<q_j <r$, and a codegree $r-q_j$ monomial of the first appearance of $x_1$.

\item"{(5)}" an odd number of products of a codegree $d+q_j$ monomial of $\mu_1$, $0<q_j <r$, and a product of a codegree $r-q_j$ monomial of $s_1$ with
the top monomial of $w_1$.

\item"{(6)}" an odd number of products of a codegree $d+q_j$ monomial of $\mu_1$, $q_j <r$, and a product of a codegree $r-m_j$ monomial of $s_1$ with
a codegree $m_j-q_j$ monomial of $w_1$, where $q_j<m_j<r$.

\item"{(7)}" an odd number of products of a  codegree $p_j$ monomial of $v_1$, $0<p_j < r$, and a codegree $r-p_j$ monomial of the first appearance of $x_1$ and
the top monomial of $\mu_2$.

\item"{(8)}" a product of the top monomial of $v_1$, a codegree $r$ monomial of the first appearance of $x_1$ and
the top monomial of $\mu_2$.

\item"{(9)}" an odd number of products of a  codegree $p_j$ monomial of $v_1$, $0<p_j < r$, and a codegree $m_j$ monomial of the first appearance of $x_1$, 
$0<m_j$, $p_j+m_j < r$ and
a codegree $d+r-p_j-m_j$ monomial  of $\mu_2$.

\item"{(10)}" an odd number of products of a codegree $d+q_j$ monomial of $\mu_1$, a codegree $m_j$ monomial of the first appearance of $x_1$, $0<m_j$, $q_j+m_j<r$,
and a codegree $r-m_j-q_j$ monomial of $u_2$.
\endroster

If (1) or (2) occur, (8)  can not occur, and  (6) occurs if and only if (7) occurs as well. If (1) occurs, (3) can not occur. 
Suppose that (1) occurs. If in addition only (2) occurs, we add a codegree $d+r$ monomial to $\mu_1$.  
If in addition to (1) only (4) and (5) occur, we also add a codegree $d+r$ monomial to $\mu_1$. If in addition to (1) 
only (5), (6) and (7) occur, we add a codegree $d+r$ monomial to $\mu_1$.
If (1) occurs, (9) and (11) can not occur. 


Suppose that (2) occurs. If in addition only (3) occurs (and in addition possibly (4), (6) and (7)) we add a codegree $r$ monomial to $s_1$.
If in addition to (2) only (4) and (5) occur, we do not add anything. If in addition to (2) 
only (5), (6) and (7) occur, we do the same.
If (2) occurs, (8)-(10) can not occur. 

Suppose that (3) occurs. The codegree $r$ monomial of $x_1$ can not be presented both as a product of the top monomial of $s_1$ with a codegree $r$ monomial
of $w_1$, and as a codegree $r$ monomial of $s_1$ with the top monomial of $w_1$. We look at all the possible ways to present the codegree $r$ monomial of
$x_1$ as a product of a codegree $q_j$ monomial of $s_1$ with a codegree $r-q_j$ monomial of $w_1$, for $0<q_j<r$. If the number of such products is odd we don't add anything. If the
number is even,
we either add a codegree $r$ monomial to $s_1$ or a codegree $r$ monomial to $w_1$ (but not both). The validity of this addition of codegree $r$ monomial to either $s_1$ or
$w_1$ can be verified by going over the possible cancellation of the given codegree $d+r$ product with all the other possible forms of such a product.  

\smallskip
This concludes the adaptation of $s_1$, $\mu_1$, and $w_1$ to include codegree $r$ monomials. The same adaptation works for $t_2$, $\mu_3$ and $w_2$.
It is still left to analyze $u_2$ and $v_2$ in order to add codegree $r$ monomials to $\mu_2$ and $\tau_2$, such that the equalities that by induction hold
for the top codegree $r-1$ parts of these elements will hold for the top codegree $r$ part.

To analyze the structure of $u_2$ and $v_2$ (and hence, of $\mu_2$, $\tau_2$, $t_1$  and $s_2$) we start by observing the following:
\roster
\item"{(i)}"  the  codegree $d+r$ products that contain either a positive codegree  monomial of $u_1$ or a positive codegree monomial of the first appearance of $x_1$,
a codegree $d+q_j$ monomial of $\tau_2$, $q_j<r$, a  monomial of the second appearance of $x_1$, 
and a monomial of $u_3$, cancel with codegree $d+r$ products that contain
either a positive codegree  monomial of $v_1$  or a  positive codegree monomial of the first appearance of $x_1$,
a codegree $p_j$ monomial of $v_2$, $p_j<r$, a  monomial of the second appearance of $x_1$,  
and a monomial of $\mu_3$. 

\item"{(ii)}"  the  codegree $d+r$ products that contain a monomial of $v_1$, a  monomial of the first appearance of $x_1$,
a codegree $d+q_j$ monomial of $\mu_2$, $q_j<r$, and either a positive codegree monomial of the second appearance of $x_1$, 
or a positive codegree monomial of $v_3$, cancel with codegree $d+r$ products that contain
a monomial of $\mu_1$, a   monomial of the first appearance of $x_1$,
a codegree $p_j$ monomial of $u_2$, $p_j<r$, either a positive codegree monomial of the second appearance of $x_1$,  
or a positive codegree monomial of $u_3$. 
\endroster

Hence, to analyze the structure of $u_2$ and $v_2$  we only need to consider codegree $d+r$  products that contain:
\roster
\item"{(i)}"  a codegree $d+r$ monomial of $u_2$  or of $v_2$. 

\item"{(ii)}" the top monomial of $\mu_1$ and a codegree $r$ monomial of $u_2$ or
a codegree $r$ monomial of $v_2$ and the top monomial of $\mu_3$. 

\item"{(iii)}" a codegree $d+q_j$ monomial of $\tau_2$, $q_j <r$, and a codegree $r-q_j$ monomial of the second appearance of $x_1$ or
a codegree $r-q_j$ of the first appearance of $x_1$ and a codegree $d+q_j$ monomial of $\mu_2$, $q_j <r$.

\item"{(iv)}" the top monomial of $\mu_1$, a codegree $r-p_j$ monomial of the first appearance of $x_1$  and a codegree $p_j$ monomial of $u_2$, $p_j < r$, or
a codegree $p_j$ monomial of $v_2$, $p_j < r$ and
and a codegree $r-p_j$ monomial of the second appearance of $x_1$  and the top monomial of $\mu_3$.
\endroster

If there are two products of codegree $d+r$ of type (i), they cancel each other, and we can ignore them in analyzing codegree $d+r$ products. 
Therefore, to analyze the structure of $u_2$, $v_2$, $s_2$, $t_1$, $\mu_2$ and $\tau_2$, we consider only those codegree $d+r$ products that can be presented  either in 
form (i) or (ii), that we denote (1) and (2) in the sequel, or in codegree $d+r$ products in the form:
\roster
\item"{(3)}" a product of the top monomial of $\tau_2$, and a codegree $r$ monomial of the second appearance of $x_1$.

\item"{(4)}" a codegree $r$ monomial of the first appearance of $x_1$, and the top monomial of $\mu_2$.

\item"{(5)}" an odd number of products of a codegree $d+q_j$ monomial of $\tau_2$, $0<q_j <r$, and a codegree $r-q_j$ monomial of the second appearance of $x_1$.

\item"{(6)}" an odd number of products of a codegree $d+q_j$ monomial of $\tau_2$, $0<q_j <r$, and a product of a codegree $r-q_j$ monomial of $s_2$ with
the top monomial of $w_2$.

\item"{(7)}" an odd number of products of a codegree $d+q_j$ monomial of $\tau_2$, $q_j <r$, and a product of a codegree $r-m_j$ monomial of $s_2$ with
a codegree $m_j-q_j$ monomial of $w_2$, where $q_j<m_j<r$.

\item"{(8)}" an odd number of products of a  codegree $p_j$ monomial of $v_2$, $0<p_j < r$, and a codegree $r-p_j$ monomial of the second appearance of $x_1$ and
the top monomial of $\mu_3$.

\item"{(9)}" a product of the top monomial of $v_2$, a codegree $r$ monomial of the second appearance of $x_1$ and
the top monomial of $\mu_3$.
\endroster

And similarly from the other sides of the equation:
\roster
\item"{(10)}" an odd number of products of a codegree $r-q_j$ monomial of the first appearance of $x_1$, and a codegree $d+q_j$ monomial of $\mu_2$, $0<q_j <r$. 

\item"{(11)}" an odd number of products of the top monomial of $w_1$,  a codegree $r-q_j$ monomial of $t_1$, 
 and a codegree $d+q_j$ monomial of $\mu_2$, $0<q_j <r$. 

\item"{(12)}" an odd number of products of a  
a codegree $m_j-q_j$ monomial of $w_1$, 
a  codegree $r-m_j$ monomial of $t_1$,
a codegree $d+q_j$ monomial of $\mu_2$, $q_j <r$,
$q_j<m_j<r$.

\item"{(13)}" an odd number of products of the top monomial of $\mu_1$, a  codegree $r-p_j$ monomial of the first appearance of $x_1$, and
a  codegree $p_j$ monomial of $u_2$, $0<p_j < r$. 

\item"{(14)}" a product of  the top monomial of $\mu_1$, a codegree $r$ monomial of the first appearance of $x_1$, and
the top monomial of $u_2$. 
\endroster

Suppose that (1) occurs. If only one of the possibilities in (2) occurs, we add a codegree $d+r$ monomial to $\mu_2$ or $\tau_2$, depending
which of the two possibilities in (2) occurs. If (1) occurs, (3) and (4)  can not occur.
If in addition to (1) only (5) occurs, then (6) or (7) must occur and not both. If only (5) and (6) occur, we add a codegree $d+r$ monomial to $\tau_2$. If in addition
to (1),   (5) and (7) occur, then
(8) must occur as well, and hence at least an additional possibility must occur. If in addition to (1), (8) occurs, then (5) and (7) must occur as well, so
an additional possibility must occur. If (1) occurs, (9) can not occur. The possibilities (10)-(14) are parallel to (5)-(9) and are dealt accordingly. 

Suppose that (1) and the two possibilities in (2) occur. If in addition only (5) and (6) occur, we add a codegree $d+r$ monomial only to $\mu_2$, and if only (10) and (11)
occur, we add a codegree $d+r$ monomial to $\tau_2$. Suppose that  (1) and only one of the products in the form (2) occur, wlog the product from the $v_i$ side, i.e.,
the one that contains $\mu_3$.
If in addition (5), (6), (10) and
(11)  occur, we add a codegree $d+r$ monomial to $\mu_2$. 

Suppose that one of the possibilities (2) occurs, wlog the one from the $v_i$ side. If the only additional product that cancels with it is also a product in form (2) from
the $u_i$ side of the equation, we add a codegree $d+r$ monomial to both $\tau_2$ and $\mu_2$. If in addition to the form (2) only possibility (3) occurs,  we add a codegree
$r$ monomial to $s_2$. (4) can not occur. If only (5) and (6) occur, we do not add anything. If  (5) and (7) occur, (8) must occur as well. (9) can not occur. If in
addition (10) and
(11) occur, we add a codegree $d+r$ monomial to both $\tau_2$ and $\mu_2$. If only (3), (5), (6), (10) and (11) occur, we add a codegree $r$ monomial to $s_2$, and a
codegree $d+r$ monomial to both $\mu_2$ and $\tau_2$.    

Suppose that the two possibilities in part (2) occur. In that case (3) can not occur. If in addition (5), (6), (11) and (12) occur, we do not add anything.
Suppose that (3) occurs. In that case (4) can not occur. If in addition only (5) and (6) occur, we add a codegree $r$ monomial to $s_2$. If in addition to
(3) only (9) occurs, we add a codegree $r$ monomial to $w_2$. a If (3) occurs, then (10)-(14)
can not occur. If (4) occurs the analysis is analogous to the case in which (3) occurs.

Suppose that (5) and (6) occur. In that case (9) can not occur. If (10) and (11) occur as well, we add a codegree $d+r$ monomial to both $\tau_2$ and $\mu_2$.

\smallskip
This concludes our treatment of codegree $d+r$ products for $r< d$. So far we proved that:
\roster
\item"{(1)}" $\mu_1 s_1 =u_1=v_1$ in the abelian group:
$G^{deg(u_1)}/G^{deg(u_1)-d}$. $deg(s_1)=d$. $u_1=v_1+\mu_1$ in the abelian group:
$G^{deg(u_1)}/G^{deg(u_1)-2d}$. 

\item"{(2)}" $t_1 \mu_2 =v_2$ and $\tau_2 s_2=v_2$ in the abelian group:
$G^{deg(u_2)}/G^{deg(u_2)-d}$. $deg(s_2)=deg(t_1)=d$. $deg(\mu_2)=deg(\tau_2)=deg(u_2)-d$.

\item"{(3)}" $b_2+\tau_2=u_2$ and $b_2+\mu_2=v_2$ in the abelian group:
$G^{deg(u_2)}/G^{deg(u_2)-(2d)}$. 

\item"{(4)}" $x_1=s_1w_1=w_1t_1=s_2w_2=w_2t_2$ in the abelian group:
$G^{deg(x_1)}/G^{deg(x_1)-d}$. 
\endroster

We continue by analyzing codegree $2d$ products. The analysis of codegree $2d$ products is similar to 
the analysis of codegree $d+r$ products for $r<d$. In their analysis we use the following observations:
\roster
\item"{(i)}" all the codegree $2d$ products that contain monomials of codegree smaller than $d$ from the elements: $u_i$, $v_i$, and $x$ in its two
appearances, cancel in pairs.

\item"{(ii)}" all the codegree $2d$ products that contain a monomial of codegree bigger than $d$, from $b_1$, $b_2$ or $b_3$, cancel in pairs.
\endroster

Hence, we need to analyze only those codegree $2d$ products that contain monomials from either $\mu_1, \mu_2, \tau_2, \mu_3$, or monomials of codegree 
$d$ from $b_1,b_2,b_3$. To analyze the elements: $u_1,v_1,b_1,\mu_1,s_1$ and $w_1$, we need to analyze codegree $2d$ products that
contain one of the following:  
\roster
\item"{(i)}"  a codegree $2d$ monomial of $u_1$  that does not appear in $v_1$ and vice versa. 

\item"{(ii)}" a codegree $d$ monomial of $v_1$ and the top monomial of $\mu_2$. 

\item"{(iii)}" a codegree $d+q_j$ monomial of $\mu_1$, $q_j <d$, and a codegree $d-q_j$ monomial of the first appearance of $x_1$.

\item"{(iv)}" a codegree $p_j$ monomial of $v_1$, $p_j < d$, and a codegree $d-p_j$ monomial of the first appearance of $x_1$ and
the top monomial of $\mu_2$.

\item"{(v)}" note that the codegree $2d$ product that contains the top monomials of $\mu_1$ and $\tau_2$ cancels with the product that contains the top
monomials of $\mu_2$ and $\mu_3$. Also the codegree $2d$ products that contain a codegree $d$ monomial of $u_1$ which is from $b_1$ (i.e., also a monomial of $v_1$),
and the top monomial of $\tau_2$, cancel with the products that contain the same codegree $d$ monomial from $v_1$, and the top monomial of $\mu_3$.
\endroster

Because of (v), the analysis of codegree $2d$ monomials of $u_1$ and $v_1$ is identical to the analysis of codegree $d+r$ monomials of these elements.
This concludes the construction of the element $s_1$, and adds codegree $2d$ monomials to $\mu_1$, and codegree $d$ monomials to $w_1$. The analysis of 
the elements $u_3,v_3,b_3,\mu_3$ and $w_2$ is identical.

We continue by analyzing the codegree $2d$ monomials in $u_2,v_2,\tau_2$ and $\mu_2$.  The observations (i) and (ii) that we used in analyzing the codegree $d+r$
monomials of these elements for $r<d$ remain valid for codegree $2d$ monomials. In addition by part (v) in the analysis of codegree $2d$ monomials of $u_1$ and $v_1$,
it follows that
the codegree $2d$ product that contains the top monomials of $\mu_1$ and $\tau_2$ cancels with the product that contains the top
monomials of $\mu_2$ and $\mu_3$. Hence, the rest of the analysis of codegree $2d$ monomials of $u_2$ and $v_2$ is identical to the analysis of
codegree $d+r$ monomials of these elements for $r<d$.
 
\smallskip
We continue by analyzing higher codegree products and monomials. We assume inductively for $r>0$  that:
\roster
\item"{(1)}" $\mu_1 (s_1+1) =u_1$ and $\mu_1 s_1=v_1$ in the abelian group:
$G^{deg(u_1)}/G^{deg(u_1)-(d+r)}$. $deg(s_1)=d$. $u_1=v_1+\mu_1$ in the abelian group:
$G^{deg(u_1)}/G^{deg(u_1)-(2d+r)}$. 

\item"{(2)}" $(t_1+1) \mu_2 =v_2$, $\tau_2 (s_2+1)=u_2$ 
and $t_1\mu_2=\tau_2s_2$,
in the abelian group:
$G^{deg(u_2)}/G^{deg(u_2)-(d+r)}$. $deg(s_2)=deg(t_1)=d$. $deg(\mu_2)=deg(\tau_2)=deg(u_2)-d$.

\item"{(3)}" $b_2+\tau_2=u_2$ and $b_2+\mu_2=v_2$ in the abelian group:
$G^{deg(u_2)}/G^{deg(u_2)-(2d+r)}$. 

\item"{(4)}" $x_1=s_1w_1=w_1t_1=s_2w_2=w_2t_2$ in the abelian group:
$G^{deg(x_1)}/G^{deg(x_1)-(d+r)}$. 
\endroster

And we continue by analyzing codegree $2d+r$ products. The analysis is similar to the analysis of codegree $d+r$ and codegree $2d$ products.
We use the following observations:
\roster
\item"{(i)}" all the codegree $2d+r$ products that contain monomials of $u_i,v_i$, $i=1,2,3$, that are all of codegree smaller than $d$,
cancel in pairs. In particular, all the codegree $2d+r$ products that contain a monomial of $x$ of codegree bigger than $d+r$, in one of its two appearances,
cancel in pairs.

\item"{(ii)}" all the codegree $2d+r$ products that contain  monomials  from all $b_1$, $b_2$ and $b_3$, cancel in pairs.

\item"{(iii)}" a codegree $2d+r$ product that contains a monomial from $\mu_1$ of codegree more than $d$, and a monomial from the first appearance of $x$,
such that the sum of their codegrees is less than $2d+r$, an element from $b_2$ and an element from $b_3$, cancels with an a product that contains an element
from $b_1$, an element from the first appearance of $x$, an element from $\mu_2$, and the same element from $b_3$. The same holds for products that
contain monomials from $b_1$, $b_2$ and $\mu_3$ with parallel restrictions.

\item"{(iv)}" a codegree $2d+r$ product that contains a monomial from $b_1$, a monomial from $\tau_2$ of codegree bigger than $d$,  and a monomial from $b_3$, 
such that the sum of the
codegrees of the monomial from $\tau_2$ and the monomial from the second appearance of $x$ is smaller than $2d+r$, cancels with a product that contains the same
monomials of $b_1$ and the first appearance of $x$, a monomial from $b_2$ and a  monomial from $\mu_3$. 
The same holds for products that
contains monomials from $b_1$, $\mu_2$ and $b_3$ with parallel restrictions.

\item"{(v)}" a codegree $2d+r$ product that contains a monomial from $\mu_1$ of codegree bigger than $d$, a monomial from $\tau_2$,  and a monomial from $b_3$, 
cancels with a product that contains a monomial from $b_1$, a monomial from $\mu_2$, and a monomial from $\mu_3$.
The same holds for products that
contain monomials from $b_1$, $\mu_2$ and $\mu_3$ with parallel restrictions.
\endroster

Hence, like in the analysis of codegree $2d$ products, 
to analyze the elements: $u_1,v_1,b_1,\mu_1$ and $w_1$, we need to analyze codegree $2d+r$ products that
contain one of the following:  
\roster
\item"{(i)}"  a codegree $2d+r$ monomial of $u_1$  that does not appear in $v_1$ and vice versa. 

\item"{(ii)}" a codegree $d+r$ monomial of $v_1$ (which is a monomial of $b_1$) and the top monomial of $\mu_2$. 

\item"{(iii)}" a codegree $d+q_j$ monomial of $\mu_1$, $q_j <d+r$, and a codegree $d+r-q_j$ monomial of the first appearance of $x_1$.

\item"{(iv)}" a codegree $p_j$ monomial of $v_1$ (which is a monomial of $b_1$), $p_j < d+r$, and a codegree $d+r-p_j$ monomial of the first appearance of $x_1$ and
the top monomial of $\mu_2$.
\endroster

Hence, the analysis of codegree $2d+r$ monomials of $u_1$ and $v_1$ is identical to the analysis of codegree $d+r$ and $2d$  monomials of these elements. Note
that in analyzing products of codegree greater than $2d+r$, the element $s_1$ is already fixed, and we only add codegree $2d+r$ monomials to $\mu_1$ and $b_1$, and
codegree $d+r$ monomials to $w_1$.
The analysis of 
the elements $u_3,v_3,b_3,\mu_3$ and $w_2$ is identical.

We continue by analyzing the codegree $2d$ monomials in $u_2,v_2,\tau_2$ and $\mu_2$.  The observations (i)-(v) that we used in analyzing the codegree
$2d+r$ monomials of $b_1$ and $\mu_1$, imply that  analyzing codegree $2d+r$ monomials of $b_2, \tau_2$ and $\mu_2$, is similar to the analysis of the codegree $d+r$
monomials of these elements. Hence, we can finally deduce that:
\roster
\item"{(1)}" $\mu_1 (s_1+1) =u_1$ and $\mu_1 s_1=v_1$. 
$deg(s_1)=d$ and $u_1=v_1+\mu_1$. 

\item"{(2)}" $(t_1+1) \mu_2 =v_2$, $\tau_2 (s_2+1)=u_2$ 
and $t_1\mu_2=\tau_2s_2$.
$deg(s_2)=deg(t_1)=d$ and $deg(\mu_2)=deg(\tau_2)=deg(u_2)-d$.

\item"{(3)}" $b_2+\tau_2=u_2$ and $b_2+\mu_2=v_2$.

\item"{(4)}" $x_1=s_1w_1=w_1t_1=s_2w_2=w_2t_2$ in the abelian group:
$G^{deg(x_1)}/G^{deg(x_1)-deg(u_1u_2u_3)}$. 
\endroster

This proves the structure of the coefficients in the statement of proposition 4.11. Suppose that there exists a solution $x_2$ to the given equation, and:
$$deg(x_2)>2(2+2^{deg(s)+2}+
deg(u_1)+deg(u_2)+deg(u_3)).$$ 
As in the analysis of the same equation in case there are shifts between the appearances of the element $x_2$,
we can continue the analysis of higher codegree
 monomials of the solution $x_2$, and get that there exist elements $w_i$, $i=1,2$,
 that satisfy: $s_iw_i=w_it_i=x_1$ in the abelian group:
$G^{deg(x_2)}/G^{deg(s_1)-1}$. By the argument that was used to prove lemma 4.2, it follows that there exists a solution $\hat x$ to the equations: $s_ix=xt_i$,
$i=1,2$.

$x_2$ satisfies $s_1x_2=x_2t_1$ in the abelian group:  
$G^{deg(s_1x_2)}/G^{2deg(s_1)-1}$. Hence, there exists an element $\hat x_2$, which is a solution of the equation $s_1x=xt_1$, and $x_2+\hat x_2=r$, where
$deg(r) \leq   
2+2^{deg(s_1)+2}$.

$x_2$ is a solution to the equation: $v_1xv_2xv_3=u_1xu_2xu_3$, where: $v_1=\tau_1 s_1$, $u_1=\tau_1(s_1+1)$, $v_2=(t_1+1)\mu_2$, $u_2=\tau_2(s_2+1)$, $v_3=(t_2+1)\mu_3$,
$u_3=t_2 \mu_3$, and $\tau_2s_2=t_1 \mu_2$.
Hence: 
$$\tau_1(s_1+1)(\hat x_2+r)\tau_2(s_2+1)(\hat x_2+r)t_2\mu_3=\tau_1s_1(\hat x_2+r)(t_1+1)\mu_2(\hat x_2+r)(t_2+1)\mu_3$$
Therefore:
$$(s_1+1)r\tau_2(s_2+1)\hat x_2t_2+
(s_1+1)\hat x_2\tau_2(s_2+1)rt_2=$$
$$=s_1r(t_1+1)\mu_2\hat x_2(t_2+1)+
s_1\hat x_2(t_1+1)\mu_2r (t_2+1)$$
in the abelian group: 
$G^{deg(
s_1r\tau_2s_2\hat x_2t_2}/
G^{deg(
s_1r\tau_2s_2rt_2}$.

Since $s_1 \hat x_2=\hat x_2 t_1$, this implies:
$$((s_1+1)r\tau_2(s_2+1)s_2+
s_1r(t_1+1)\mu_2(s_2+1))\hat x_2=
\hat x_2((t_1+1)\tau_2(s_2+1)rt_2+
t_1(t_1+1)\mu_2r (t_2+1))$$
in the same abelian group.

Therefore, 
$$(s_1+1)r\tau_2(s_2+1)s_2+
s_1r(t_1+1)\mu_2(s_2+1)=p(s_1)$$  and:
$(t_1+1)\tau_2(s_2+1)rt_2+
t_1(t_1+1)\mu_2r (t_2+1)=p(t_1)$ for some polynomial $p$.

This implies that:
$r \tau_2 s_2+s_1r\mu_2$ is a polynomial in $s_1$, and: 
$\tau_2rt_2+t_1\mu_2r$ is a polynomial in $t_1$. Hence, $(rt_1+s_1r) \mu_2$ is a polynomial in $s_1$, and:
$\tau_2(rt_2+s_2r)$ is a polynomial in $t_1$.

Since we assumes that the top monomials of the coefficient do not contain periodicity, it can not be that the top monomials of $\tau_2$ and $\mu_2$ are
equal, and equal to the top monomials of $t_1$ and $s_1$, hence,  
$rt_1+s_1r=rt_2+s_2r \neq 1$.

If $deg(\tau_2)=deg(s_1)$, then the top monomials of $s_1$ and $t_1$ are equal, and the top monomials of $u_2$ and $v_2$ have periodicity, a contradiction.
The top monomial of $\tau_2$ has no periodicity, so $deg(\tau_2)<2deg(s_1)$.  
If $deg(\tau_2)> deg(s_1)$, then necessarily the top monomials of $u_2$ and $v_2$ contain periodicity, a contradiction. 

Suppose that $deg(\tau_2)<deg(s_1)$. If the top monomial of $\tau_2$ is the same as the top monomial of $\mu_2$, then the top monomials of the two sides of the
equation contain periodicity, a contradiction. If the top monomials of $\mu_2$ and $\tau_2$ are distinct, then the top monomials of $u_2$ and $v_2$ contain
periodicity, a contradiction. 

Therefore, $rt_1+s_1r=0$, so $r$ is a solution of the equation: $s_1x=xt_1$ and so is $x_2=\hat x_2+r$, and the conclusion of proposition 4.11 follows.

\line{\hss$\qed$}

Proposition 4.11 and its proof enable us to prove theorem 4.7 in case there are no shifts, i.e., in case the degrees of the
elements $u_i,v_i$ satisfy: $deg(u_i)=deg(v_i)$ for all indices $i$.

\vglue 1.5pc
\proclaim{Proposition 4.12} 
Let $u_1,\ldots,u_n,v_1,\ldots,v_n \in  FA$ and suppose that the equation:
$$u_1xu_2xu_3 \ldots u_{n-1}xu_n \, = \, v_1xv_2xv_3 \ldots v_{n-1}xv_n$$ has a solution $x_1$ of degree bigger than 
$2(deg(u_1)+ \ldots + deg(u_n))^2$. Suppose further  that:
\roster
\item"{(1}" for every index $i$, $1 \leq i \leq n$, $deg(u_i)=deg(v_i)$.

\item"{(2)}"  the top homogeneous parts of $u_i$ and $v_i$ are monomials
with no periodicity.

\item"{(3)}" for some index $i$, $u_i \neq v_i$.

\item"{(4)}"  all the periodicity
in the top monomials, that are associated with  the top monomials of the  two sides of the equation after substituting the solution $x_1$, 
is contained in the periodicity of the top
monomial of the solution $x_1$.
\endroster

Then there exist some elements $s,t \in FA$, $deg(s)=deg(t) < min deg(u_i)$, such that:
\roster
\item"{(1)}" every solution of the equation $sx=xt$ is a solution of the given equation.

\item"{(2)}" every  solution  $x_2$ of the given equation that satisfies:
$deg(x_2)>2(2+2^{deg(s)+2}+deg(u_1)+ \ldots +deg(u_n))$, is also a solution of the equation $sx=xt$.

\item"{(3)}" for every index $i$, $1 \leq i \leq r$, for which $u_i \neq v_i$, there exist elements $\tau_i, \mu_i$, such that the elements
$u_i,v_i$ are either: $\tau_i (s_i+1)$ or $(t_{i-1}+1) \mu_i$ or $\tau_i s_i$ or $t_{i-1} \mu_i$, where the elements $s_i$ are either $s$ or $s+1$,
and the elements $t_i$ are either $t$ or $t+1$, and $t_{i-1} \mu_i= \tau_i s_i$.
\endroster
\endproclaim

\nfp The proof of the structure of the coefficients is similar to the proof of proposition 4.11. Given the structure of the coefficients, it is clear that 
every solution of the equation $sx=xt$ is a solution of the given equation. It is left to prove that every long enough solution of the given
equation is a solution of the equation $sx=xt$.

Suppose that $x_2$ is a solution of the given equation that satisfies: 
$deg(x_2)>2(2+2^{deg(s)+2}+deg(u_1)+ \ldots +deg(u_n))$. By the argument that we used in proposition 4.11, it follows that the equation $sx=xt$ has a solution,
and that $x_2= \hat x_2+r$, where $\hat x_2$ is a solution to the equation $sx=xt$, and  $deg(r) \leq  
2+2^{deg(s)+2}$.

In that case we get the equality:
$$\tau_1(s_1+1)(\hat x_2+r)\tau_2(s_2+1)(\hat x_2+r) \ldots \tau_{n-1}(s_{n-1}+1)(\hat x_2+r)t_{n-1}\mu_n=
tau_1s_1(\hat x_2+r)(t_1+1)\mu_2(\hat x_2+r) \ldots (t_{n-2}+1)\mu_{n-1}(\hat x_2+r)(t_{n-1}+1)\mu_n$$
and since $\hat x_2$ is a solution of the equation: $s_1x=xt_1$, we get the equality:
$$(s_1+1)r\tau_2(s_2+1)\hat x_2 \ldots \tau_{n-1}(s_{n-1}+1)\hat x_2t_{n-1}+ \ldots +
(s_1+1)\hat x_2\tau_2(s_2+1)\hat x_2 \ldots \tau_{n-1}(s_{n-1}+1)rt_{n-1}=$$
$$s_1r(t_1+1)\mu_2\hat x_2 \ldots (t_{n-2}+1)\mu_{n-1}\hat x_2(t_{n-1}+1)+ \ldots +
s_1\hat x_2(t_1+1)\mu_2\hat x_2 \ldots (t_{n-2}+1)\mu_{n-1}r(t_{n-1}+1)$$
in the abelian group:
$G^{m_1}/G^{m_2}$, where: 
$$m_1=
deg((s_1+1)r\tau_2(s_2+1)\hat x_2 \ldots \tau_{n-1}(s_{n-1}+1)\hat x_2t_{n-1})$$
and $m_2=m_1-deg(\hat x_2)+deg(r)$.

This implies the equality:
$$\hat x_2((t_1+1)
\tau_2(s_2+1)\hat x_2 \ldots \tau_{n-1}(s_{n-1}+1)rt_{n-1}+
t_1(t_1+1)\mu_2\hat x_2 \ldots (t_{n-2}+1)\mu_{n-1}r(t_{n-1}+1))+$$
$$((s_1+1)r\tau_2(s_2+1)\hat x_2 \ldots \tau_{n-1}(s_{n-1}+1)s_{n-1}+
s_1r(t_1+1)\mu_2\hat x_2 \ldots (t_{n-2}+1)\mu_{n-1}(s_{n-1}+1))\hat x_2+
\hat x_2w\hat x_2$$
for some element  $w \in FA$, and in the same abelian group: $G^{m_1}/G^{m_2}$. Therefore:
$$(t_1+1)
\tau_2(s_2+1)\hat x_2 \ldots \tau_{n-1}(s_{n-1}+1)rt_{n-1}+
t_1(t_1+1)\mu_2\hat x_2 \ldots (t_{n-2}+1)\mu_{n-1}r(t_{n-1}+1)$$
is a polynomial in $t_1$ and:
$$(s_1+1)r\tau_2(s_2+1)\hat x_2 \ldots \tau_{n-1}(s_{n-1}+1)s_{n-1}+
s_1r(t_1+1)\mu_2\hat x_2 \ldots (t_{n-2}+1)\mu_{n-1}(s_{n-1}+1)$$
is a polynomial in $s_1$ in the same abelian group:
$G^{m_1}/G^{m_2}$. 
 
We examine the polynomial in $t_1$:
$$
\tau_2(s_2+1)\hat x_2 \ldots \tau_{n-1}(s_{n-1}+1)rt_{n-1}+
t_1\mu_2\hat x_2 \ldots (t_{n-2}+1)\mu_{n-1}r(t_{n-1}+1)$$

Note that:
$$
\tau_2(s_2+1)\hat x_2 \ldots \tau_{n-1}s_{n-1}=
t_1\mu_2\hat x_2 \ldots (t_{n-2}+1)\mu_{n-1}$$
Hence:
$$
\tau_2(s_2+1)\hat x_2 \ldots \tau_{n-1}(s_{n-1}+1)rt_{n-1}+
t_1\mu_2\hat x_2 \ldots (t_{n-2}+1)\mu_{n-1}r(t_{n-1}+1)=$$
$$\tau_2(s_2+1)\hat x_2 \ldots \tau_{n-1}rt_{n-1}+
t_1\mu_2\hat x_2 \ldots (t_{n-2}+1)\mu_{n-1}r$$

Since:
$$\tau_2(s_2+1)\hat x_2 \ldots \tau_{n-1}s_{n-1}=
t_1\mu_2\hat x_2 \ldots (t_{n-2}+1)\mu_{n-1}$$
It follows that:
$$\tau_2(s_2+1)\hat x_2 \ldots \tau_{n-1}(rt_{n-1}+s_{n-1}r)$$
is a polynomial in $t_1$ in the abelian group:
$G^{m_1}/G^{m_2}$. 
 
If $rt_{n-1}+s_{n-1}r \neq 0$, then as we argued in the proof of proposition 4.11, and in case we assumed that there is a shift between
the appearances of the elements $x_2$ from the two sides of the given equation, either the top monomials of the coefficients $u_i$ and $v_i$ are periodic, or
the top monomials in the two sides of the given equation has periodicity that is not contained in the appearances of the top monomial of the solution
$x_2$. Both contradict our assumptions and the proposition follows.

\line{\hss$\qed$}

At this point we need to consider equations in which some of the appearances of the elements $x$ are shifted, and some are not.

\vglue 1.5pc
\proclaim{Lemma 4.13} Let $u_1,u_2,u_3,v_1,v_2,v_3 \in FA$ satisfy $u_1 \neq v_1$, $deg(u_1)=deg(v_1)$, $deg(u_2)>rk(v_2)$, $rk(v_3)>rk(u_3)$.
Suppose that the top homogeneous parts of $u_i$ and $v_i$ 
are monomials
(for $i=1,2,3)$) with no non-trivial periodicity. If there exists a solution $x_1$ to the equation $u_1xu_2xu_3=v_1xv_2xv_3$,
and the only non-trivial periodicity in the top monomials of the two sides of the equation is contained in the top monomial of the solution $x_1$,  
and $deg(x_1)> 2(deg(u_1)+deg(u_2)+deg(u_3)^2)$, then there
exist elements $s,t \in FA$, such that either:
\roster
\item"{(1)}" there exists $\mu_1$ for which $u_1=\mu_1 (s+1)$ and $v_1=\mu_1 s$.

\item"{(2)}" there exist $\mu_2$ and $s_2,t_2$, for which: $(t+1) \mu_2=v_2$ and $t \mu_2s_2=u_2$. Furthermore, $v_3=t_2u_3$ and 
the pair $(s_2,t_2)$ is either $(s,t)$ or $(s+1,t+1)$.
\endroster

or:
\roster
\item"{(1)}" there exists $\mu_1$ for which $u_1=\mu_1 s$ and $v_1=\mu_1 (s+1)$.

\item"{(2)}" there exist $\mu_2$ and $s_2,t_2$, for which: $(t+1) \mu_2=u_2$ and $v_2s_2=t \mu_2$. Furthermore, $v_3=t_2u_3$ and 
the pair $(s_2,t_2)$ is either $(s,t)$ or $(s+1,t+1)$.
\endroster

As in the conclusion of theorem 4.7, every solution of the equation $sx=xt$ is a solution of the given equation $u_1xu_2xu_3=v_1xv_2xv_3$. 
Every solution $x_2$ of the given equation:
$u_1xu_2xu_3=v_1xv_2xv_3$, that satisfies:
$deg(x_2)>2(2+2^{deg(s_1)+2}+deg(u_1)+deg(u_2)+deg(u_3))$, is also a solution of the equation: $sx=xt$.
\endproclaim

\nfp The proof is similar to the proof of proposition 4.11.

\line{\hss$\qed$}

At this point we can complete the proof of theorem 4.7. We already analyzed the case 
in which there are non-trivial shifts between (the top monomials of) pairs of appearances of the variable $x$ in the two sides of the equation.
Propositions 4.11 and 4.12 analyze the case in which there are no shifts between pairs of appearances of the variable $x$ in the two sides of
the equation, and proposition 4.13 analyzes the case $n=3$ in which there is a pair with no shift and a pair with a shift.

By the techniques that were used in proving proposition 4.11 and in analyzing the case in which there are non-trivial shifts between pairs of
appearances of the variable $x$, if there is a pair of coefficients, $u_i,v_i$, such that $u_i=v_i$ and the $i-1$ (hence, also the $i$-th)
pair of appearances of the variable $x$ has no shift, then the equation breaks into two equations, the first contains the coefficients: $u_1,\ldots,u_{i-1},v_1,\ldots,v_{n-1}$,
and the second contains the coefficients: $u_{i+1},\ldots,u_n,v_{i+1},\ldots,v_n$. Therefore, in the sequel we may assume that there is no such pair of
coefficients: $u_i,v_i$.

Then there exist some elements $s,t \in FA$, $deg(s)=deg(t) < min deg(u_i)$, and elements: $s_1,\ldots,s_{n-1},t_1,\ldots,t_{n-1}$, such that:
\roster
\item"{(1)}" for every index $i$, the pair $(s_i,t_i)$ is either: $(s,t)$ or $(s+1,t+1)$.

\item"{(2)}" for every pair of coefficients, $u_i,v_i$, for which the two pairs of appearances of the variable $x$ from the two sides of the pair
of coefficients
have no non-trivial shift, either $u_i=v_i$, or there exist elements $\tau_i$ and $\mu_i$ such that either: $u_i= \tau_i s_i$ and $v_i=\tau_i(s_i+1)$ (or vice versa), or
$u_i= t_{i-1} \mu_i$ and $v_i=(t_{i-1}+1) \mu_i$ (or vice versa), or
$u_i= (t_{i-1}+1) \mu_i$ and $v_i= \tau_i (s_i+1)$ (or vice versa). 

\item"{(3)}" if $deg(u_1)=deg(v_1)$, either $u_1=v_1$, or there exists $\tau_1$ such that:
$u_1= \tau_1 s_1$ and $v_1=\tau_1(s_1+1)$ (or vice versa). 
If $deg(u_n)=deg(v_n)$, either $u_n=v_n$, or there exists an element $\mu_n$ such that:   
$u_n= t_{n-1} \mu_n$ and $v_n=(t_{n-1}+1) \mu_n$ (or vice versa). 

\item"{(4)}" for every pair $u_i,v_i$, for which the two pairs of appearances of the variable $x$ from the two sides of the pair of coefficients have
non-trivial shifts, $u_is_i=t_{i-1}v_i$ (or vice versa), or $u_i=t_{i-1}v_is_i$ (or vice versa).

\item"{(5)}" if $deg(u_1) \neq deg(v_1)$,  
then: $u_1= v_1 s_1$ or vice versa.
If $deg(u_n) \neq deg(v_n)$, then: 
$u_n= t_{n-1} v_n$ or vice versa. 

\item"{(6)}" suppose that  $deg(u_i) \neq deg(v_i)$, $1 < i < n$, there is no shift between the $i-1$ appearances of the variable $x$, and there is a 
non-trivial shift between the $i$th appearances of the variable $x$ from the two sides of the equation. The either: $u_is_i=v_i$ or vice versa, in which case the original
equation can be broken into two equations, the first contains the first $i-1$ pairs of coefficients, and the second contains the last $n+1-i$ pairs of coefficients.
Or $v_i=(t_{i-1}+1) \mu_i$ and $u_is_i=t_{i-1} \mu_i$, or vice versa. Or $u_i=(t_{i-1}+1) \mu_i$ and $v_i=t_{i-1} \mu_i s_i$, or vice versa. 

\item"{(7)}" suppose that  $deg(u_i) \neq deg(v_i)$, $1 < i < n$, there is no shift between the $i$th appearances of the variable $x$, and there is a 
non-trivial shift between the $i-1$ appearances of the variable $x$ from the two sides of the equation. Then either: $t_{i-1} u_i=v_i$ or vice versa, in which case the original
equation can be broken into two equations, the first contains the first $i$ pairs of coefficients, and the second contains the last $n-i$ pairs of coefficients.
Or $v_i=\tau_i(s_{i}+1)$ and $t_{i-1}u_i=\tau_i s_i$, or vice versa. Or $u_i=\tau_i (s_i+1)$ and $v_i=t_{i-1} \tau_i s_i$, or vice versa. 
\endroster

This description of the coefficients in a general equation with one variable, in which the coefficients have no periodicity, and the top homogeneous parts of the
coefficients are monomials, finally implies:
\roster
\item"{(1)}" every solution of the equation $sx=xt$ is a solution of the given equation.

\item"{(2)}" every solution  $x_2$ of the given equation,  that satisfies: 
$deg(x_2)>2(2+2^{deg(s)+2}+deg(u_1)+ \ldots +deg(u_n)$, is also a solution of the equation: $sx=xt$.
\endroster

The proof of (1) follows from the structure of the coefficients, and the proof of (2) follows by the argument that was used to prove (2) for the case in which there
are no shifts between the various appearances of the top monomial of the solution $x_2$ in the two sides of the given equation in proposition 4.12.

This concludes the proof of theorem 4.7.

\line{\hss$\qed$}


\vskip 2em

\end